\documentclass{crmp-l}
\usepackage{amssymb}
\usepackage{amsbsy}
\usepackage{amscd}
\usepackage[mathscr]{eucal}
\usepackage{verbatim}

\makeatletter
\@ifclasslater{amsproc}{1999/11/24}{}{
\renewcommand*\subjclass[2][1991]{%
  \def\@subjclass{#2}%
  \@ifundefined{subjclassname@#1}{%
    \ClassWarning{\@classname}{Unknown edition (#1) of Mathematics
      Subject Classification; using '1991'.}%
  }{%
    \@xp\let\@xp\subjclassname\csname subjclassname@#1\endcsname
  }%
}
\renewcommand{\subjclassname}{%
  \textup{1991} Mathematics Subject Classification}
\@xp\let\csname subjclassname@1991\endcsname \subjclassname
\@namedef{subjclassname@2000}{%
  \textup{2000} Mathematics Subject Classification}
}
\makeatother

\usepackage{pstricks}
\usepackage{curves}

\makeatletter
\def\cal{\mathcal}
\def\Bbb{\mathbb}

\newenvironment{pf*}[1]{\proof[#1]}{\endproof}
\newcommand{\rom}{\textup}
\makeatother
\newenvironment{aenume}{%
  \begin{enumerate}%
  }{\end{enumerate}}
\newenvironment{NB}{
{\bf NB}. \footnotesize
}{}
 \renewenvironment{NB}{%
   \comment
 }{\endcomment}
\newtheorem{Theorem}[equation]{Theorem}
\newtheorem{Corollary}[equation]{Corollary}
\newtheorem{Lemma}[equation]{Lemma}
\newtheorem{Proposition}[equation]{Proposition}

\theoremstyle{definition}

\newtheorem{Notation}[equation]{Notation}

\theoremstyle{remark}
\newtheorem{Remark}[equation]{Remark}

\numberwithin{equation}{section}

\newcommand{\thmref}[1]{Theorem~\ref{#1}}
\newcommand{\secref}[1]{\S\ref{#1}}

\newcommand{\propref}[1]{Proposition~\ref{#1}}
\newcommand{\corref}[1]{Corollary~\ref{#1}}
\newcommand{\subsecref}[1]{\S\ref{#1}}

\newcommand{\defeq}{\overset{\operatorname{\scriptstyle def.}}{=}}
\newcommand{\C}{{\Bbb C}}
\newcommand{\Z}{{\Bbb Z}}
\newcommand{\Q}{{\Bbb Q}}
\newcommand{\R}{{\Bbb R}}
\newcommand{\proj}{{\Bbb P}}
\newcommand{\CP}{\proj}

\newcommand{\SU}{\operatorname{\rm SU}}
\newcommand{\GL}{\operatorname{GL}}

\newcommand{\U}{\operatorname{\rm U}}
\newcommand{\SO}{\operatorname{\rm SO}}

\newcommand{\End}{\operatorname{End}}
\newcommand{\Hom}{\operatorname{Hom}}

\newcommand{\rank}{\operatorname{rank}}

\newcommand{\tr}{\operatorname{tr}}

\newcommand{\pd}[2]{\frac{\partial#1}{\partial#2}}

\newcommand{\id}{\operatorname{id}}
\newcommand{\ve}{\varepsilon}

\newcommand{\linf}{\ell_\infty}
\newcommand{\shfO}{\mathcal O}
\newcommand{\dslash}{/\!\!/} 
\newcommand{\bp}{{\widehat\proj}^2}
\newcommand{\bM}{{\widehat M}}
\newcommand{\Mreg}{M_0^{\operatorname{reg}}}

\newcommand{\Supp}{\operatorname{Supp}}
\newcommand{\ch}{\operatorname{ch}}

\newcommand{\Fin}{F^{\text{\rm inst}}}

\newcommand{\Fz}{F_0}
\newcommand{\Fo}{F_1}
\newcommand{\q}{\mathfrak q}
\newcommand{\hT}{\widetilde T}

\newcommand{\Zin}{Z^{\text{\rm inst}}}
\newcommand{\bZin}{\widehat{Z}^{\text{\rm inst}}}

\newcommand{\bE}{\widehat{\mathcal E}}

\begin{document}
\title{Lectures on Instanton Counting}
\author{Hiraku Nakajima}
\address{Department of Mathematics, Kyoto University, Kyoto 606-8502,
Japan}
\email{nakajima@math.kyoto-u.ac.jp}
\thanks{The first author is supported by the Grant-in-aid
for Scientific Research (No.15540023), JSPS}

\author{K\={o}ta Yoshioka}
\address{Department of Mathematics, Faculty of Science, Kobe University,
Kobe 657-8501, Japan}
\email{yoshioka@math.kobe-u.ac.jp}

\dedicatory{Dedicated to Professor Akihiro Tsuchiya on his sixtieth
birthday}

\begin{abstract}
These notes have two parts. The first is a study of Nekrasov's
deformed partition functions $Z(\ve_1,\ve_2,\vec{a};\q,\vec{\tau})$ of
$N=2$ SUSY Yang-Mills theories, which are generating functions of the
integration in the equivariant cohomology over the moduli spaces of
instantons on $\mathbb R^4$. The second is review of geometry of the
Seiberg-Witten curves and the geometric engineering of the gauge
theory, which are physical backgrounds of Nekrasov's partition
functions.

The first part is continuation of our previous paper \cite{part1},
where we identified the Seiberg-Witten prepotential with
$Z(0,0,\vec{a};\q,0)$.
We put higher Casimir operators to the partition function and clarify
their relation to the Seiberg-Witten $u$-plane. We also determine the
coefficients of $\ve_1\ve_2$ and $(\ve_1^2+\ve_2^2)/3$ (the genus $1$
part) of the partition function, which coincide with two measure
factors $A$, $B$ appeared in the $u$-plane integral.
The proof is based on the blowup equation which we derived in
\cite{part1}.

These notes are based on authors' lectures at Centre de Recherches
Math\'ema\-tiques, Universit\'e de Montr\'eal, July 2003.
\end{abstract}

\maketitle
\tableofcontents


\subsection*{Acknowledgement}
The authors are grateful to N.~Nekrasov for discussions, especially on
perturbation terms and the genus $1$ correction. The authors would
like to their sincere gratitude to Professor Akihiro Tsuchiya on his
encouragement for us and organization of many workshops over
years. This work cannot be done without his enthusiasm for
understanding physics.

\section{Introduction}

In this long introduction, we review a history of Donaldson invariants
and Seiberg-Witten geometry, which leads to the Nekrasov's deformed
partition function. This section contains no mathematically rigorous
results, but provides the motivation for our study in later sections.

\subsection{Donaldson invariants : a mathematical definition}
Let $X$ be a smooth, compact, oriented, $4$-manifold with a Riemannian
metric $g$ with $b^+\ge 1$ and odd. We also assume $\pi_1(X) = 1$ for
brevity.
Let $P\to X$ be an $\SO(3)$-bundle over $X$. Let $M(P)$ be the moduli
space of irreducible anti-self-dual connections on $P$. This is a
manifold with dimension $-2 p_1(P) - 3(1+b^+)$ for a generic metric
$g$. Let $\mathcal P\to X\times M(P)$ be the universal bundle.  Then
the Donaldson invariant is a polynomial on $H_0(X)\oplus H_2(X)$
defined by
\begin{equation*}
   D_P( p^a S^b) = \int_{M(P)} \mu(p)^a \mu(S)^b, \qquad
   p\in H_0(X), S\in H_2(X),
\end{equation*}
where $\mu\colon H_i(X)\to H^{4-i}(M(P))$ is given by the slant
product $\mu(\bullet) = -\frac14 p_1(\mathcal P)/\bullet$. Since
$M(P)$ is not compact, we must justify the definition of the
integration, and this can be done by using the Uhlenbeck
compactification, as one can find in textbooks on the Donaldson
theory \cite{DK,FM}.

We then formulate a generating function
\begin{equation*}
   D_\xi(p,S) = \sum_P \sum_{m,n\ge 0}
   D_P(\frac{S^m}{m!}\frac{p^n}{n!}),
\end{equation*}
where $\xi = w_2(P)$ is fixed.
When $b^+ > 1$, $D_\xi$ is independent of the metric and defines
invariants of the differentiable structure of $X$. When $b^+=1$, it is
piecewise constant as a function of $g$.

Although the invariants $D_\xi$ can be defined, their calculation was
difficult in general. This was because it is difficult to describe the
moduli spaces $M(P)$ explicitly. The situation was changed when
Kronheimer-Mrowka \cite{KM} proved a structure theorem for $D_\xi$ in
1994: Although $D_\xi$ involves infinitely many moduli spaces, it is
determined by finite data, if $D_\xi$ satisfies a so-called {\it
simple type\/} condition.

Soon afterward, Fintushel-Stern obtained the `blowup formula' which
describe the relation between $D_\xi$ on $X$ and that on the blowup
$\widehat X$ \cite{FS}. The formula involves an elliptic function. The
underlying elliptic curve is related to the structure theorem so that
the simple type condition means that it degenerates to a rational
curve. The blowup formula will play a fundamental role in this paper.

\subsection{Seiberg-Witten geometry}

In 1988, Witten described Donaldson invariants as correlation
functions of certain operators in a twisted version of $\mathcal N=2$
SUSY (supersymmetric) Yang-Mills theory \cite{Wit}. We do not explain
what this statement means here, but mention that it is an infinite
dimensional analogue of the Chern-Weil formula \cite{AL}.

Shortly after \cite{KM} was appeared, Seiberg-Witten analyzed the
original $\mathcal N=2$ SUSY Yang-Mills theory with gauge group
$\SU(2)$ \cite{SW}. The original theory is formulated on $\R^4$, and
was no mathematically rigorous definition of the `prepotential', which
they calculated, at that time.  Giving such a definition is one of the
main purpose of these notes.  (See \thmref{thm:main}.) But we present
an `informal' definition here.

Let $H^*_{\SU(2)}(\mathrm{pt})$ be the $\SU(2)$-equivariant cohomology
of a point with complex coefficients. It is naturally identified with
the Weyl group (in this case $\{ \pm 1\}$) invariant part of the
symmetric product of the dual of the (complexified) Cartan subalgebra
$\mathfrak h$ (in this case $\C$). It is the coordinate ring of
$\mathfrak h/W$. This space $\mathfrak h/W$ is the classical limit of
the so-called {\it $u$-plane}, a family of `vacuum states', which
plays the most important role in the Seiberg-Witten geometry.

The coordinate ring $A(\mathfrak h/W) = H^*_{\SU(2)}(\mathrm{pt})$ has
a generator
\(
   -\frac12 \tr \left(
     \begin{smallmatrix}
      -a & 0 \\ 0 & a
     \end{smallmatrix}
     \right)^2 = -a^2,
\)
where $a$ is considered as a coordinate on $\mathfrak h$.
Let us denote it by $u_{\mathrm{cl}}$ since it is a coordinate of the
classical limit of the $u$-plane. We make a `quantum correction' $u$
of the function $u_{\mathrm{cl}}$ by using the framed moduli space
$M(2,n)$ of instantons on $S^4$. The precise definition will be given
below, but it is roughly given by
\begin{equation}\label{eq:expect_u}
   u = - \left.\sum_{n\ge 0} \Lambda^{4n} \int_{M(2,n)} \mu(p)
   \right/
   \sum_{n\ge 0} \Lambda^{4n} \int_{M(2,n)} 1.
\end{equation}
Here $\Lambda$ is a formal variable, the integration is done in the
equivariant homology group, and $\mu$ is defined by the same formula
as in Donaldson invariants. The moduli space $M(2,n)$ has an
$\SU(2)$-action given by the change of the framing. The classical part
is the term $n=0$, then $M(2,0)$ is a single point, so the integration
is just an identity operator. In this case, $\mu(p)\in
H^*_{\SU(2)}(\mathrm{pt})$ is nothing but the generator
$-u_{\mathrm{cl}}$. Thus $u_{\mathrm{cl}}$ is the classical limit of
$u$ as we explained.

When $n > 0$, the moduli space $M(2,n)$ is {\it noncompact\/} and we
need to justify the integration. Here the problem is not a technical
one, and has a very different nature from the noncompactness appeared
in the definition of Donaldson invariants, which was overcome by
Uhlenbeck compactification. In fact, if $M(2,n)$ had a suitable
compactification, the integration of $1$ would be $0$ by the degree
reason. The integration will be defined via the localization theorem
in the equivariant homology group. The precise formulation will be
given in \subsecref{subsec:integ}. As the upshot, the integral does
not have the value $H^*_{\SU(2)}(\mathrm{pt})$, but in its {\it
fractional field}. (In fact, we need to consider extra two dimensional
torus as below. Or we should consider $u$ as an operator as in
\subsecref{subsec:reform}.)
Thus $u$ is a rational function on $\mathfrak h/W$. In the
Seiberg-Witten geometry, the role of $u$ and $u_{\mathrm{cl}}$ is
reversed. We define the {\it $u$-plane\/} as the parameter space for
$u$, i.e., $u$ is the coordinate of the $u$-plane. Then we consider
$u_{\mathrm{cl}}$ (and $a$) as a rational function on the $u$-plane.

Other than the function $u$, there are several important geometric
objects on the $u$-plane. They are defined via the integration over
the instanton moduli spaces. One of the most important objects is the
{\it prepotential\/}, which has a form:
\begin{equation}\label{eq:prepotential}
   \mathcal F_0 = \mathcal F_0^{\mathrm{pert}} + 
   \mathcal F_0^{\mathrm{inst}}
   ,
\end{equation}
where $\mathcal F_0^{\mathrm{pert}}$ is the {\it perturbative part\/}
of the prepotential, which is an explicit rational function on
$\mathfrak h/W$. The part $\mathcal F_0^{\mathrm{inst}}$ is the {\it
instanton part\/}, and is a power series in $\Lambda^4$. The
coefficient of $\Lambda^{4n}$ is given by integration over $M(2,n)$.
The $u$-plane is a {\it special K\"ahler\/} manifold, where the
prepotential is included in its definition. For example, the K\"ahler
metric is the imaginary part of the second derivative of the
prepotential. See \cite{Fr} for more detail.

The main result of \cite{SW} is the determination of the $u$-plane and
the prepotential $\mathcal F_0$. As a result, the $u$-plane is the
parameter space for elliptic curves:
\begin{equation*}
   y^2 = (z^2 + u - 2\Lambda^2)(z^2 + u + 2\Lambda^2).
\end{equation*}
The prepotential $\mathcal F_0$ is given by using certain elliptic
integrals.
The original method used for the determination was a highly nontrivial
physical argument. One of the most essential ingredients is
understanding of its behavior under the `duality' transformation
$\tau \mapsto -1/\tau$, where $\tau$ is the period of the above
elliptic curve, which is given by the second derivative of $\mathcal
F_0$ with respect to the coordinate $a$. This is rather mysterious
transformation in view of the definition \eqref{eq:prepotential}. In
our approach, we will see theta functions quite naturally. So the
duality will come from the Poisson summation formula, but we do not
really understand its geometric origin.

Note that this picture is very similar to that of the {\it mirror
symmetry}. The prepotential above is a counterpart of the
Gromov-Witten invariants and is on the `symplectic'
side. The elliptic curves (Seiberg-Witten curves) are on the `complex'
side. In fact, this is not just analogy. The geometric engineering
which will be reviewed in \subsecref{subsec:engineer} explains the
result as a special case of the mirror symmetry.

For a later purpose, we give some functions explicitly. Let $\tau$ be
the period of the Seiberg-Witten elliptic curve. Then
\begin{equation}\label{eq:rank2}
\begin{gathered}[c]
  u =
  - \frac{\theta_{00}^4+\theta_{10}^4}
  {\theta_{00}^2\theta_{10}^2}\Lambda^2,
\\
  \frac{du}{da} = \frac{2\sqrt{-1}}{\theta_{00}\theta_{10}}\Lambda,
\\
  a = \sqrt{-1}\frac
  {2E_2+\theta_{00}^4+\theta_{10}^4}{3\theta_{00}\theta_{10}}\Lambda.
\end{gathered}
\end{equation}
Here $\theta_* = \theta_*(0|\tau)$ is the theta function and $E_2 =
E_2(\tau)$ is the (normalized) second Eisenstein series. The reader
should be careful when he/she compares these with the formulas in
\cite{MW}. Our $u$ (resp.\ $a$) is multiplied by $-2$ (resp.\ 
$2\sqrt{-1}$).

Finally note that the elliptic curve becomes singular at $u =
\pm2\Lambda^2$. In the classical limit $\Lambda\to 0$, these fall into
a single point $0$, which is the singular point in the classical
$u$-plane $\mathfrak h/W$.

\subsection{The $u$-plane integral}\label{subsec:u-plane_integral}
We return back to Donaldson invariants. Witten \cite{Wi} explained
that $D_\xi$ has three contributions:
\begin{equation*}
   D_\xi(p,S) = Z_u(p,S) + Z_{+}(p,S) + Z_{-}(p,S).
\end{equation*}
The parts $Z_{\pm}(p,S)$ come from the measure supported on the
singularity $\pm2\Lambda^2$ of the $u$-plane. These are given by
invariants defined via the moduli spaces of monopoles, called
Seiberg-Witten invariants. As for application to topology, $Z_u$ is
irrelevant as it depends only on $H^2(X,\Z)$. Furthermore, $Z_u$
vanishes when $b_+ > 1$. But we are interested in structures of
instanton moduli spaces which are reflected in $Z_u$.

When $b_+=1$, more precise description of $D_\xi$ was given by
Moore-Witten \cite{MW}. (See also \cite{LNS1,LNS2} for similar
results.) We briefly recall their description, since some parts are
closely related to our study.
The parts $Z_{\pm}(p,S)$ are written by the Seiberg-Witten invariants
summed over various choices of $\mathrm{Spin}^c$ structures. See
\cite[\S7]{MW} for the explicit expression. The remaining part $Z_u$
is the integration with respect to a smooth volume form. It is called
the {\it u-plane integral}. 
We choose and fix a harmonic self-dual two form $\omega$ with $\int_X
\omega\wedge\omega = 1$. This is unique up to sign, and the choice of
$\omega$ is related to the orientation of the moduli space.
We also put $\Lambda = 1$. Then
\begin{equation}\label{eq:u-plane_int}
  Z_u(p,S) = \int_{\text{$u$-plane}} da d\overline{a}
  A(u)^{\chi} B(u)^{\sigma} e^{pu + S^2 T} \Psi,
\end{equation}
with
{\allowdisplaybreaks
\begin{gather*}
   A(u) = \alpha \left(\frac{du}{da}\right)^{1/2}, \qquad
   B(u) = \beta (u^2 - 4)^{1/8},
\\
  T = \frac1{24} \left(\frac{du}{da}\right)^2 E_2(\tau) - \frac16 u,
\\
\begin{aligned}[t]
   \Psi = 
   - &\frac{\sqrt{-2}}{4y^{1/2}}
   \frac{d\overline{\tau}}{d\overline{a}}
   \exp\left[
     \frac1{8\pi y}\left(\frac{du}{da}\right)^2 S_+^2
   \right]
   e^{2\pi\sqrt{-1}\lambda_0^2}
   \sum_{\lambda\in H^2 + \frac12 \xi}
   (-1)^{(\lambda - \lambda_0)\cdot w_2(X)}
\\     
   &\times\left[(\lambda,\omega)-\frac{1}{4\pi y}\frac{du}{da}
       (S,\omega)\right]
     \exp\left[
       -\sqrt{-1}\pi\overline{\tau}\lambda_+^2
       - \sqrt{-1}\pi\tau\lambda_-^2
       + \frac{du}{da}(S,\lambda_-)
     \right].
\end{aligned}
\end{gather*}
Here} $\chi$ (resp.\ $\sigma$) is the Euler number (resp.\ signature)
of $X$, $\alpha$, $\beta$ are universal constants independent of $X$,
$\tau = x + iy$, $\lambda_0$ is a fixed element in $\frac12 \xi +
H^2(X,\Z)$, and $(\bullet)_\pm$ denotes the self-dual and anti-self-dual
part of $\bullet$ respectively.

Since this is a divergent integral, and we must regularize it. See the
original paper \cite{MW} how it is done.

The term $T$ is called a {\it contact term}. Its determination can be
done by several ways. In \cite{LNS1} by equating the answers givin by
various ways, a nontrivial equation was derived. This is the {\it
contact term equation\/}, which will be important for our study. See
\thmref{thm:contact} and \thmref{thm:main}.
The terms $A$, $B$ come from a Riemannian metric $g$.

Let us analyze the effect of the blowup $\widehat X\to X$ on the
$u$-plane integral since it is closely related to our study.
Let $C$ be the exceptional curve. We want to evaluate $Z_u(p,S+tC)$,
where $S\in H_2(X)$ is considered as a class of $H_2(\widehat X)$ via
the projection.

Since $\chi(\widehat X) = \chi(X) + 1$, $\sigma(\widehat X) =
\sigma(X) - 1$, the factor $A(u)^\chi B(u)^\sigma$ is multiplied by
\begin{equation*}
   \frac{A(u)}{B(u)} = \frac\alpha\beta (u^2-4)^{-1/8}
   \left(\frac{du}{da}\right)^{1/2}
   = \frac1{\theta_{01}
   } \text{(up to constant)}.
\end{equation*}
We work in a chamber $C_+ = 0$, so we have
\begin{equation*}
\begin{split}
   & \frac{\Psi_{\widehat X}}{\Psi_X}
   = \sum_{n\in\Z+\frac12 w_2(\widetilde P)\cdot C} (-1)^n
   \exp\left[\sqrt{-1}\pi\tau n^2 - nt \frac{du}{da}\right]
\\
   = \; &
   \theta_{*}(\frac{t\sqrt{-1}}{2\pi} \frac{du}{da}|\tau),
\end{split}
\end{equation*}
where $* = 01$ or $11$ according to $w_2(\widetilde P)\cdot C = 0$ or
$1$. Therefore we get
\begin{equation*}
   \frac{Z_u(p,S+tC)}{Z_u(p,S)}
   = \exp(-T t^2)
   \frac{\theta_{*}(\frac{t\sqrt{-1}}{2\pi} \frac{du}{da}|\tau)}
   {\theta_{01}(0|\tau)}
\end{equation*}
up to a constant multiple. The constant turns out to be $1$ as left
hand side is $1$ at $t=0$ when $* = 01$.

\subsection{Nekrasov's deformed partition function}

As we explained, the prepotential $\mathcal F_0$ was given as
integration over instanton moduli spaces. Before Nekrasov gave an
explicit expression \cite{Nek}, it was written in terms of
differential forms on moduli spaces. So it was difficult to calculate,
understand its meaning... (See \cite{DHKM}.) Nekrasov's idea was to
use an extra $2$-dimensional torus action and apply the localization
theorem in the equivariant homology. Technically it was also important
that the Uhlenbeck (partial) compactification of the moduli space has
a nice resolution of singularities introduced by the first author
\cite{Na-R4}. (The latter space will be denoted by $M(2,n)$ in the
main body of the paper.) Let $\ve_1,\ve_2$ be two generators of
$H^*_{T^2}(\mathrm{pt})$. Then we define
\begin{equation}
\label{eq:prepotential2}
    F = \ve_1\ve_2 F^{\mathrm{pert}} + 
   \ve_1\ve_2 \log
   \left(\sum_{n\ge 0} \Lambda^{4n} \int_{M(2,n)} 1\right),
\end{equation}
where $F^{\mathrm{pert}}$ is a certain two parameter deformation of
$\mathcal F_0^{\mathrm{pert}}$. Each coefficient of $\Lambda^{4n}$ is
a rational function in $\ve_1,\ve_2$, and is a mathematically
rigorously defined. Nekrasov conjectured $F|_{\ve_1,\ve_2=0}$ is equal
to $\mathcal F_0$, given by the Seiberg-Witten curve. This is
mathematically meaning full statement. This conjecture was proved by
\cite{part1} and \cite{NO} by totally different methods.
%

The method used in \cite{NO} was geometric and a standard technique in
the study of Donaldson invariants. We consider the instanton moduli
spaces $\bM(2,c_1,n)$ on the blowup, introduce an operator $\mu(C)$ in
this equivariant setting, and compute this equivariant analog of
Donaldson invariants. From the explicit expression given by the
localization theorem, it is very easy to derive the blowup formula in
a combinatorial form. On the other hand, by a simple dimension
counting argument shows that $\int_{\bM(2,0,n)} \mu(C)^2 = 0$. This
vanishing give a differential equation satisfied by the original $F$.
We call it {\it the blowup equation}. (See
\eqref{eq:blowup_formula2}.) It characterizes $F$. When we put
$\ve_1=\ve_2 = 0$, this equation turns out to be the contact term
equation, which we mentioned. Since the contact term equation can be
derived from the Seiberg-Witten curve in a mathematically rigorous way
(see \secref{sec:SW}), this gives a proof of Nekrasov's conjecture.

\subsection{Gravitational corrections}
After identifying $F|_{\ve_1,\ve_2=0}$ with the Seiberg-Witten
prepotential \eqref{eq:SWprep}, it becomes natural to ask the meaning
of higher order terms in the expansion
\begin{equation*}
    F = F_0 + (\ve_1 + \ve_2) H + \ve_1\ve_2 A
    + \frac{\ve_1^2+\ve_2^2}3 B + \cdots .
\end{equation*}
Nekrasov asserted that these are {\it gravitational corrections\/} to
the gauge theory \cite[\S4]{Nek}. Using the differential equation
mentioned above, we prove these $A$, $B$ coincides with those $A$, $B$
appeared in the $u$-plane integrand. ($H$ turns out to be a simple
function.) (The calculation was done jointly with N.~Nekrasov.)

Moreover, by the geometric engineering \cite{KKV} (see
\subsecref{subsec:engineer}), we can expect these terms are certain
limits of higher genus Gromov-Witten invariants for a noncompact
Calabi-Yau $3$-fold, in this case the canonical bundle of
$\proj^1\times\proj^1$. More precisely, we put $\ve_1 = -\ve_2 =
\hbar$ and consider
\begin{equation*}
   F = F_0 + F_1\hbar^2 + F_2 \hbar^4 + \cdots.
\end{equation*}
Then $F_g$ is a limit of the genus $g$ Gromov-Witten invariants. Since
$1/\ve_1\ve_2 F$ is more fundamental (see \eqref{eq:prepotential2}),
we should write this as
\begin{equation*}
   \frac1{\hbar^2} F = \sum_{g=0}^\infty \hbar^{2g-2} F_g.
\end{equation*}
This is more natural as $2-2g$ is the Euler number of a genus $g$
Riemann surface. It probably explains the singularity $1/\hbar^2$.

Recently many Gromov-Witten invariants for noncompact toric Calabi-Yau
have been calculated (see \cite{AKMV} and the references
therein). These are identified with the Jones-Witten invariants via
the geometric transition (called `large $N$ duality'), as first
proposed by Gopakumar-Vafa \cite{GV}. A first of such examples is the
identification of Gromov-Witten for the resolved conifold and the
$\SU(N)$ Jones-Witten invariant for $S^3$. These identifications have
been proved in a mathematical rigorous way in a number of examples
(see \cite{OP,Zhou}).

In the case of $K_{\proj^1\times\proj^1}$, the Jones-Witten side is
$\SU(N)$-invariants for the Hopf link. Using the calculation by
Morton-Lukac \cite{ML}, Iqbal+Kashani-Poor show that the invariants of
the Hopf link has the same combinatorial expression as that of $F$
given by the localization formula \cite{IK}. (See also \cite{EK}.)

Note that these results identify the $n$-instanton correction with the
Gromov-Witten invariants of degree $n$ (with respect to one of the
factors of $\proj^1\times\proj^1$) for each $n$. Thus they do not say
much about the structure of the generating function $F$, which is
studied in this paper. Therefore it is interesting to understand the
blowup equation from the Gromov-Witten side.

\section{Seiberg-Witten curves}\label{sec:SW}
In this section we introduce the Seiberg-Witten curves, give the
definition of the prepotential, and derive the renormalization
equation and the contact term equation, which will characterize the
prepotential.

We give some details, though one can find most of them in physics
literature. The reason is that we must carefully choose cycles on the
Seiberg-Witten curve to determine a characteristic of the theta
function in a mathematically rigorous way. It is a standard exercise
but we cannot find the argument in the literature.

The material discussed here is a minimum of the Seiberg-Witten
geometry. We omit many things, such as monodromies, Picard-Fuchs
equations, relations to integrable systems, etc. Even for the
differential equations satisfied by the prepotential, our treatment is
a minimum. The Whitham hierarchy underlying these equations will not
be discussed. The reader may wonder where these equations come from,
though the authors' approach through the instanton moduli spaces will
be explained in \secref{sec:blowupeq}. For the original approaches,
see \cite{Mar,Marshakov} and the references therein.

There is a nice survey article \cite{Donagi} for mathematicians which
describes relation between integrable systems and the Seiberg-Witten
geometry, as well as background on physics. We recommend it to our
reader since it has no overlaps with this paper.

Note that we multiply $a_\alpha$ by $-\sqrt{-1}$ from the conventional
one in order to match with one in the instanton counting.

\subsection{Definition of the Seiberg-Witten prepotential}
We consider a family of curves (Riemann surfaces) parametrized by
$\vec{u} = (u_2,\dots,\linebreak[0]u_r)$:
\begin{equation*}
   C_{\vec{u}} : \Lambda^r\left(w + \frac1w\right)
   = P(z) = z^r + u_2 z^{r-2} + u_3 z^{r-3} + \cdots + u_r.
\end{equation*}
We call them {\it Seiberg-Witten curves}. The projection
$C_{\vec{u}}\ni (w,z)\mapsto z\in \proj^1$ gives a structure of
hyperelliptic curves. The hyperelliptic involution $\iota$ is given by
$\iota(w) = 1/w$.

If we introduce $y = \Lambda^r(w - \frac1w)$, we have 
\[
   y^2 = P(z)^2 - 4\Lambda^{2r} 
   = (P(z) - 2\Lambda^r)(P(z) + 2\Lambda^r).
\]
This special form of the right hand side will play a crucial role later.

The parameter space $\{ \vec{u} \in \C^{r-1}\}$ is called the {\it
$u$-plane}. Here $\Lambda$ is also a parameter, but we treat it
separately from $\vec{u}$. The parameter $\Lambda$ is called the {\it
renormalization scale\/} in physics. When $\Lambda = 0$, the theory
goes to the {\it classical limit}.
We consider  $\vec{u} = (u_2,\dots,u_r)$ as a coordinate system on the
$u$-plane. This is a global coordinate.

Let $z_1,\dots,z_r$ be the solutions of $P(z) = 0$. We will work on a
region of the $u$-plane where $|z_\alpha - z_\beta|$, $|z_\alpha|$ are
much larger than $|\Lambda|$, and then analytically continue. In
particular $z_\alpha$'s are distinct.  The vector $\vec{z} =
(z_1,\dots,z_r)$ ($\sum_\alpha z_\alpha = 0$) is a local coordinate on
the $u$-plane. The relation between $\vec{z}$ and $\vec{u}$ is very
simple. The former is a coordinate on $\C^{r-1}$ while the latter is
on $\C^{r-1}/S_r\approx\C^{r-1}$, where $S_r$ is the symmetric group
of $r$ letters. In other words, $(-1)^p u_p$ is the $p$th elementary
symmetric function in $z_1$, \dots, $z_r$. It is better to keep this
simple relation in mind, since this coordinate system $\vec{z}$ is a
{\it quantum correction\/} of another coordinate system $\vec{a}$
introduced below.

We can find $z_\alpha^\pm$ near $z_\alpha$ such that
$P(z_\alpha^\pm) = \pm 2\Lambda^r$ when $|u| \gg |\Lambda|$. These
are the $2r$-branched points of the projection
$C_{\vec{u}}\to\proj^1$. The infinity is not a branched point, and its
inverse image consists of $\infty_+$ ($w=\infty$) and $\infty_-$
($w=0$). The genus of $C_{\vec{u}}$ is $r-1$.
In the classical limit $\Lambda\to 0$, both $z_\alpha^\pm$ go to
$z_\alpha$, and the curves develop singularities. 

Let us define the {\it quantum discriminant\/} by
\begin{equation}
\label{eq:desc}
    \Delta
    = \left(4\Lambda^r\right)^{2r}
      \prod_{\alpha < \beta} (z_\alpha^+ - z_\beta^+)^2
    (z_\alpha^- - z_\beta^-)^2.
\end{equation}
On the locus $\Delta = 0$, the Seiberg-Witten curves develop
singularities. As we mentioned, we study a region away from this
locus.
\begin{NB}
I miss the term $\left(4\Lambda^r\right)^{2r}$ in the definition in
Oct.\ 25, this comes from $\prod (z_\alpha^+-z_\beta^-)^2$.
\end{NB}

\begin{figure}[htbp]
\centering
\psset{xunit=.5mm,yunit=.5mm,runit=.5mm}
\begin{pspicture}(0,0)(166.00,90.00)
\psellipse[linewidth=1.00,linecolor=black](80.00,45.00)(80.00,45.00)
\psellipse[linewidth=1.00,linecolor=black](45.00,45.00)(25.00,15.00)
\psellipse[linewidth=1.00,linecolor=black](115.00,45.00)(25.00,15.00)
\psline[linewidth=1.00,linecolor=black]{-}(70.00,45.00)(71.00,50.00)
\psline[linewidth=0.45,linecolor=black]{-}(90.00,45.00)(89.00,50.00)
\psline[linewidth=1.00,linecolor=black]{-}(140.00,45.00)(141.00,50.00)
\psellipse[linewidth=0.50,linecolor=black,linestyle=dotted,dotsep=1.00](150.00,45.00)(10.00,4.00)
\psline[linewidth=1.00,linecolor=black]{-}(20.00,45.00)(19.00,50.00)
\psellipse[linewidth=0.50,linecolor=black,linestyle=dotted,dotsep=1.00](10.00,45.00)(10.00,4.00)
\psline[linewidth=1.00,linecolor=black]{-}(90.05,44.60)(89.05,49.60)
\psellipse[linewidth=0.50,linecolor=black,linestyle=dotted,dotsep=1.00](80.00,45.00)(10.00,4.00)
\psellipse[linewidth=0.50,linecolor=black](45.50,45.00)(30.50,19.00)
\psellipse[linewidth=0.50,linecolor=black](80.00,45.00)(70.00,35.00)
\rput(-5.00,45.00){$z_1^+$}
\rput(27.50,45.00){$z_1^-$}
\rput(63.75,45.00){$z_2^-$}
\rput(97.50,45.00){$z_2^+$}
\rput(60.00,45.00){}
\rput(132.50,45.00){$z_3^+$}
\rput(140.00,45.00){}
\rput(166.00,45.00){$z_3^-$}
\rput(165.00,45.00){}
\rput(7.00,54.00){$A_1$}
\rput(80.00,54.00){$A_2$}
\rput(152.00,54.00){$A_3$}
\rput(45.00,68.00){$B_2$}
\rput(80.00,85.00){$B_3$}
\psline[linewidth=0.50,linecolor=black]{->}(46.25,26.03)(46.26,26.03)
\psline[linewidth=0.50,linecolor=black]{->}(80.99,10.00)(81.00,10.00)
\psline[linewidth=0.50,linecolor=black]{<-}(79.61,41.21)(79.62,41.21)
\psline[linewidth=0.50,linecolor=black]{<-}(6.42,41.25)(6.85,41.21)
\psline[linewidth=0.50,linecolor=black]{<-}(145.91,41.42)(146.34,41.38)

\linethickness{0.60pt}
{\renewcommand{\xscale}{1}
\renewcommand{\yscale}{.4}
\scaleput(149.50,111.9){\arc[10](10.00,0.00){180}} 
}
\linethickness{0.60pt}
{\renewcommand{\xscale}{1}
\renewcommand{\yscale}{.4}
\scaleput(10.00,111.9){\arc[10](10.00,0.00){180}} 
}
\linethickness{0.60pt}
{\renewcommand{\xscale}{1}
\renewcommand{\yscale}{.4}
\scaleput(80.00,111.9){\arc[10](10.00,0.00){180}} 
}
\end{pspicture}
\caption{Seiberg-Witten curve and cycles ($r=3$)}
\label{fig:SW}
\end{figure}
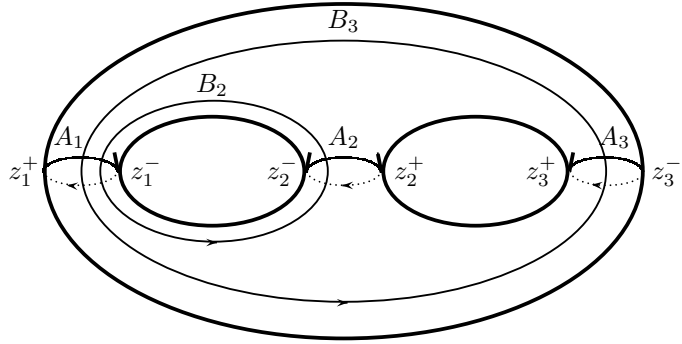

The hyperelliptic curve $C_{\vec{u}}$ is made of two copies of the
Riemmann sphere, glued along the $r$-cuts between $z_\alpha^-$ and
$z_\alpha^+$ ($\alpha = 1,\dots, r$), as usual.
Let $A_\alpha$ be the cycle encircling the cut between $z_\alpha^-$
and $z_\alpha^+$. We have $\sum_\alpha A_\alpha = 0$. We draw
$C_{\vec{u}}$ as in Figure~\ref{fig:SW}. The hyperelliptic involution
$\iota$ is the rotation by $\pi$ about the axis passing through the
branched points $z_\alpha^\pm$. Then we choose cycles $B_\alpha$
($\alpha=2,\dots,r$) as in Figure~\ref{fig:SW} so that $\{ A_\alpha,
B_\alpha \mid \alpha=2,\dots, r\}$ form a symplectic basis of
$H_1(C_{\vec{u}},\Z)$, i.e., $A_\alpha \cdot A_\beta = 0 =
B_\alpha\cdot B_\beta$, $A_\alpha\cdot B_\beta = \delta_{\alpha\beta}$
for $\alpha,\beta=2,\dots, r$. (The cycle $A_1$ is omitted.)
In the figure the branched points are lined as $z_1^+, z_1^-, z_2^-,
z_2^+, \cdots$ from the left. That is $z_\alpha^+$ is on the left
(resp.\ right) of $z_\alpha^-$ for $\alpha$ odd (resp.\ odd).

Note that we {\it cannot\/} take $A$, $B$-cycles {\it globally\/} on
the $u$-plane. The cycles are transformed by monodromies around the
locus $\Delta = 0$. In fact, the study of monodromies is important as
it has been used for constancy checks of the Seiberg-Witten curves to
some physically expected properties of the prepotential (introduced
below).  However we do not study monodromy behavior here except that
around $\Lambda = 0$. We first fix a small region in the $u$-plane and
then analytically continue. We choose a region containing the part
that $z_\alpha$'s are real and satisfy $z_1 > z_2 > \dots > z_r$. We
also assume $\Lambda$ is a positive real number. Since we assume
$\Lambda$ small, we have $z_1^+ > z_1 > z_1^- > z_2^- > z_2 > z_2^+ >
\cdots$.  This choice determines $A$, $B$-cycles as in
Figure~\ref{fig:SW}. Note the branched points are lined from the right
by the order in Figure~\ref{fig:SW}. Thus the choice is natural in
this region. Note also that we choose the inverse image of the region
with respect to the projection $\C^{r-1} \to \C^{r-1}/S_r$. The
permutation ambiguity is less important than the monodromies, but we
use the choice as we want to specify what is $a_\alpha$.

Let us define the {\it Seiberg-Witten differential\/} by 
\begin{equation*}
   dS = - \frac1{2\pi 
                 }z \frac{dw}w
   = -\frac1{2\pi 
     }
     \frac{z P'(z) dz}{\sqrt{P(z)^2  - 4\Lambda^{2r}}}
   = -\frac1{2\pi 
       } \frac{z P'(z) dz}y
   .
\end{equation*}
It is a meromorphic differential having poles at $\infty_\pm$. We
define functions $a_\alpha$, $a^D_\beta$ on the $u$-plane
($|u|\gg|\Lambda|$) by
\begin{equation}\label{eq:def_a}
   a_\alpha = \int_{A_\alpha} dS, \qquad
   a^D_\beta = 2\pi \sqrt{-1} \int_{B_\beta} dS,
   \qquad \alpha=1,\dots,r,\ \beta=2,\dots,r.
\end{equation}

Let us study the behavior of the function $a_\alpha$ around $\Lambda
= 0$. We move the cycle $A_\alpha$ so that $P(z)$ and $1/P(z)$ are
bounded there. In particular, we are in a sheet where $\sqrt{P(z)^2 -
4\Lambda^{2r}}$ is single-valued. We choose the sheet so that it is
approximated by $P(z)$ on the $A_\alpha$-cycle. We suppose $A_\alpha$
has the counterclockwise rotation in the sheet. In
Figure~\ref{fig:SW}, the sheet is the part lower than the plane
containing $A_\alpha$'s. (See also the proof of \propref{prop:F_beh}
below.)
Then we have the following expansion:
\begin{equation}\label{eq:z_expand}
\begin{split}
   a_\alpha &= -\frac1{2\pi
                      }
     \int_{A_\alpha} dz\, z \frac{P'(z)}{P(z)}
     \sum_{n=0}^\infty \binom{\frac12}{n}
     \left(\frac{4\Lambda^{2r}}{P(z)}\right)^n
   = -\sqrt{-1} z_\alpha + O(\Lambda^{2r}),
\end{split}
\end{equation}
where $\binom{a}{n}$ is the binomial coefficient. In particular,
$\vec{a} = (a_1,\dots, a_r)$ ($\sum_\alpha a_\alpha = 0$) is a local
coordinate system for small $\Lambda$. As we mentioned before, the
coordinate $\vec{z}$ is the quantum correction of $\vec{a}$. The
$a_\alpha$ is a function in $u_p$, but conversely we consider $u_p$ as
a function in $a_\alpha$ (and also in $\Lambda$).

We differentiate the Seiberg-Witten differential $dS$ by setting $w$
to be constant:
\begin{equation*}
   \left.\frac{\partial}{\partial u_p} dS\right|_{w=\mathrm{const}}
   = \frac1{2\pi 
              } \frac{z^{r-p}}{P'(z)}\frac{dw}w
   = \frac1{2\pi
              } \frac{z^{r-p} dz}y.
\end{equation*}
\begin{NB}
There was a sign mistake in the earlier version (Oct.\ 16).
\end{NB}
It is well-known that these form a basis of holomorphic differentials
on $C_{\vec{u}}$ for $p=2,\dots,r$ (see e.g., \cite[\S2.3]{GH}). In
other words, the Seiberg-Witten differential is a `{\it potential\/}'
for holomorphic differentials.
Let $(\sigma_{\alpha p})$ be the matrix given by
\begin{equation*}
   \sigma_{\alpha p}
   = \pd{a_\alpha}{u_p} 
   = \frac1{2\pi 
         } \int_{A_\alpha} \frac{z^{r-p}}{P'(z)}\frac{dw}w,
   \qquad \alpha, p =2,\dots,r.
\end{equation*}
If $(\sigma^{p\alpha})$ is the inverse matrix, the normalized
holomorphic $1$-forms 
\[
  \omega_\beta
  = \frac1{2\pi 
        } \sum_p \sigma^{p\beta}
    \frac{z^{r-p}}{P'(z)}\frac{dw}w
  = \left.\pd{}{a_\beta} dS\right|_{w=\mathrm{const}}
\]
satisfies $\int_{A_\alpha} \omega_\beta =
\delta_{\alpha\beta}$. Therefore the {\it period matrix\/} $\tau =
(\tau_{\alpha\beta})$ of the curve $C_{\vec{u}}$ is given by
\begin{equation*}
  \tau_{\alpha\beta}
  = \int_{B_\alpha} \omega_\beta
  = \frac1{2\pi\sqrt{-1}}\pd{a^D_\alpha}{a_\beta}.
\end{equation*}
Since $(\tau_{\alpha\beta})$ is symmetric (see e.g.,
\cite[\S2.2]{GH}), there exists a locally defined function $\mathcal
F_0$ on the $u$-plane such that
\begin{equation}\label{eq:SWprep}
   a_\alpha^D = - \frac{\partial\mathcal F_0}{\partial a_\alpha}.
\end{equation}
It is unique up to constant. We fix the constant so that $\mathcal
F_0$ is homogeneous of degree $2$:
\begin{equation}\label{eq:Euler}
   \left(\sum a_\alpha \frac{\partial}{\partial a_\alpha}
   + \Lambda \pd{}{\Lambda} \right) \mathcal F_0
 = 2\mathcal F_0.
\end{equation}
This function $\mathcal F_0$ is called the {\it Seiberg-Witten
prepotential}. We may also write $\mathcal F_0(\vec{a})$ or $\mathcal
F_0(\vec{a};\Lambda)$.
From the definition we have
\begin{equation}\label{eq:period}
   \tau_{\alpha\beta}
   = - \frac1{2\pi\sqrt{-1}}
     \frac{\partial^2 \mathcal F_0}{\partial a_\alpha \partial a_\beta}.
\end{equation}
We put the subscript $0$ because this will be identified with the
genus $0$ part of the Nekrasov's deformed partition function.

\subsection{The logarithmic singularities of the Seiberg-Witten
   prepotential}

Let us study the behavior of $\mathcal F_0$ as a
function in $\Lambda$, following \cite{HKP}. Our aim is to show
\begin{Proposition}\label{prop:F_beh}
We have
\begin{equation*}
  \begin{split}
    \mathcal F_0
   & = \sum_{\alpha \neq \beta}
    \gamma_0(a_\alpha-a_\beta;\Lambda) + O(\Lambda^{2r})
\\
   & =  \sum_{\alpha \neq \beta}\left[ \frac12 (a_\alpha-a_\beta)^2
    \log\left(\frac{a_\alpha-a_\beta}{\Lambda}\right)
    -\frac34 (a_\alpha-a_\beta)^2\right] + O(\Lambda^{2r})
   ,
  \end{split}
\end{equation*}
where $\gamma_0(x;\Lambda) = \frac12
x^2\log\left(\frac{x}\Lambda\right) - \frac34 x^2$ is the coefficient
of $1/\ve_1\ve_2$ in $-\gamma_{\ve_1,\ve_2}(x;\Lambda)$ as in
\eqref{eq:pert_expand}.
\end{Proposition}

The part $\sum_{\alpha \neq \beta}\gamma_0(a_\alpha-a_\beta;\Lambda)$
is called the {\it perturbative part\/} of the prepotential $\mathcal
F_0$, and denoted by $\mathcal F_0^{\mathrm{pert}}$. The remaining
part is called the {\it instanton part}, and denoted by $\mathcal
F_0^{\mathrm{inst}}$. It is a power series in $\Lambda^{2r}$:
\(
   \mathcal F_0^{\mathrm{inst}} = f_1 \Lambda^{2r} + f_2 \Lambda^{4r}
   + \cdots + f_n \Lambda^{2rn} + \cdots.
\)
The coefficient $f_n$ is called the {\it $n$th instanton correction\/}
to the prepotential. This is because we will identify $f_n$ something
defined via the $n$-instanton moduli space.

The choice of the branch of $\log$ is as follows. Suppose that
$z_\alpha$, $z_\alpha^\pm\in \R$, $\Lambda\in\R_{>0}$ and $z_1 > z_2 >
\dots > z_r$ as above.
We choose a path $A_\alpha$ encircling $z_\alpha^-$ and $z_\alpha^+$
so that it is invariant under the complex conjugation $z\mapsto
\overline{z}$. Then $a_\alpha$ is purely imaginary. We have
\(
    \sqrt{-1} a_1 > \sqrt{-1} a_2 > \cdots > \sqrt{-1} a_r.
\)
We choose the branch of $\log$ so that
\begin{equation*}
   \frac12 
    \log\left(\frac{a_\alpha-a_\beta}{\Lambda}\right)
   +
   \frac12 
    \log\left(\frac{a_\beta-a_\alpha}{\Lambda}\right)
   =
    \log\left(\frac{\sqrt{-1}(a_\alpha-a_\beta)}{\Lambda}\right)
\end{equation*}
is real for $\alpha < \beta$. In what follows, we assume this choice
of $z_\alpha$, etc. It is enough to consider this case by analytic
continuation.
\begin{NB}
Since we have fixed $A$, $B$-cycles, the functions $a_\alpha$,
$a^D_\beta$ are determined. Thus $\mathcal F_0$ should be determined
without ambiguity. Thus the branch of $\log$ is automatically
determined. In the version Oct.\ 25, I wrote that I fix the branch
only to check the choice of cycles. This may be confusing.
\end{NB}

\begin{proof}[Proof of \propref{prop:F_beh}]
%
First we study
\begin{equation*}
   a^D_\alpha
   = 
   2 \pi \sqrt{-1} \int_{B_\alpha} dS
   .
\end{equation*}
Locally, this is a function in $\Lambda^{2r}$. But it is multi-valued,
as the cycle $B_\alpha$ transforms to $B_\alpha+A_\alpha-A_1$ when we
analytically continue from $\Lambda^{2r}$ to
$e^{2\pi\sqrt{-1}}\Lambda^{2r}$. Therefore
\(
   a^D_\alpha + \sqrt{-1}(z_\alpha - z_1) \log \Lambda^{2r}
\)
is a single-valued function in $\Lambda$. This kind of the monodromy
behavior is quite important in the conventional arguments.

For $\beta=2,\dots, r$, let $C'_\beta$ be the straight line from
$z_{\beta-1}^\pm$ to $z_{\beta}^\pm$ ($+$ for $\beta$ odd, $-$ for
$\beta$ even) in one sheet and define the cycle $C_\beta$ as
$C'_\beta$ followed by $-\iota C'_\beta$. This is a cycle rounding the
hole in Figure~\ref{fig:SW}. We have
\begin{equation*}
   B_\alpha = \sum_{\beta=2}^\alpha C_\beta,
\end{equation*}
and
\begin{equation*}
   -2\pi\int_{C_\beta} dS
   = 2 \int^{z_\beta^\pm}_{z_{\beta-1}^\pm} 
      \frac{z P'(z) dz}{\sqrt{P(z)^2  - 4\Lambda^{2r}}}.
\end{equation*}
In the last expression, $z$ is real. But we should be careful for the
choice of the branch of $\sqrt{P(z)^2  - 4\Lambda^{2r}}$. This is not
necessarily $\ge 0$ contrary to the usual convention for the real
function. It is determined by the analytic continuation. We choose the
sheet so that $\sqrt{P(z)^2 - 4\Lambda^{2r}}$ has the same sign as
$P(z)$ on each interval $[z_{\beta-1}^\pm,z_\beta^\pm]$. This is
the same sheet used in \eqref{eq:z_expand} (i.e., the lower half) and
we have the right orientation so that $A_\alpha\cdot B_\alpha = 1$.
Note also that $a^D_\alpha$ is pure imaginary as the integrals are
real.

Fix $\delta > 0$ small with $|\Lambda|\ll |\delta|$ and rewrite the
integral as
\begin{equation}\label{eq:a_D}
   -\pi \int_{C_\beta} dS
   = 
   \left(
   \int^{z_\beta-\delta}_{z_{\beta-1}+\delta}
   +
   \int_{z_\beta-\delta}^{z_\beta^\pm}
   - 
   \int_{z_{\beta-1}+\delta}^{z_{\beta-1}^\pm}
   \right)
     \frac{z P'(z) dz}{\sqrt{P(z)^2  - 4\Lambda^{2r}}}.
\end{equation}
The first integral is regular at $\Lambda=0$:
\begin{equation*}
  \begin{split}
   & \int^{z_\beta-\delta}_{z_{\beta-1}+\delta} \frac{z P'(z) dz}{P(z)}
      + O(\delta)
   = \int^{z_\beta-\delta}_{z_{\beta-1}+\delta}
     \sum_{\gamma=1}^r (1 + \frac{z_\gamma}{z - z_\gamma}) dz
      + O(\delta)
\\
   =\; & 
   \begin{aligned}[t]
   & z_{\beta} (r + \log\delta) - z_{\beta-1} ( r + \log \delta) 
\\
   & \qquad\qquad
    + \sum_{\gamma\neq\beta} z_\gamma \log|z_\beta - z_\gamma|
   - \sum_{\gamma\neq\beta-1} z_\gamma \log|z_{\beta-1} - z_\gamma|
   + O(\delta).
   \end{aligned}
  \end{split}
\end{equation*}
Here we choose the branch of $\log$ so that all the above expressions
are real.
\begin{NB}
There is a mistake in \cite[2.14]{HKP}, where $r$ is replaced by $1$.
\end{NB}

The second integral of \eqref{eq:a_D} is equal to
\begin{equation}\label{eq:II}
    z_\beta \int_{z_\beta-\delta}^{z_\beta^\pm}
    \frac{P'(z) dz}{\sqrt{P(z)^2  - 4\Lambda^{2r}}}
    + \int_{z_\beta-\delta}^{z_\beta^\pm}
    \frac{(z - z_\beta) P'(z) dz}{\sqrt{P(z)^2  - 4\Lambda^{2r}}}.
\end{equation}
The first term is
\begin{equation*}
    z_\beta \int_{w = w_\beta(\delta)}^{w=\pm 1} \frac{d w}w
    = - z_\beta \log |w_\beta(\delta)|,
\end{equation*}
where
\begin{equation*}
   w_\beta(\delta) 
   = \frac1{2\Lambda^r}\left(P(z_\beta-\delta) +
     \sqrt{P(z_\beta-\delta)^2 - 4\Lambda^{2r}}\right)
    .
\end{equation*}
\begin{NB}
Since $\Lambda^r(w + \frac1w) = P(z)$, we have $P(z_\beta^\pm) =
\pm 2\Lambda^r = \Lambda^r (\pm 1 + \frac1{\pm 1})$. Therefore
$z=z_\alpha^\pm \leftrightarrow w=\pm 1$.
\end{NB}
We have
\begin{equation*}
  \begin{split}
   & \log |w_\beta(\delta)|
   = \log \frac{|P(z_\beta-\delta)|}{2\Lambda^r}
     + \log\left( 1 + \sqrt{1 - 
       \frac{4\Lambda^{2r}}{P(z_\beta-\delta)^2}}\right)
\\
   = 
   \; &
   \log\left(\frac{\left|\delta
      \prod_{\gamma\neq\beta}(z_\beta-z_\gamma)\right|}{\Lambda^r}\right)
      + O(\delta).
  \end{split}
\end{equation*}
Let us consider the second term of \eqref{eq:II}. Let $z - z_\beta =
\prod_{\gamma\neq\beta} (z_\beta- z_\gamma)^{-1} P(z) + E(z)$. We
have $E(z) = O(\delta^2)$ in the range of the integration. But the
above calculation of the first part shows that the integral of $E(z)$
yields $O(\delta^2)O(\log\delta) = O(\delta)$. Therefore the second
term is
\begin{equation*}
  \begin{split}
    & \frac1{\displaystyle\prod_{\gamma\neq\beta} (z_\beta- z_\gamma)}
    \int^{z_\beta-\delta}_{z_\beta^\pm}
    \frac{P(z) P'(z) dz}{\sqrt{P(z)^2  - 4\Lambda^{2r}}} 
    + O(\delta)
\\
    = \; &
    \frac1{\displaystyle\prod_{\gamma\neq\beta} (z_\beta- z_\gamma)}
    \left[ \sqrt{P(z)^2  - 4\Lambda^{2r}}
       \right]^{z_\beta-\delta}_{z_\beta^\pm}
    + O(\delta).
  \end{split}
\end{equation*}
But as $P(z_\beta^\pm) = 2\Lambda^{r}$, $P(z_\beta-\delta) =
O(\delta)$, the contribution is $O(\delta)$.

The third integral has a similar expression with $z_\beta$, $\delta$
replaced by $z_{\beta-1}$, $-\delta$ respectively. Altogether we get
\begin{equation*}
  \begin{aligned}[b]
   & -\pi \int_{C_\beta} dS 
   - r(z_\beta - z_{\beta-1})\left(1 + \log\Lambda\right)
\\
   & \qquad
   + \sum_{\substack{\gamma\neq\beta}} 
   (z_\beta-z_\gamma)\log|z_\beta-z_\gamma|
   - \sum_{\substack{\gamma\neq \beta-1}} 
   (z_{\beta-1}-z_\gamma)\log|z_{\beta-1}-z_\gamma|
  \end{aligned}
 = O(\delta).
\end{equation*}
But the left hand side is independent of $\delta$.
This means that the left hand side is, in fact, $O(\Lambda^{2r})$.
Combining with \eqref{eq:z_expand}, we have
\begin{equation*}
   \frac12 a^D_\alpha
   =
   \begin{aligned}[t]
   & r (a_\alpha - a_1)\left(1 + \log\Lambda\right) - 
    \sum_{\substack{\beta\neq\alpha}} 
   (a_\alpha-a_\beta)\log\left|\sqrt{-1}(a_\alpha-a_\beta)\right|
   \\
   & \qquad
   + \sum_{\substack{\beta\neq 1}} 
   (a_1-a_\beta)\log\sqrt{-1}(a_1-a_\beta) + O(\Lambda^{2r}).
   \end{aligned}
\end{equation*}
In the last part, we do not take the absolute value of
$\sqrt{-1}(a_1-a_\beta)$ since $\sqrt{-1}a_1 > \sqrt{-1}a_\beta$.

Now let us differentiate $\mathcal F_0^{\mathrm{pert}}$ in the
statement. Let $\overline{\gamma}_0(x;\Lambda) = x^2
\log\left(\frac{\sqrt{-1}x}\Lambda\right) - \frac32 x^2$. (Remember
our choice of the branch of $\log$.) We have the following
\begin{equation*}
\begin{split}
   & - \pd{\mathcal F_0^{\mathrm{pert}}}{a_\alpha}
   = - \sum_{\beta < \gamma} \pd{}{a_\alpha}
   \overline{\gamma}_0(a_\beta-a_\gamma;\Lambda)
\\
   & =
   \begin{aligned}[t]
   &
   - \sum_{\substack{\alpha < \beta}} 
   \overline{\gamma}_0'(a_\alpha-a_\beta;\Lambda)
   + \sum_{\substack{\beta < \alpha}} 
    \overline{\gamma}_0'(a_\beta-a_\alpha;\Lambda)
   + \sum_{\substack{1 \neq \beta}} 
    \overline{\gamma}_0'(a_1 - a_\beta;\Lambda).
   \end{aligned}
  \end{split}
\end{equation*}
\begin{NB}
More detail:
\(
   \mathrm{LHS} =
   \begin{aligned}[t]
   & - \sum_{\substack{\alpha \neq \gamma}}
      \overline{\gamma}_0'(a_\alpha-a_\gamma;\Lambda)
     + \sum_{\substack{1\neq \beta < \alpha}}
     \overline{\gamma}_0'(a_\beta-a_\alpha;\Lambda)
   \\
   & \qquad\qquad
   + \sum_{\substack{\gamma \neq \alpha \\ \gamma\neq 1}} 
     \overline{\gamma}_0'(a_1 - a_\gamma;\Lambda)
   + 2 \overline{\gamma}_0'(a_1-a_\alpha;\Lambda)
   \end{aligned}
   = \mathrm{RHS}
\)
\end{NB}
We substitute
\(
    \overline{\gamma}_0'(x) = 2x\log\frac{\sqrt{-1}x}\Lambda - 2x
\)
to get
\begin{equation*}
\begin{split}
   - \frac12 \pd{\mathcal F_0^{\mathrm{pert}}}{a_\alpha} &=
     \begin{aligned}[t]
    &
   r (a_\alpha - a_1)
  -
     \sum_{\substack{\beta\neq\alpha}} 
   (a_\alpha-a_\beta)\log\frac{\left|\sqrt{-1}(a_\alpha-a_\beta)\right|}\Lambda
   \\
   & \qquad
   + \sum_{\substack{\beta\neq 1}} 
   (a_1-a_\beta)\log\frac{\sqrt{-1}(a_1-a_\beta)}\Lambda
    .
   \end{aligned}
\end{split}
\end{equation*}
\begin{NB} More detail:
\(
   \mathrm{LHS} =
   \begin{aligned}[t]
     & 
   \sum_{\substack{\alpha \neq \beta}} 
    (a_\alpha-a_\beta)
   - \sum_{\substack{\beta\neq 1}} 
    (a_1-a_\beta)
   -  
   \sum_{\substack{\alpha \neq \beta}} 
    (a_\alpha-a_\beta)
   \log\frac{\left|\sqrt{-1}(a_\alpha-a_\beta)\right|}{\Lambda}
   \\
   & \qquad
   + \sum_{\substack{\beta\neq 1}} 
    (a_1-a_\beta) \log\frac{\sqrt{-1}(a_1-a_\beta)}{\Lambda}
   \end{aligned}
   = \mathrm{RHS}.
\)
\end{NB}
This coincides with the above expression. The proof of
\propref{prop:F_beh} is completed.
\end{proof}

\subsection{A renormalization group equation}

We prove the so-called `renormalization group equation' following
\cite{STY} in this subsection:
\begin{Proposition}\label{prop:renorm}
\begin{equation*}
   \Lambda\pd{}{\Lambda}\mathcal F_0 = - 2r u_2.
\end{equation*}
\begin{NB}
There was a sign mistake in the earlier version (Oct. 16).
\end{NB}
\end{Proposition}

This equation was found earlier by \cite{Matone} for $\SU(2)$, and
independently by \cite{EY}. See also \cite{HKP1}.

\begin{proof}
We differentiate the Euler equation \eqref{eq:Euler}:
\begin{equation*}\label{eq:diffEuler}
\begin{split}
   &\pd{}{u_p}\left(\Lambda\pd{}{\Lambda}\mathcal F_0\right)
   = 2\pd{\mathcal F_0}{u_p} - \sum_\alpha 
   \pd{}{u_p}\left(a_\alpha\pd{\mathcal F_0}{a_\alpha}\right)
   = -\sum_\alpha \left(\pd{a_\alpha}{u_p}a^D_\alpha
     - a_\alpha\pd{a^D_\alpha}{u_p}\right)
\\
  =\; &
  - 2\pi\sqrt{-1}
  \sum_\alpha \left[\int_{A_\alpha} \pd{}{u_p} dS \int_{B_\alpha} dS
  - \int_{A_\alpha} dS \int_{B_\alpha} \pd{}{u_p} dS
  \right]
  .
\end{split}
\end{equation*}
Let us make a change of variable $x = 1/z$. We expand the
Seiberg-Witten differential and its differential around $x=0$:
\begin{equation*}
\begin{split}
   & dS = -\frac1{2\pi 
     }
     \frac{z P'(z) dz}{\sqrt{P(z)^2  - 4\Lambda^{2r}}}
     = \left(s_{-2} x^{-2} + s_0 + s_1 x + \cdots\right) dx,
\\
  & \pd{}{u_p} dS = \frac1{2\pi 
     }
     \frac{z^{r-p}dz}{\sqrt{P(z)^2  - 4\Lambda^{2r}}}
   = \left(\omega^p_0 + \omega^p_1 x + \cdots\right) dx.
\end{split}
\end{equation*}
(Recall that $dS$ is a meromorphic differential having poles only at
$\infty_\pm$.) By the Riemann bilinear relation (see e.g.,
\cite[\S2.3]{GH}), we have
\begin{equation*}
   \pd{}{u_p}\left(\Lambda\pd{}{\Lambda}\mathcal F_0\right)
   = 8\pi^2\sum_n \frac{s_{-n}\omega^p_{n-2}}{n-1}
   = 8\pi^2 s_{-2} \omega_0^p.
\end{equation*}
Since
\(
   s_{-2} = \frac{r}{2\pi},
\)
\(
   \omega^p_0 = -\frac1{2\pi}\delta_{p2},
\)
we get
\begin{equation*}
   \pd{}{u_p}\left(\Lambda\pd{}{\Lambda}\mathcal F_0\right)
 = -2r \delta_{p2}.
\end{equation*}
Integrating out, we get the assertion. Here the integration constant
is zero thanks to the homogeneity of $\mathcal F_0$.
\end{proof}

\subsection{The contact term equation}

In this subsection we show the following partial differential equation:
\begin{Theorem}\label{thm:contact}
We have
\begin{equation*}
   \Lambda\pd{}{\Lambda}u_p
   =
   \frac{2r}{\pi\sqrt{-1}} \sum_{\alpha,\beta=2}^r
   \pd{u_p}{a_\alpha} \pd{u_2}{a_\beta}
   \pd{}{\tau_{\alpha\beta}}\log\Theta_E(0|\tau)
\end{equation*}
for $p=2,3,\dots,r$. Here $E$ is the even half-integer characteristic
given by
\(
\left[
\begin{smallmatrix}
   \vec{0} \\ \vec{\Delta}
\end{smallmatrix}
\right]
\)
in \eqref{eq:characteristic}.
\end{Theorem}

This equation was first derived by Losev-Nekrasov-Shatashvili
\cite{LNS1,LNS2} during their study of the topologically twisted
version of $\mathcal N=2$ SUSY Yang-Mills theory, i.e., the physical
counterpart of the Donaldson theory. More precisely, they derived the
equation by studying the effect of the blowup on the so-called
`contact terms'. So we call the equation the {\it contact term\/}
equation. Later Gorsky, Marshakov, Mironov and Morozov \cite{GM3}
derived the contact term equation in the framework of the
Seiberg-Witten curve. We give the proof following their approach in
this subsection. In fact, they did not determine the characteristic.
It was determined in \cite{LNS1,LNS2}, but the argument involves a
physical intuition. Here we can give a mathematically rigorous proof
thanks to our precise definition of the $B$-cycles used in the
definition of the prepotential.

For a later purpose, we give a remark.
Recall that $(-1)^p u_p$ is the $p$th elementary symmetric function in
in variables $z_1,\cdots, z_r$. Let $c_p$ be the $p$th power sum
multiplied by $\frac{(-\sqrt{-1})^p}{p!}$:
\begin{equation}\label{eq:c}
   c_p = \frac{(-\sqrt{-1})^p}{p!} \sum_{\alpha=1}^r {z_\alpha^p}, \quad
   c_1 = 0, c_2 = u_2, \cdots, \text{etc}.
\end{equation}
Then $(c_2,c_3,\dots,c_r)$ is another coordinate system on the
$\vec{u}$-plane. Since $c_p$ is a polynomial in $u_q$'s, it is also a
solution of the contact term equation:
\begin{equation*}
   \Lambda\pd{}{\Lambda}c_p
   =
   \frac{2r}{\pi\sqrt{-1}} \sum_{\alpha,\beta=2}^r
   \pd{c_p}{a_\alpha} \pd{u_2}{a_\beta}
   \pd{}{\tau_{\alpha\beta}}\log\Theta_E(0|\tau).
\end{equation*}

Before giving the proof of \thmref{thm:contact}, we give a corollary
which will play an important role later.
\begin{Corollary}\label{cor:contact}
\begin{equation*}
   \left(\Lambda\pd{}{\Lambda}\right)^2 \mathcal F_0
   =
   \frac{-1}{\pi\sqrt{-1}} \sum_{\alpha,\beta=2}^r
   \pd{}{a_\alpha}\left(\Lambda\pd{}\Lambda \mathcal F_0\right)
   \,
   \pd{}{a_\beta}\left(\Lambda\pd{}\Lambda \mathcal F_0\right)
   \pd{}{\tau_{\alpha\beta}}\log\Theta_E(0|\tau).
\end{equation*}
\begin{NB}
There was a sign mistake in the earlier version (Oct. 16).
\end{NB}
\end{Corollary}
This equation together with the description of the perturbative part
(\propref{prop:F_beh}) completely determines the prepotential
$\mathcal F_0$. See the proof of \thmref{thm:main} and
\subsecref{subsec:recursive} below. This observation was due to
\cite{EMM}. (See also \cite{Matone} for an earlier result for
$\SU(2)$.)

\begin{proof}[Proof of \thmref{thm:contact}]
%
Recall that we consider $u_p$ as functions of $a_\alpha$, $\Lambda$.
We differentiate \eqref{eq:def_a} by $\log\Lambda$ to get
\begin{equation*}
   \sum_p \pd{u_p}{\log\Lambda}\int_{A_\alpha} \pd{}{u_p} dS
   + \int_{A_\alpha} \pd{}{\log\Lambda} dS = 0.
\end{equation*}
Therefore
\begin{equation*}
\begin{split}
   & 
   \sum_p \pd{u_p}{\log\Lambda} \pd{a_\alpha}{u_p}
   = -\int_{A_\alpha} \pd{}{\log\Lambda} dS
   = \frac{r}{2\pi}\int_{A_\alpha} \frac{P(z)}{P'(z)}\frac{dw}{w}
   = \frac{r}{2\pi}\int_{A_\alpha} \frac{P(z) dz}{y}
\\
   =
   \; &
   \frac{r}{2\pi}\int_{A_\alpha} \frac{(P(z) + y)dz}{y}
   .
\end{split}
\end{equation*}
The last expression can be given by the Szeg\"o kernel (see
\eqref{eq:Szego}) as
\begin{equation*}
   \Psi_E^2(z_1,\infty_\pm)
   = -\frac{P(z_1)\pm y(z_1)}{2y(z_1)} dz_1
   \left.d\left(\frac1{z_2}\right)\right|_{z_2=\infty_\pm}.
\end{equation*}
Here we choose the leftmost point $z_1^+$ as the base point for the
Abel-Jacobi map. And the even half-integer characteristic $E$
corresponds to the partition of the branched points into
\begin{equation*}
   \{ z_\alpha^+ \mid \alpha=1,\dots,r\}
   \sqcup
   \{ z_\alpha^- \mid \alpha=1,\dots,r\}.
\end{equation*}

On the other hand, we have
\begin{equation*}
   \left.\omega_\beta\right|_{z=\infty_\pm}
   = -\frac1{2\pi} \sum_p \pd{u_p}{a_\beta}
   \left.\frac{z^{r-p} dz}{y}\right|_{z=\infty_\pm}
   = \frac1{2\pi} \pd{u_2}{a_\beta}
   \left.d\left(\frac1z\right)\right|_{z=\infty_\pm},
\end{equation*}
where we have used $y\sim z^{r}$ at $z=\infty_\pm$ in the second equality.
Therefore by Fay's identity \eqref{eq:Fay} we have
\begin{equation*}
  \begin{split}
    & \frac{r}{2\pi}\int_{A_\alpha} \frac{(P(z) + y)dz}{y}
   = - \frac{r}{2\pi^2} \pd{u_2}{a_\beta}
      \frac{\partial^2}{\partial \xi_\alpha\partial \xi_\beta}
     \log\Theta_E(0|\tau)
\\
   =\; &
   \frac{2r}{\pi\sqrt{-1}}\pd{u_2}{a_\beta}
   \pd{}{\tau_{\alpha\beta}}\log\Theta_E(0|\tau).
  \end{split}
\end{equation*}
In the second equality we have used the heat equation and the fact
that $E$ is an even half-integer characteristic, and hence the
derivative of $\Theta_E(\vec{\xi}|\tau)$ at $\vec\xi=0$ vanishes.
Multiplying the matrix $\left(\pd{u_p}{a_\alpha}\right)$ to both hand
sides, we get the differential equation as in the assertion.

Finally we determine the even half-integer characteristic $E$
explicitly. Looking at Figure~\ref{fig:SW}, we find that the partition
corresponding to characteristic $0$ is
\begin{equation*}
   \left\{ z_1^+, z_2^-, z_3^+, \cdots \right\}
   \sqcup
   \left\{ z_1^-, z_2^+, z_3^-, \cdots \right\}.
\end{equation*}
Namely $z_\alpha^+$ and $z_\alpha^-$ for $\alpha$ even are
interchanged from $E$. Since
\(
   \int_{z_\alpha^-}^{z_\alpha^+} \omega_\beta
   = \frac12 \int_{A_\alpha} \omega_\beta
   = \frac12 \delta_{\alpha\beta},
\)
we find that the characteristic $E$ is
\(
\left[
\begin{smallmatrix}
   \vec{0} \\ \vec{\Delta}
\end{smallmatrix}
\right]
\)
in \eqref{eq:characteristic}.
\end{proof}

\begin{Remark}
In \cite{GM3} `time variables' $T_1$, $T_2$, \dots, $T_{r-1}$ are
introduced in the framework of Whitham hierarchy. Then the contact
term equations are the specialization of the equations at $T_2 = T_3 =
\cdots = 0$ ($T_1$ is essentially $\log\Lambda$). On the other hand,
we will introduce {\it infinitely many\/} variables $\tau_1, \tau_2,
\dots$ in \secref{sec:Nek}. We will show
\(
   \pd{}{T_p} = \pd{}{\tau_p}
\)
for $p=1,2,\dots,r-1$ when it is restricted to
\(
   \tau_1 = \tau_2 = \cdots = 0
\)
in \thmref{thm:main}(2). However the equations for the
$\tau_p$-derivatives are not explicitly written down, and are
different from the equations for $T_p$-derivatives outside this
subspace.
\end{Remark}

\subsection{Rank $2$ case}\label{subsec:rank2}

When $r=2$, i.e., the Seiberg-Witten curve is an elliptic curve, we
have the expressions \eqref{eq:rank2} for $u$, $a$ in terms of theta
functions and Eisenstein series. Here we write $u = u_2$, $a = a_2$.
The derivation of the expressions are left to the reader as an
exercise. One can prove \thmref{thm:contact} using the
expressions. See \cite[Appendix]{GM3}.

\begin{NB}
In order to check the sign/typos we give the computation.
\begin{proof}
We may assume $\Lambda = 1$ by replacing $u\mapsto u/\Lambda^2$. The
Seiberg-Witten curve is
\begin{equation*}
   y^2 = (z^2 + u - 2)(z^2 + u + 2)
   = \left\{\left(\frac{z}{\sqrt{2-u}}\right)^2 - 1\right\}
   \left\{\left(\frac{z}{\sqrt{-2-u}}\right)^2 - 1\right\}
   (u^2 - 4).
\end{equation*}
When $u \ll 0$, the roots $\sqrt{\pm 2-u}$ are real and $\sqrt{2-u} >
\sqrt{-2-u} > -\sqrt{-2-u} > -\sqrt{2-u}$. The $A$-cycle encircles
$-\sqrt{-2-u}$ and $-\sqrt{2-u}$, while the $B$-cycle encircles
$\sqrt{-2-u}$ and $-\sqrt{-2-u}$.
This is the standard form of a quartic curve in \cite[\S4]{Ume} with
$y\mapsto \frac{y}{\sqrt{u^2-4}}$, $x = \frac{z}{\sqrt{-2-u}}$ and
$\kappa = \sqrt{\frac{-2-u}{2-u}}$. Note that the $A$-cycle in
\cite[\S4]{Ume} encircles $\pm 1$, so $A$, $B$-cycles are
interchanged. Note also that the curve \cite[\S4.5]{Ume} has period
$2\tau$, so the theta functions there should be $\theta_*(z|-2/\tau)$.
The change $\tau \mapsto -1/\tau$ makes $\kappa\mapsto \kappa'$. We
specify the period during the proof. We recall some formulas:
\begin{align*}
   & \theta_{00}(\tau)^4 = \theta_{01}(\tau)^4 +
   \theta_{10}(\tau)^4
   \tag*{\cite[(3.35)]{Ume}}
\\
   & \theta_{00}(\tau/2)\theta_{01}(\tau/2) = \theta_{01}(\tau)^2
   \tag*{\cite[(3.78)]{Ume}},
\\
   & \theta_{10}(\tau/2)^2 = 2\theta_{00}(\tau)\theta_{10}(\tau)
   \tag*{\cite[(3.79)]{Ume}},
\\
   & \theta_{00}(\tau/2)^2 + \theta_{01}(\tau/2)^2
   = 2\theta_{00}(\tau)^2.
   &\tag*{\cite[(3.90)]{Ume}}
\end{align*}

By \cite[(4.5)]{Ume} we have
\begin{equation*}
   \frac{-2-u}{2-u} = \kappa^{2}
   = \frac{\theta_{10}(-2/\tau)^4}{\theta_{00}(-2/\tau)^4}
   = \frac{\theta_{01}(\tau/2)^4}{\theta_{00}(\tau/2)^4}.
\end{equation*}
Therefore
\begin{equation*}
   u = 4\left(\frac{\theta_{01}(\tau/2)^4}{\theta_{00}(\tau/2)^4}
     - 1\right)^{-1} + 2
   = 
   \frac{- 4\theta_{00}(\tau/2)^4 + 2\theta_{10}(\tau/2)^4}
   {\theta_{10}(\tau/2)^4},
\end{equation*}
where we have used \cite[(3.35)]{Ume} for $\tau/2$. The denominator is
\begin{equation*}
   {4\theta_{00}^2(\tau)\theta_{10}^2(\tau)}
\end{equation*}
by {\cite[(3.79)]{Ume}}. The numerator is
\begin{equation*}
\begin{split}
   & - 4\theta_{00}(\tau/2)^4 + 2\theta_{10}(\tau/2)^4
   \overset{\text{\cite[(3.35)]{Ume}}}{=}
   -2\theta_{00}(\tau/2)^4 - 2\theta_{01}(\tau/2)^4
\\
   = \; &
   -2 \left(\theta_{00}(\tau/2)^2 + \theta_{01}(\tau/2)^2\right)^2
   + 4 \theta_{00}(\tau/2)^2\theta_{01}(\tau/2)^2
   \overset{\text{\cite[(3.90),(3.78)]{Ume}}}{=}
   -8 \theta_{00}(\tau)^4 + 4\theta_{01}(\tau)^4
\\
   \overset{\text{\cite[(3.35)]{Ume}}}{=} &
   - 4\theta_{00}(\tau)^4 - 4\theta_{10}(\tau)^4.
\end{split}
\end{equation*}
Therefore we get the first assertion.

The second assertion follows as
\begin{equation*}
\begin{split}
   &\frac{da}{du} = -\frac1{2\pi}\int_{A} \frac{dz}{\sqrt{P(z)^2-4}}
   = -\frac1{\pi}\; \frac{\sqrt{-2-u}}{\sqrt{u^2-4}}
   \int_{-1}^{-1/\kappa}
   \frac{dx}{\sqrt{(1-x^2)(1-\kappa^2 x^2)}}
\\
   {\overset{\text{\cite[Lem.~4.1]{Ume}}}{=}}
   \; &
   -\frac1{\pi}\; \frac1{\sqrt{2-u}}
   \sqrt{-1} K'
   = -\frac1{2\pi}\; \frac{\theta_{10}(\tau/2)^2}{\theta_{00}(\tau/2)^2}
   \sqrt{-1}K(\kappa'(-2/\tau))
\\
   =
   \; &
   - \frac1{2\pi}\; \frac{\theta_{10}(\tau/2)^2}{\theta_{00}(\tau/2)^2}
   \sqrt{-1}K(\kappa(\tau/2))
   {\overset{\text{\cite[Th.~4.1]{Ume}}}{=}}
   -\frac{\sqrt{-1}}4 \theta_{10}(\tau/2)^2
\\
   \overset{\text{\cite[(3.79)]{Ume}}}{=}
   \; &
   -\frac{\sqrt{-1}}2 \theta_{00}(\tau)\theta_{10}(\tau).
\end{split}
\end{equation*}
Since the branch of the square root must be chosen carefully, we fix
the sign later. (See p.195 for \cite{Ume}.)

Let us show the third assertion. We have
\begin{equation*}
\begin{split}
  a &= -\frac1{2\pi} \int_A \frac{z P'(z) dz}{\sqrt{P(z)^2 - 4}}
  = -\frac2{\pi} \frac{(-2-u)^{3/2}}{\sqrt{u^2-4}}
  \int_{-1}^{-1/\kappa} \frac{x^2 dx}{\sqrt{(1-x^2)(1-\kappa^2 x^2)}}
\\
  &\overset{x = 1/\kappa w}{=}
  -\frac2{\pi}\frac{-2-u}{\sqrt{2-u}} \frac1{\kappa^2}
  \int_{1}^{1/\kappa} \frac{w^{-2}dw}{\sqrt{(1-w^2)(1-\kappa^2 w^2)}}
\\
  &\underset{\text{See the proof of \cite[Lem. 4.1]{Ume}}}
  {\overset{w=(1-\kappa^{\prime2}t^2)^{-1/2}}{=}}
  \frac2{\pi} \sqrt{-1} \sqrt{2-u}
  \int_{0}^{1} \frac{(1-\kappa^{\prime2}t^2)dt}
  {\sqrt{(1-t^2)(1-\kappa^{\prime2}t^2)}} 
\\  
   &= \sqrt{-1}\frac4\pi \frac{E(\kappa')}{\kappa'}
   = \sqrt{-1}\frac4\pi \frac{E(\kappa'(-2/\tau))}{\kappa'(-2/\tau)}
   = \sqrt{-1}\frac4\pi \frac{E(\kappa(\tau/2))}{\kappa(\tau/2)}
\\
   & \overset{\text{\cite[\S22.73]{WW}}}{=}
   \sqrt{-1}
   \frac4\pi \frac{\theta_{00}(\tau/2)^2}{\theta_{10}(\tau/2)^2}
   \left(
    \frac{\theta_{00}(\tau/2)^4 + \theta_{01}(\tau/2)^4}
    {3\theta_{00}(\tau/2)^4} \frac\pi2 \theta_{00}(\tau/2)^2
    - \frac1{6\pi} \frac1{\theta_{00}(\tau/2)^{2}}
    \frac{\theta_{11}'''(\tau/2)}{\theta_{11}'(\tau/2)}
    \right).
\end{split}
\end{equation*}
Again we check the sign later.
We further use
\begin{align*}
   & -\theta_{11}'''(\tau)/\theta_{11}'(\tau) = \pi^2 E_2(\tau),
   \tag*{\cite[\S21.41]{WW}}
\\
   & 2E_2(\tau)-E_2(\tau/2)=
   \frac{\theta_{00}^4(\tau/2)+\theta_{01}^4(\tau/2)}2
   = \theta_{00}^4(\tau)+\theta_{10}^4(\tau).
   \tag*{??}
\end{align*}
Then
\begin{equation*}
\begin{split}
  a &= \sqrt{-1}
    \frac{2(\theta_{00}(\tau/2)^4 + \theta_{01}(\tau/2)^4 + E_2(\tau/2))}
    {3\theta_{10}(\tau/2)^2}
\\
  &= \sqrt{-1}\frac{2(\theta_{00}(\tau)^4 + \theta_{10}(\tau)^4) +
     2E_2(\tau)-\theta_{00}^4(\tau)-\theta_{10}^4(\tau)}
       {3\theta_{00}(\tau)\theta_{10}(\tau)}
    = \sqrt{-1}
    \frac{2E_2(\tau)+\theta_{00}^4(\tau)+\theta_{10}^4(\tau)}
       {3\theta_{00}(\tau)\theta_{10}(\tau)}.
\end{split}
\end{equation*}

In order to determine the sign, we suppose $P(z)$ satisfies
assumptions as in the proof of \propref{prop:F_beh}. Then we study the
asymptotic behavior at $u\to -\infty$.
We have $a = a_2 \sim -\sqrt{-1} z_2 \sim \sqrt{-1} \sqrt{2-u}$. We
also have $\kappa\to 1$, $\kappa'\to 0$, so $E(\kappa')\to \pi/2$.
Therefore the sign of $a$ is correct. (Or, just check that
$a\in\sqrt{-1}\R_{>0}$.)

We also have
\(
  \frac{da}{du} \sim -\frac{\sqrt{-1}}{2\sqrt{2-u}}.
\)
Since $K(\kappa')\to\pi/2$, the sign of $da/du$ is also correct.
\end{proof}
\end{NB}

Let us record the following formula for the later purpose.
\begin{equation}\label{eq:contact_term}
\begin{split}
   & \Lambda\pd{u}{\Lambda}
   = 2u - a \frac{du}{da}
   = 2u + \sqrt{-1}\frac
  {2E_2+\theta_{00}^4+\theta_{10}^4}{3\theta_{00}\theta_{10}}\Lambda
  \frac{du}{da}
\\
  = \; &
  2u -\frac13 E_2 \left(\frac{du}{da}\right)^2 
  + \frac23 \frac{\theta_{00}^4+\theta_{10}^4}
  {\theta_{00}^2\theta_{10}^2}\Lambda^2
  =  -\frac13 E_2 \left(\frac{du}{da}\right)^2 + \frac43 u,
\end{split}
\end{equation}
where the first equality follows from the homogeneity of $u$.

For the reader who wants to compare the formulas with ones in the
literature, we record how parameters differ.
Let us make a change of variable by
\begin{equation*}
   w - \frac{u}3 = -\Lambda^2
   \frac{z-\sqrt{-u+2\Lambda^2}}{z+\sqrt{-u+2\Lambda^2}}.
\end{equation*}
The branched points
$z=\sqrt{-u+2\Lambda^2}$, $-\sqrt{-u+2\Lambda^2}$,
$\sqrt{-u-2\Lambda^2}$, $-\sqrt{-u-2\Lambda^2}$ are mapped to
$w=u/3$, $\infty$, $(-u-3\sqrt{u^2-4\Lambda^4})/6$,
$(-u+3\sqrt{u^2-4\Lambda^4})/6$ respectively. And the curve
\(
   y^2=(z^2+u)^2-4\Lambda^4
\)
is isomorphic to a Weierstrass form
\(
   y^2=4w^3-g_2w-g_3
\)
with
\begin{equation*}
   g_2 = 4 \left(\frac{1}{3}u^2-\Lambda^4\right), \qquad
   g_3 = -\frac{1}{27} u\left(8u^2-36\Lambda^4\right).
\end{equation*}
This is the form of curves appeared in \cite{FS} with $u$ replaced by
$-x$ therein.
If we make a further change of variable as
\(
   w = 2x + u/3,
\)
we get
\begin{equation*}
   \frac1{32} y^2
   = x\left(x^2 + \frac{u}2 x + \frac{\Lambda^4}4\right).
\end{equation*}
\begin{NB}
\begin{equation*}
\begin{split}
  & 4 w^3-g_2w-g_3
  = 4 (2x + \frac{u}3)^3 - 4 \left(\frac{1}{3}u^2-\Lambda^4\right)
  (2x + \frac{u}3) + \frac{1}{27} u\left(8u^2-36\Lambda^4\right)
\\
  = \; &
  4\cdot 8 x^3  + 4\cdot 4 u x^2 + \frac83 u^2 x + \frac4{27} u^3
  - \frac83 u^2 x - \frac4{9}u^3 + 8\Lambda^4 x + \frac43 u\Lambda^4
  + \frac{8}{27}u^3 - \frac43 u \Lambda^4
\\
  = \; &
  32 x ( x^2 + \frac{x}2 v + \frac{\Lambda^4}4).
\end{split}
\end{equation*}
\end{NB}
This is the form of the curves in \cite{MW} after the replacement 
$u\mapsto -2u$, $y \mapsto 4\sqrt{2}y$.

\section{Instanton moduli spaces}

\subsection{Basic definitions}
In this and next subsections we briefly recall properties of framed
moduli spaces of instantons (resp.\ torsion-free sheaves) on $S^4$
(resp.\ $\proj^2$) and the corresponding moduli spaces on blowup. For
more detail, see \cite[\S1]{part1} and \cite[Chapters~2,3]{Lecture}
and the references therein.

Let $M(r,n)$ be the framed moduli space of torsion free sheaves on
$\proj^2$ with rank $r$ and $c_2 = n$, which parametrizes isomorphism
classes of $(E,\Phi)$ such that 
\begin{enumerate}
\item $E$ is a torsion free sheaf of $\rank E =r$, $\langle c_2(E),
[\proj^2]\rangle =n$ which is locally free in a neighbourhood of
$\linf$,
\item $\Phi \colon E|_{\linf} {\overset{\sim}{\to}}
\shfO_{\linf}^{\oplus r}$ is an isomorphism called `framing at infinity'.
\end{enumerate}
Here $\linf=\{[0:z_1:z_2] \in \proj^2 \} \subset \proj^2$ is the line
at infinity.
Notice that the existence of a framing $\Phi$ implies $c_1(E)=0$.

This is known to be nonsingular of dimension $2nr$.

Let $\Mreg(r,n)$ be the open subset consisting of locally free
sheaves. 
By a result of Donaldson~\cite{Don} it can be identified with the
framed moduli space of instantons on $S^4$ which parametrizes
anti-self-dual connections $A$ on a principal $\SU(r)$-bundle $P$ with
$\langle c_2(P),[S^4]\rangle = n$ modulo gauge transformations $\gamma$
with $\gamma_\infty=\id$.

Let $M_0(r,n)$ be the Uhlenbeck (partial) compactification of
$\Mreg(r,n)$. Set theoretically it is defined by
\begin{equation*}
     M_0(r,n) = 
   \bigsqcup_{k=0}^n M_0^{\operatorname{reg}}(r,n-k)\times S^{k}\C^2,
\end{equation*}
where $S^k\C^2$ is the $k$th symmetric product of $\C^2$. We can endow
a structure of an affine algebraic variety to $M_0(r,n)$ so that there
is a projective morphism
\begin{equation*}
    \pi\colon M(r,n)\to M_0(r,n).
\end{equation*}
The corresponding map between closed points can be identified with
\begin{equation*}\label{eq:map_pi}
   (E,\Phi) \longmapsto
   ((E^{\vee\vee},\Phi), \operatorname{Supp}(E^{\vee\vee}/E))\in
   M_0^{\operatorname{reg}}(r,n')\times S^{n-n'}\C^2.
\end{equation*}
where $E^{\vee\vee}$ is the double dual of $E$ and
$\operatorname{Supp}(E^{\vee\vee}/E)$ is the support of
$E^{\vee\vee}/E$ counted with multiplicities. Note that $E^{\vee\vee}$
is a locally free sheaf. For moduli spaces on general projective
surfaces, such morphisms from moduli spaces of sheaves to Uhlenbeck
compactifications were constructed by J.~Li \cite{Li} and Morgan
\cite{Mor}.

Let $T$ be the maximal torus of $\GL_r(\C)$ consisting of diagonal
matrices and let $\hT = \C^*\times\C^* \times T$. We define an action
of $\hT$ on $M(r,n)$ as follows: For 
$(t_1,t_2)\in \C^*\times\C^*$, let $F_{t_1,t_2}$ be an automorphism of 
$\proj^2$ defined by
\[
    F_{t_1,t_2}([z_0: z_1 : z_2]) = [z_0: t_1 z_1 : t_2 z_2].
\]
For $\operatorname{diag}(e_1,\dots,e_r)\in T$ let $G_{e_1,\dots,e_r}$
denotes the isomorphism of $\shfO_{\linf}^{\oplus r}$ given by
\[
    \shfO_{\linf}^{\oplus r}\ni (s_1,\dots, s_r) \longmapsto
     (e_1 s_1, \dots, e_r s_r).
\] 
Then for $(E,\Phi)\in M(r,n)$, we define
\begin{equation}\label{eq:action}
    (t_1,t_2,e_1,\dots,e_r)\cdot (E,\Phi)
    = \left((F_{t_1,t_2}^{-1})^* E, \Phi'\right),
\end{equation}
where $\Phi'$ is the composite of homomorphisms
\begin{equation*}
   (F_{t_1,t_2}^{-1})^* E|_{\linf} 
   \xrightarrow{(F_{t_1,t_2}^{-1})^*\Phi}
   (F_{t_1,t_2}^{-1})^* \shfO_{\linf}^{\oplus r}
   \longrightarrow \shfO_{\linf}^{\oplus r}
   \xrightarrow{G_{e_1,\dots, e_r}} \shfO_{\linf}^{\oplus r}.
\end{equation*}
Here the middle arrow is the homomorphism given by the action.

In a similar way, we have a $\hT$-action on $M_0(r,n)$. The map
$\pi\colon M(r,n)\to M_0(r,n)$ is equivariant.

The fixed points $M(r,n)^{\hT}$ consist of $(E,\Phi) =
(I_1,\Phi_1)\oplus\cdots\oplus (I_r,\Phi_r)$ such that
\begin{aenume}
\item $I_\alpha$ is an ideal sheaf of $0$-dimensional subscheme
$Z_\alpha$ contained in $\C^2 = \proj^2\setminus\linf$.
\item $\Phi_\alpha$ is an isomorphism from $(I_\alpha)_{\linf}$ to the $\alpha$th factor of $\shfO_{\linf}^{\oplus r}$.
\item $I_\alpha$ is fixed by the action of $\C^*\times\C^*$, coming
  from that on $\proj^2$.
\end{aenume}

On the other hand, the fixed points $M_0(r,n)^{\hT}$ consist of the single
point $n[0]\in S^n\C^2 \subset M_0(r,n)$.

\begin{figure}[htbp]
\centering
\psset{xunit=.8mm,yunit=.8mm,runit=.8mm}
\begin{pspicture}(0,0)(65.00,65.00)
\psline[linewidth=0.40,linecolor=black]{<-}(5.00,65.00)(5.00,5.00)
\psline[linewidth=0.40,linecolor=black]{->}(5.00,5.00)(65.00,5.00)
\psframe[linewidth=0.15,linecolor=black](10.00,10.00)(20.00,20.00)
\psframe[linewidth=0.15,linecolor=black](20.00,10.00)(30.00,20.00)
\psframe[linewidth=0.15,linecolor=black](30.00,10.00)(40.00,20.00)
\psframe[linewidth=0.15,linecolor=black,fillcolor=black,fillstyle=hlines,hatchwidth=0.28,hatchsep=1.42,hatchangle=45.00,hatchcolor=black](30.00,20.00)(40.00,30.00)
\psframe[linewidth=0.15,linecolor=black](10.00,50.00)(20.00,60.00)
\psframe[linewidth=0.15,linecolor=black](10.00,20.00)(20.00,30.00)
\psframe[linewidth=0.15,linecolor=black,fillcolor=lightgray,fillstyle=solid](20.00,20.00)(30.00,30.00)
\psframe[linewidth=0.15,linecolor=black](40.00,10.00)(50.00,20.00)
\psframe[linewidth=0.15,linecolor=black](10.00,30.00)(20.00,40.00)
\psframe[linewidth=0.15,linecolor=black](20.00,20.00)(30.00,30.00)
\psframe[linewidth=0.15,linecolor=black](50.00,10.00)(60.00,20.00)
\psframe[linewidth=0.15,linecolor=black,fillcolor=black,fillstyle=vlines,hatchwidth=0.28,hatchsep=1.42,hatchangle=45.00,hatchcolor=black](20.00,30.00)(30.00,40.00)
\psframe[linewidth=0.15,linecolor=black](10.00,40.00)(20.00,50.00)
\psframe[linewidth=0.15,linecolor=black,fillcolor=black,fillstyle=vlines,hatchwidth=0.28,hatchsep=1.42,hatchangle=45.00,hatchcolor=black](20.00,40.00)(30.00,50.00)
\psframe[linewidth=0.15,linecolor=black](30.00,30.00)(40.00,40.00)
\rput(65.00,15.17){$x^5$}
\rput(64.48,14.48){}
\rput(45.00,25.00){$x^3y$}
\rput(25.00,55.00){$xy^4$}
\rput(15.00,65.00){$y^5$}
\psline[linewidth=0.15,linecolor=black]{<-}(25.17,45.86)(40.17,55.86)
\rput(42.41,52.33){}
\rput(46.17,57.64){$a_Y(s)$}
\rput(42.16,52.33){}
\psline[linewidth=0.15,linecolor=black]{<-}(35.00,25.00)(50.00,35.00)
\rput(55.64,36.86){$l_Y(s)$}
\rput(26.72,24.57){}
\rput(25.00,25.00){$s$}
\rput(35.00,45.00){$x^2y^3$}
\end{pspicture}
\caption{Young diagram and ideal}
\label{fig:young}
\end{figure}
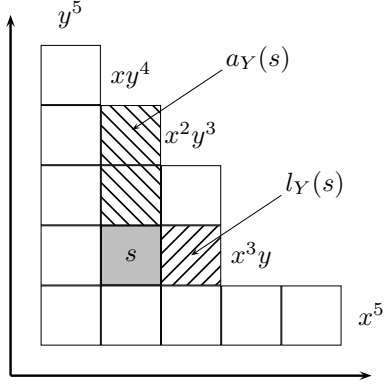

We parametrize the fixed point set $M(r,n)^{\hT}$ by a $r$-tuple of
Young diagrams $\vec{Y} = (Y_1,\dots,Y_r)$ so that the ideal
$I_\alpha$ is spanned by monomials $x^i y^j$ placed at $(i-1,j-1)$
outside $Y_\alpha$ as illustrated in Figure~\ref{fig:young}. The
constraint is that the total number of boxes
\(
    |\vec{Y}| \defeq \sum_\alpha |Y_\alpha|
\)
is equal to $n$.

Let $Y = (\lambda_1\ge \lambda_2 \ge \cdots)$ be a Young diagram,
where $\lambda_i$ is the length of the $i$th column.
Let $Y' = (\lambda'_1\ge \lambda_2' \ge \dots)$ be the transpose of
$Y$. Thus $\lambda'_j$ is the length of the $j$th row of $Y$.
Let $l(Y)$ denote the number of columns of $Y$, i.e., $l(Y) =
\lambda'_1$.
Let
\begin{alignat*}{2}
& a_Y(i,j) = \lambda_i - j, & \qquad & a'(i,j) = j - 1 \\
& l_Y(i,j) = \lambda'_j - i, &\qquad & l'(i,j) = i - 1.
\end{alignat*}
Here we set $\lambda_i = 0$ when $i > l(Y)$. Similarly $\lambda'_j =
0$ when $j > l(Y')$.
When the square $s = (i,j)$ lies in $Y$, these are called {\it
arm-length}, {\it arm-colength}, {\it leg-length}, {\it leg-colength}
respectively, and we usually consider in this case. But our formula
below involves these also for squares outside $Y$. So these
take negative values in general. Note that $a'$ and $l'$ does not
depend on the diagram, and we do not write the subscript $Y$.

\begin{Theorem}\label{thm:M(r,n)weights}
Let $(E,\Phi)$ be a fixed point of $\hT$-action corresponding to
$\vec{Y} = (Y_1,\dots, Y_r)$. Then the $\hT$-module structure of
$T_{(E,\Phi)} M(r,n)$ is given by
\begin{equation*}
   \sum_{\alpha,\beta=1}^r N_{\alpha,\beta}^{\vec{Y}}(t_1,t_2),
\end{equation*}
where
\[
   N_{\alpha,\beta}^{\vec{Y}}(t_1,t_2)
   = e_\beta\, e_\alpha^{-1}\times
   \left\{
   \sum_{s \in Y_\alpha}
      \left( t_1^{-l_{Y_\beta}(s)} t_2^{a_{Y_\alpha}(s)+1}\right)
%
    + \sum_{t\in Y_\beta} 
        \left(t_1^{l_{Y_\alpha}(t)+1}t_2^{-a_{Y_\beta}(t)}\right)
        \right\}
.
\]
\end{Theorem}

Here we have used the following notation.
\begin{Notation}\label{not:module}
We denote by $e_\alpha$ ($\alpha=1,\dots, r$) the one dimensional
$\hT$-module given by
\begin{equation*}
   \hT\ni (t_1,t_2, e_1, \dots, e_r) \mapsto e_\alpha.
\end{equation*}
Similarly, $t_1$, $t_2$ denote one-dimensional $\hT$-modules. Thus
the representation ring $R(\hT)$ is isomorphic to $\Z[t_1^\pm,
t_2^\pm, e_1^\pm, \dots, e_r^\pm]$, where $e_\alpha^{-1}$ is the dual
of $e_\alpha$.
\end{Notation}

\subsection{Moduli spaces on the blowup}

Let $\bp$ be the blowup of $\proj^2$ at $[1:0:0]$. Let $p\colon
\bp\to\proj^2$ denote the projection. The manifold $\bp$ is the closed
subvariety of $\proj^2\times\proj^1$ defined by
\begin{equation*}
   \{ ([z_0 : z_1 : z_2], [z : w] \in \proj^2\times\proj^1\mid 
   z_1 w = z_2 z \},
\end{equation*}
where the map $p\colon \bp\to\proj^2$ is the projection to the first
factor. Let us denote the inverse image of $\linf$ under
$\bp\to\proj^2$ also by $\linf$ for brevity. It is given by the
equation $z_0 = 0$. The complement $\bp\setminus\linf$ is the blowup
${\widehat\C}^2$ of $\C^2$ at the origin.
Let $C$ denote the exceptional set. It is given by
$z_1 = z_2 = 0$.

In this subsection, $\shfO$ denotes the structure sheaf of $\bp$,
$\shfO(C)$ the line bundle associated with the divisor $C$,
$\shfO(mC)$ its $m$th tensor product.

Let $\bM(r,k,n)$ be the framed moduli space of torsion free sheaves
$(E,\Phi)$ on $\bp$ with rank $r$, $\langle c_1(E),[C]\rangle = -k$
and $\langle c_2(E) - \frac{r-1}{2r} c_1(E)^2, [\bp]\rangle = n$.
This is also nonsingular of dimension $2nr$. (Remark that $n$ may not
be integer in general.)

\begin{Theorem}\label{thm:blowup}
There is a projective morphism $\widehat\pi\colon\bM(r,k,n)\to
M_0(r,n-\frac1{2r}k(r-k))$ \textup($0\le k<r$\textup) defined by
\begin{equation*}
   (E,\Phi) \mapsto \left(((p_* E)^{\vee\vee}, \Phi),
   \Supp(p_*E^{\vee\vee}/p_*E) + \Supp(R^1p_* E)\right).
\end{equation*}
\end{Theorem}

Let us define an action of the $(r+2)$-dimensional torus $\hT =
\C^*\times\C^* \times T$ on $\bM(r,k,n)$ by modifying the action on
$M(r,n)$ as follows. For $(t_1,t_2)\in \C^*\times\C^*$, let
$F'_{t_1,t_2}$ be an automorphism of $\bp$ defined by
\[
    F'_{t_1,t_2}([z_0: z_1 : z_2], [z : w]) 
    = ([z_0: t_1 z_1 : t_2 z_2], [t_1 z : t_2 w]).
\]
Then we define the action by replacing $F_{t_1,t_2}$ by $F'_{t_1,t_2}$
in \eqref{eq:action}. The action of the latter $T$ is exactly the same 
as before. The morphism $\widehat\pi$ is equivariant.

Note that the fixed point set of $\C^*\times\C^*$ in ${\widehat\C}^2 =
\bp \setminus\linf$ consists of two points $([1 : 0 : 0], [1 : 0])$,
$([1 : 0 : 0], [0 : 1])$. Let us denote them $p_1$ and $p_2$.

Let us define an action of the $(r+2)$-dimensional torus $\hT =
\C^*\times\C^* \times T$ on $\bM(r,k,n)$ by modifying the action on
$M(r,n)$ as follows. For $(t_1,t_2)\in \C^*\times\C^*$, let
$F'_{t_1,t_2}$ be an automorphism of $\bp$ defined by
\[
    F'_{t_1,t_2}([z_0: z_1 : z_2], [z : w]) 
    = ([z_0: t_1 z_1 : t_2 z_2], [t_1 z : t_2 w]).
\]
Then we define the action by replacing $F_{t_1,t_2}$ by $F'_{t_1,t_2}$
in \eqref{eq:action}. The action of the latter $T$ is exactly the same 
as before. The morphism $\widehat\pi$ is equivariant.

The fixed points $\bM(r,k,n)^{\hT}$ consist of
\( 
   (E,\Phi) = (I_1(k_1 C),\Phi_1) \oplus \cdots \oplus
   (I_r(k_r C),\Phi_r)
\)
such that
\begin{aenume}
\item $I_\alpha(k_\alpha C)$ is the tensor product
$I_\alpha\otimes\shfO(k_\alpha C)$, where $k_\alpha\in\Z$ and
$I_\alpha$ is an ideal sheaf of $0$-dimensional subscheme $Z_\alpha$
contained in ${\widehat\C}^2 = \bp\setminus\linf$.
\item $\Phi_\alpha$ is an isomorphism from $(I_\alpha)_{\linf}$ to the $\alpha$th factor of $\shfO_{\linf}^{\oplus r}$.
\item $I_\alpha$ is fixed by the action of $\C^*\times\C^*$, coming
from that on $\bp$.
\end{aenume}

\begin{figure}[htbp]
\centering
\psset{xunit=.8mm,yunit=.8mm,runit=.8mm}
\begin{pspicture}(0,0)(60.00,60.00)
\psarc[linewidth=0.50,linecolor=black]{-}(5.00,5.00){25.00}{-10.00}{100.00}
\psline[linewidth=0.50,linecolor=black]{->}(25.00,5.00)(60.00,5.00)
\psline[linewidth=0.50,linecolor=black]{->}(5.00,25.00)(5.00,60.00)
\rput(25.00,25.00){$C$}
\rput(45.00,10.00){$(z_1/z_0,w/z)$}
\rput(20.00,35.00){$(z/w,z_2/z_0)$}
\rput(1.38,26.38){$p_2$}
\rput(25.52,2.07){$p_1$}
\end{pspicture}
\caption{blowup and fixed points}
\label{fig:blowup}
\end{figure}
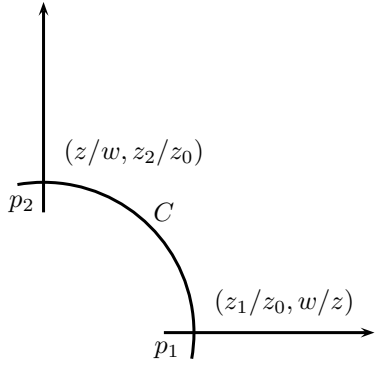

The support of $Z_\alpha$ must be contained in
the fixed point set in ${\widehat\C}^2$, i.e., $\{ p_1, p_2 \}$.
Thus $Z_\alpha$ is a union of $Z_\alpha^1$ and $Z_\alpha^2$,
subschemes supported at $p_1$ and $p_2$ respectively.
If we take a coordinate system $(x, y) = (z_1/z_0, w/z)$ (resp.\ $=
(z/w, z_2/z_0)$) around $p_1$ (resp.\ $p_2$), then $Z_\alpha^1$
(resp.\ $Z_\alpha^2$) is generated by monomials $x^i y^j$. (See
Figure~\ref{fig:blowup}.) Then $Z_\alpha^1$ (resp.\ $Z_\alpha^2$)
corresponds to a Young diagram $Y_\alpha^1$ (resp.\ $Y_\alpha^2$) as
before. Therefore the fixed point set is parametrized by $r$-tuples
\(
   (\vec{k},\vec{Y^1},\vec{Y^2}) =
   ((k_1,Y_1^1, Y_1^2),\linebreak[2]\dots,\linebreak[1]
   (k_r, Y_r^1, Y_r^2)),
\)
where $k_\alpha\in\Z$ and $Y_\alpha^1$, $Y_\alpha^2$ are
Young diagrams. The constraint is
\begin{equation}\label{eq:cond}
    \sum_\alpha k_\alpha = k, \qquad
    |\vec{Y}^1| + |\vec{Y}^2| + 
    \frac1{2r} \sum_{\alpha < \beta} |k_\alpha - k_\beta|^2 = n.
\end{equation}
We will use the convention for $\vec{k}$ in \secref{sec:root}.

Note that the fixed point data have three parts, $\vec{k}$,
$\vec{Y^1}$ and $\vec{Y^2}$. This will be reflected in the blowup
formula \eqref{eq:blowup_formula} below. Also, the appearance of
$\vec{k}\in \Z^r$ explain the reason why many formulas below contain
the theta function.

\begin{Theorem}\label{thm:bM(r,k,n)weights}
Let $(E,\Phi)$ be a fixed point of $\hT$-action corresponding to
$(\vec{k}, \vec{Y}^1,\linebreak[0]\vec{Y}^2)$. Then the $\hT$-module
structure of $T_{(E,\Phi)} \bM(r,k,n)$ is given by
\begin{equation*}
   \!\!\sum_{\alpha,\beta=1}^r\!\! L_{\alpha,\beta}^{\vec{k}}(t_1,t_2) 
     + t_1^{k_\beta - k_\alpha}
       N^{\vec{Y^1}}_{\alpha,\beta}(t_1,t_2/t_1)
     + t_2^{k_\beta - k_\alpha}
       N^{\vec{Y^2}}_{\alpha,\beta}(t_1/t_2,t_2),
\end{equation*}
where
\begin{equation*}
   L_{\alpha,\beta}^{\vec{k}}(t_1,t_2) = 
   e_\beta\, e_\alpha^{-1}\times
   \begin{cases}
     {\displaystyle
     \sum_{\substack{i,j\ge 0\\i+j \le k_\alpha-k_\beta-1}}
          t_1^{-i} t_2^{-j}}
       & \text{if $k_\alpha > k_\beta$}, \\
     {\displaystyle
     \sum_{\substack{i,j\ge 0\\i+j \le k_\beta-k_\alpha-2}}
          t_1^{i+1} t_2^{j+1}}
       & \text{if $k_\alpha + 1 < k_\beta$}, \\
     0 & \text{otherwise}.
   \end{cases}
\end{equation*}
\end{Theorem}

The reason for the change of weights $(t_1,t_2/t_1)$, $(t_1/t_2,t_2)$
is clear. It comes from the action on the coordinate system around
$p_1$ and $p_2$.

\subsection{Topology of moduli spaces}
Thanks to the existence of the torus action, the homology groups of
$M(r,n)$, $\bM(r,k,n)$ enjoy nice properties. In particular, we can
calculate their Betti numbers whose generating functions have
beautiful formulas. The results of this and next subsections will not
be used in the other parts of this paper. The reader in a hurry may
skip this and next subsections.

\begin{Theorem}\label{thm:topology}
\rom{(1)} $\pi^{-1}(n[0])$ is isomorphic to the punctual quot-scheme
pa\-rameterizing zero dimensional quotients $\shfO_{\proj^2}^{\oplus r}
\to Q$ with $\operatorname{Supp}(Q) = n[0]$.

\rom{(2)} $M(r,n)$ is homotopy equivalent to $\pi^{-1}(n[0])$.

\rom{(3)} Both $M(r,n)$ and $\pi^{-1}(n[0])$ have $\alpha$-partitions
into affine spaces.

\rom{(4)} $H_{\mathop{\mathrm{odd}}}(M(r,n),\Z) = 0$ and
$H_{\mathop{\mathrm{even}}}(M(r,n),\Z)$ is a free abelian group. The
cycle map $A_*(M(r,n))\to H_{\mathop{\mathrm{even}}}(M(r,n),\Z)$ is an
isomorphism. The same assertions hold for $\pi^{-1}(n[0])$.
\end{Theorem}

Recall that a finite partition of a variety $X$ into locally closed
subvarieties is said to be an $\alpha$-partition if the subvarities in 
the partition can be indexed $X_1,\dots, X_n$ in such a way that
$X_1\cup X_2\cup\dots\cup X_i$ is closed in $X$ for $i=1,\dots,n$.

Here $H_*(\bullet,\Z)$, $A_*(\bullet)$ denote the Borel-Moore homology
group and the Chow group. See \secref{sec:BMhom}.

\begin{proof}
(1) By the geometric description \eqref{eq:map_pi} of the map $\pi$,
the fiber $\pi^{-1}(n[0])$ consists of $(E,\Phi)$ such that
$E^{\vee\vee} = \shfO_{\proj^2}^{\oplus r}$ and
$\operatorname{Supp}(\shfO_{\proj^2}^{\oplus r}/E) = n[0]$. Thus the
quotient $\shfO_{\proj^2}^{\oplus r}/E$ is a point in the punctual
quot-scheme.

(2) A similar result was proved in \cite[5.5]{Na-Duke1} for quiver
varieties by using a $\C^*$-action. The proof can be adapted to our
situation as follows.

Let us consider a one parameter subgroup $\C^*\ni t\mapsto \lambda(t)
= (t,t,1,\dots,1)\in \hT$. When $t\to 0$, $\lambda(t)\cdot x$ goes
to $0$ for any $x$ in $M_0(r,n)$. Now apply an argument of
Slodowy~\cite[4.3]{Slo-Four} to $\pi\colon M(r,n)\to M_0(r,n)$.

(3) Choose a generic one parameter subgroup $\lambda\colon \C^*\to
\hT$ so that the fixed point set is unchanged: $M(r,n)^{\lambda(\C^*)}
= M(r,n)^{\hT}$. Moreover we take so that $\lim_{t\to 0}\lambda(t) =
0$.
Then points fixed by $\lambda(\C^*)$ are given as above. In
particular, they form a finite set. For each fixed point $w$, we
consider $(\pm)$-attracting set:
\begin{equation*}
   S_w = \left\{ x\in M(r,n) \left|\;
       \lim_{t\to 0} \lambda(t)\cdot x = w \right\}\right., \quad
   U_w = \left\{ x\in M(r,n) \left|\;
       \lim_{t\to\infty}\lambda(t)\cdot x = w \right\}\right..
\end{equation*}
These are affine spaces by \cite{BB}. Moreover, there exists an order
on fixed points so that $\bigcup_{y\le w} S_y$ (resp.\ $\bigcup_{y\le
w} U_y$ is closed in $\bigcup_x S_x$ (resp.\ $\bigcup_x U_x$) for each 
$w$. (See e.g., \cite[\S1]{AB1}.)
We claim that $\bigcup S_w = M(r,n)$, $\bigcup U_w = \pi^{-1}(n[0])$.
Consider the corresponding action on $M_0(r,n)$. 
For this purpose, we recall the ADHM description: $M_0(r,n)$ is an
affine algebro-geometric quotient
\begin{equation*}
   \{ (B_1,B_2,i,j) \mid [B_1,B_2] + ij=0 \}\dslash \GL_n(\C),
\end{equation*}
where $B_1$, $B_2$ are $n\times n$ complex matrices, and
$i\in\Hom(\C^r,\C^n)$, $j\in\Hom(\C^n,\C^r)$. (See
\cite[Chapter~2]{Lecture} and \cite{part1}.)
The action of $\GL_n(\C)$ is given by
\(
 g \cdot (B_1,B_2,i,j)=( g B_1 g^{-1}, g B_2 g^{-1}, gi, jg^{-1}).
\)
The $\hT$-action is given by
\begin{multline*}
    (B_1, B_2, i, j) \longmapsto
     (t_1 B_1, t_2 B_2, i e^{-1}, t_1 t_2 e j), \qquad
\\
   \text{for $t_1, t_2\in \C^*$,
             $e = \operatorname{diag}(e_1,\dots, e_r)\in(\C^*)^r$}.
\end{multline*}

By \cite{Lu-qv} the coordinate ring of $M_0(r,n)$ is generated by the
following two types of functions:
\begin{aenume}
\item $\tr(B_{\alpha_N}B_{\alpha_{N-1}}\cdots B_{\alpha_1}\colon
\C^n \to \C^n)$,
where $\alpha_i = 1$ or $2$.
\item $\chi(jB_{\alpha_N} B_{\alpha_{N-1}} \cdots B_{\alpha_1} i)$,
where $\alpha_i = 1$ or $2$, and $\chi$ is a linear form on
$\End(\C^r)$.
\end{aenume}
Both types of functions have positive weights with respect to
$\lambda$. Therefore every point in $M_0(r,n)$ converges to $0$ as
$t\to 0$, and every point except $0$ goes to infinity as $t\to \infty$.
Since $\pi$ is proper, we get the claim. (We use the fact that the
orbit of a one-parameter subgroup has limit if it is contained in a
compact set.) Thus $\bigcup S_w$ (resp.\ $\bigcup U_w$) gives us an
$\alpha$-partition of $M(r,n)$ (resp.\ $\pi^{-1}(0)$) into affine
spaces.

(4) is a consequence of (3) by \cite[Lemma 1.8]{DLP}.
\end{proof}

\begin{Theorem}\label{thm:P_t(M(r,n))}
The Poincar\'e polynomials of the punctual quot-scheme
$\pi^{-1}(n[0])$ is given by
\begin{equation}\label{eq:P_t(M(r,n))}
    P_t(\pi^{-1}(n[0])) 
    = \sum_{(Y_1,\dots,Y_r)} \prod_{\alpha=1}^r
           t^{2(r|Y_\alpha| - \alpha l(Y_\alpha))},
\end{equation}
where the summation runs over the set of $r$-tuple of Young diagrams
$\vec{Y} = (Y_1,\dots, Y_r)$ with $|\vec{Y}| = n$.
\end{Theorem}

As in the calculation in \cite[the end of Chapter~6]{Lecture}, we get
the following nice expression for the generating function.
\begin{Corollary}\label{cor:genP_t(M(r,n))}
The generating function of the Poincar\'e polynomials of
$\pi^{-1}(\linebreak[3]n[0])$ is given by
\begin{equation*}
   \sum_n P_t(\pi^{-1}(n[0])) q^n 
   = \prod_{\alpha=1}^r \prod_{d=1}^\infty
    \frac1{1 - t^{2(rd - \alpha)} q^d}.
\end{equation*}
\end{Corollary}

\begin{proof}[Proof of \thmref{thm:P_t(M(r,n))}]
Since $\pi^{-1}(n[0])$ is homotopic to $M(r,n)$, it is enough to
calculate the dimensions of (ordinary) homology groups. For this
calculation, we use the $\alpha$-partition given in
\thmref{thm:topology}.

We need to specify the one parameter subgroup $\lambda\colon \C^* \to
\hT$ in the proof of \thmref{thm:topology}. Let
\begin{equation}\label{eq:lambda}
   \lambda(t) = (t^{m_1}, t^{m_2}, t^{n_1}, \dots, t^{n_r}).
\end{equation}
If we choose weights $m_1, m_2, n_\alpha$ generic, the Zariski closure 
of $\lambda(\C^*)$ is equal to the whole $\hT$, and the fixed
point set of $\lambda(\C^*)$ coincides with that of $\hT$.

Furthermore, we assume
\begin{equation}\label{eq:weights}
    m_2 \gg n_1 > n_2 > \dots > n_r \gg m_1 > 0
\end{equation}
in order to make the calculation of the index simpler.
(This choice of the weights is due to T.~Gocho, and was given as an
answer to an exercise in a preliminary version of \cite{Lecture}.)

By the proof of \thmref{thm:topology} the corresponding
$(-)$-attracting set is an $\alpha$-partition of $\pi^{-1}(0)$.

\begin{NB}
The following is the original argument.

We give a K\"ahler metric to $M(r,n)$ as a hyper-K\"ahler
quotient as in \cite[Chap.~3]{Lecture}.
Consider the corresponding moment map
\begin{equation*}
\begin{split}
   & F\left(\U(V)\cdot(B_1,B_2,i,j)\right)\\
   = \; & m_1 \| B_1 \|^2 + m_2 \| B_2 \|^2 + (m_1 + m_2) \| j\|^2
   + \tr\left((-i^\dagger i + j j^\dagger)
     \operatorname{diag}(n_1,\dots,n_r)\right).
\end{split}
\end{equation*}
By Morse theory for moment maps of an $S^1$-action (see e.g.,
\cite[5.1]{Lecture}), the critical point set of $F$ is precisely the
fixed point set of $\lambda(\C^*)$, hence is a finite set by our
choice of weights.
Again by a general theory
\begin{enumerate}
\item $F$ is non-degenerate in the sense of Bott,
\item the index of a critical point is equal to twice the sum of
dimension of negative weight spaces (as a $\lambda(\C^*)$-module) of
the tangent space,
\item $F$ is a perfect Morse function.
\end{enumerate}
In our situation, $F$ is non-degenerate in the original sense since the 
critical points are isolated. Moreover the statement (3) follows from
the Morse inequality since the indices are always even.

We show that $F$ is a proper function on $M(r,n)$. This guarantees the 
applicability of Morse theory. Combining
\begin{equation*}
   \tr\left(i^\dagger i\,\operatorname{diag}(n_1,\dots,n_r)\right)
   \le n_1 \tr(i^\dagger i)
\end{equation*}
with the equation 
\begin{equation*}
   \tr(i^\dagger i) - \tr(j^\dagger j) = n\zeta
\end{equation*}
deduced from the moment map equation $\mu_1 = \frac{\sqrt{-1}}2 \zeta$
(see \cite[(3.16)]{Lecture} for $\mu_1$), we get
\begin{equation*}
\begin{split}
    & m_1 \| B_1 \|^2 + m_2 \| B_2 \|^2 + \| i \|^2 
   + \left(m_1 + m_2 - (n_1 + 1)\right)\| j\|^2 \\
   \le \; & F\left(\U(V)\cdot(B_1,B_2,i,j)\right) + (n_1 + 1) n \zeta.
\end{split}
\end{equation*}
This implies that $F$ is proper under \eqref{eq:weights}. Then the
gradient flow for $F$ stays in a compact set, and we can use the Morse
theory for $F$ in order to calculate the Betti numbers of $M(r,n)$.

Or, we can use an argument in \cite[4.6]{Kazhdan-Lusztig} to obtain
the Bialynicki-Birula decomposition in this situation. Their result is
applicable since the $\lambda(\C^*)$-action contracts $M(r,n)$ to a
compact subset as $t$ tends to $0$ by the properness of $F$.
\end{NB}

Our remaining task is to calculate indices of critical points.  Let
$(E,\Phi)$ be a critical point, i.e., a fixed point of the
$\hT$-action. Then $T_{(E,\Phi)} M(r,n)$ has a $\hT$-module structure,
and has an induced $\C^*$-module structure via $\lambda$. By our
choice of weights \eqref{eq:weights}, negative weight spaces for
$\C^*$-action are direct sum of weight spaces for $\hT$-action such
that the one of the followings holds
\renewcommand{\descriptionlabel}[1]{\hspace\labelsep \upshape #1}
\begin{description}
\item[$(1)$] weight of $t_2$ is negative,
\item[$(2)$] weight of $t_2$ is zero and weight of $e_1$ is negative,
\item[$(3)$] weight of $t_2$, $e_1$ are zero and weight of $e_2$ is
negative,
\item[$(4)$] weight of $t_2$, $e_1$, $e_2$ are zero and weight of $e_3$ is
negative,
\item[$\cdots$]
\item[$(r+1)$] weight of $t_2$, $e_1$, $e_2$, \dots, $e_{r-1}$ are
zero and weight of $e_r$ is negative,
\item[$(r+2)$] weight of $t_2$, $e_1$, $e_2$, \dots, $e_{r}$ are
zero and weight of $t_1$ is negative.
\end{description}

Recall that we decompose the tangent space into 
\(
   \sum_{\alpha,\beta} N_{\alpha,\beta}^{\vec{Y}}(t_1,t_2)
\)
in \thmref{thm:M(r,n)weights}. We calculate the sum of dimensions of
weight spaces with the above condition in each summand separately, and
then sum up the contribution from each summand. In the summand $\alpha
= \beta$, the contribution is
\begin{equation*}
   |Y_\alpha| - l(Y_\alpha).
\end{equation*}
If $\alpha < \beta$, the above condition is equivalent to that weight
of $t_2$ is nonpositive. Hence the contribution is equal to the number
of terms in \cite[(1.18)]{part1}, i.e.,
\begin{equation*}
   |Y_\beta|.
\end{equation*}
If $\alpha > \beta$, the above condition is equivalent to that weight
of $t_2$ is negative. If we look at \cite[(1.18)]{part1}, we find that
the contribution is
\begin{equation*}
   |Y_\beta| - l(Y_\beta).
\end{equation*}
Thus the total contribution is
\begin{equation*}
   \sum_{\beta=1}^r r |Y_\beta| - (r-\beta+1)l(Y_\beta).
\end{equation*}
Changing the variable as $\alpha = r-\beta+1$, we get the
formula~\eqref{eq:P_t(M(r,n))}.
\end{proof}

We now turn to the moduli spaces on blowup.
The proof of the following is the same as that of
\thmref{thm:topology} and hence omitted.
\begin{Theorem}\label{thm:topology2}

\rom{(1)} $\bM(r,k,n)$ is homotopy equivalent to
$\widehat\pi^{-1}(n[0])$.

\rom{(2)} Both $\bM(r,k,n)$ and $\widehat\pi^{-1}(n[0])$ have
$\alpha$-partitions into affine spaces.

\rom{(3)} $H_{\mathop{\mathrm{odd}}}(\bM(r,k,n),\Z) = 0$ and
$H_{\mathop{\mathrm{even}}}(\bM(r,k,n),\Z)$ is a free abelian
group. The cycle map $A_*(\bM(r,k,n))\to
H_{\mathop{\mathrm{even}}}(\bM(r,k,n),\Z)$ is an isomorphism. The same
holds for $\widehat\pi^{-1}(n[0])$.
\end{Theorem}

The following theorem is proved in a similar way as
\thmref{thm:P_t(M(r,n))}. The detail is left to the reader.
\begin{Theorem}\label{thm:P_t(bM(r,k,n))}
The Poincar\'e polynomial of $\widehat\pi^{-1}(n[0])$ is given by
\begin{equation*}
    P_t(\widehat\pi^{-1}(n[0])) 
   = \sum
    \prod_{\alpha=1}^r
      t^{2(r|Y_\alpha^1| + r|Y_\alpha^2| - \alpha l(Y_\alpha^1))}
    \prod_{\alpha < \beta} t^{(k_\alpha - k_\beta)(k_\alpha - k_\beta+1)},
\end{equation*}
where the summation runs over the set $((k_1,Y_1^1,Y_1^2),\dots,
(k_r,Y_r^1,Y_r^2))$ with \eqref{eq:cond}.
\end{Theorem}

\begin{Corollary}\label{cor:gen}
The generating function of the Poincar\'e polynomials of
$\widehat\pi^{-1}(n[0])$ is given by
\begin{multline*}
   \sum_n P_t(\widehat\pi^{-1}(n[0])) q^n 
\\
   = 
    \left( \prod_{\alpha=1}^r \prod_{d=1}^\infty
              \frac1{1-t^{2(rd-\alpha)}q^d}\right)
    \left( \prod_{d=1}^\infty \frac1{1-t^{2rd}q^d} \right)^r 
    \times\! \sum_{\{\vec{k}\} = -\frac{k}r}\!
    t^{2\langle \vec{k},\rho\rangle}
    (t^{2r} q)^{(\vec{k},\vec{k})/2}.
\end{multline*}
\end{Corollary}

\begin{NB}
\begin{proof}[Proof of \thmref{thm:P_t(bM(r,k,n))}]
As in the proof of \thmref{thm:P_t(M(r,n))} we take a one parameter
subgroup $\lambda\colon\C^*\to T^{r+2}$ given by \eqref{eq:lambda}
with \eqref{eq:weights}. Probably it might be possible to endow
$\bM(r,k,n)$ with a complete K\"ahler metric to such that the
corresponding moment map is proper. But we use an argument in
\cite[4.6]{Kazhdan-Lusztig} instead to obtain the Bialynicki-Birula
decomposition of $\bM(r,k,n)$. In order to apply their result, we must
check that $\lambda(\C^*)$ contracts $\bM(r,k,n)$ to a compact set as
$t$ tends to $0$. But this follows from the corresponding assertion
for $M(r,n)$ since $\pi$ and $\widehat\pi$ are equivariant and proper.

Our remaining task is to compute the dimension of the sum of
negative weight spaces with respect to $\lambda$, i.e., weight spaces
of $T^{r+2}$ with either of conditions (1) $\sim$ (r+2) in the proof
of \thmref{thm:P_t(M(r,n))}.
Since the tangent space decomposes into $L_{\alpha,\beta}(t_1,t_2)$,
$t_1^{k_\beta-k_\alpha} N^{\vec{Y^1}}_{\alpha,\beta}(t_1,t_2/t_1)$,
$t_2^{k_\beta-k_\alpha} N^{\vec{Y^2}}_{\alpha,\beta}(t_1/t_2,t_2)$ as
in \thmref{thm:bM(r,k,n)weights}, we study negative weight spaces in
each summand separately.

First consider the summand $L_{\alpha,\beta}$. If $\alpha = \beta$,
then $L_{\alpha,\beta} = 0$.
If $\alpha < \beta$, the condition is equivalent to that weight of
$t_1$ is nonpositive.
If $k_\alpha > k_\beta$, then the total
dimension of the sum of negative weight spaces is
\begin{equation*}
   \frac12 (k_\alpha - k_\beta + 1)(k_\alpha - k_\beta).
\end{equation*}
If $k_\alpha + 1 < k_\beta$, then the dimension is $0$.
In the remaining cases $k_\alpha = k_\beta$, $k_\alpha + 1 = k_\beta$,
the dimension is also $0$.
If $\alpha > \beta$, the condition is equivalent to that weight of
$t_1$ is negative. If $k_\alpha > k_\beta$, then the total dimension
is
\begin{equation*}
   \frac12 (k_\alpha - k_\beta - 1)(k_\alpha - k_\beta)
   = \frac12 (k_\beta - k_\alpha + 1)(k_\beta - k_\alpha).
\end{equation*}
Otherwise, the dimension is $0$.
Assuming $\alpha < \beta$, we combine the contribution from
$L_{\alpha,\beta}$ with one from $L_{\beta,\alpha}$. We find that the
sum is equal to
\begin{equation*}
   \frac12 (k_\alpha - k_\beta + 1)(k_\alpha - k_\beta)
\end{equation*}
in all cases.

Next consider the summand $t_1^{k_\beta-k_\alpha}
N^{\vec{Y^1}}_{\alpha,\beta}(t_1,t_2/t_1)$.  Noticing that weights of
$t_1$ do not affect the condition, we find that the contribution is
the same as one in \thmref{thm:P_t(M(r,n))}.
Thus if we sum up contributions for $\alpha$, $\beta$, we get
\begin{equation*}
   \sum_{\alpha=1}^r
    \left( r |Y_\alpha^1| - \alpha l(Y_\alpha^1)\right).
\end{equation*}

Finally consider the summand $t_2^{k_\beta-k_\alpha}
N^{\vec{Y^2}}_{\alpha,\beta}(t_1/t_2,t_2)$.
If $\alpha = \beta$, then
\begin{equation*}
   N_{\alpha,\alpha}^{\vec{Y^2}}(t_1,t_2) = \sum_{s\in Y_\alpha^2}
   t_1^{-l_{Y_\alpha^2}(s)} t_2^{a_{Y_\alpha^2}(s)+l_{Y_\alpha^2}(s)+1}
   + 
   t_1^{l_{Y_\alpha^2}(s)+1} t_2^{-a_{Y_\alpha^2}(s)-l_{Y_\alpha^2}(s)-1}.
\end{equation*}
The contribution is equal to the number of terms with a negative power 
of $t_2$, thus $|Y_\alpha^2|$.

Assuming $\alpha < \beta$, we consider the summands
$t_2^{k_\beta-k_\alpha} N^{\vec{Y^2}}_{\alpha,\beta}(t_1/t_2,t_2)$ and
$t_2^{k_\alpha-k_\beta} N^{\vec{Y^2}}_{\beta,\alpha}(t_1/t_2,t_2)$
simultaneously. We separate weight spaces into four types:
\begin{enumerate}
\item a weight space in $t_2^{k_\beta-k_\alpha}
N^{\vec{Y^2}}_{\alpha,\beta}(t_1/t_2,t_2)$ with weight of $t_2$ is
nonpositive,
\item a weight space in $t_2^{k_\beta-k_\alpha}
N^{\vec{Y^2}}_{\alpha,\beta}(t_1/t_2,t_2)$ with weight of $t_2$ is
positive,
\item a weight space in $t_2^{k_\alpha-k_\beta}
N^{\vec{Y^2}}_{\beta,\alpha}(t_1/t_2,t_2)$ with weight of $t_2$ is
negative,
\item a weight space in $t_2^{k_\alpha-k_\beta}
N^{\vec{Y^2}}_{\beta,\alpha}(t_1/t_2,t_2)$ with weight of $t_2$ is
nonnegative.
\end{enumerate}
We need to calculate the dimension of sum of weight spaces of type (1)
and (3).
By the symmetry 
\(
    N^{\vec{Y^2}}_{\alpha,\beta}(t_1/t_2 ,t_2) 
    = N^{\vec{Y^2}}_{\beta,\alpha}((t_1/t_2)^{-1}, t_2^{-1}) t_1
\)
we find that the dimension of the sum of weight spaces of type (1)
(resp.\ (2)) is the same as that of type (4) (resp.\ (3)).
Thus the dimension of sum of weight spaces of type (1) and (3) is
equal to the half of that of the whole space, i.e.
$|Y_\alpha^2| + |Y_\beta^2|$.
Summing up all $\alpha$, $\beta$, we get $\sum_{\alpha} r |Y_\alpha^2|$.
Altogether, we finally obtain the assertion.
\end{proof}
\end{NB}

\subsection{A different choice of the one parameter subgroup}

This subsection is an interesting detour. We compute Betti numbers of
$\widehat\pi^{-1}(n[0])$ in a different way. A comparison with the
formula in the previous subsection gives us a nontrivial combinatorial
identity.

Let us choose weights for the one-parameter subgroup $\lambda$
in \eqref{eq:lambda} so that
\begin{equation*}
   m_1 = m_2 \gg n_1 > n_2 > \dots > n_r > 0
\end{equation*}
and $m_1$, $n_\alpha$ are generic.

Since this $\lambda$ is not generic, the fixed points are different
from those for $\hT$. But they are described similarly as
\( 
   (E,\Phi) = (I_1(k_1 C),\Phi_1) \oplus \cdots \oplus
   (I_r(k_r C),\Phi_r)
\)
such that
\begin{aenume}
\item $I_\alpha(k_\alpha C)$ is the tensor product
$I_\alpha\otimes\shfO(k_\alpha C)$, where $k_\alpha\in\Z$ and
$I_\alpha$ is an ideal sheaf of $0$-dimensional subscheme $Z_\alpha$
contained in ${\widehat\C}^2 = \bp\setminus\linf$.
\item $\Phi_\alpha$ is an isomorphism from $(I_\alpha)_{\linf}$ to the
$\alpha$th factor of $\shfO_{\linf}^{\oplus r}$.
\item $I_\alpha$ is fixed by the diagonal subgroup $\Delta\C^*$ of
$\C^*\times\C^*$, coming from that on $\bp$.
\end{aenume}

Furthermore, we can parametrize the components of the fixed point set
by $(\vec{k},\vec{Y}) = ((k_1, Y_1),\dots,(k_r,Y_r))$ with $k_\alpha$
as above and $Y_\alpha$ is a Young diagram. Here the constraint is
\begin{equation}\label{eq:cond'}
    \sum_\alpha k_\alpha = k, \qquad
    |\vec{Y}|
    + \frac1{2r} \sum_{\alpha < \beta} |k_\alpha - k_\beta|^2 = n.
\end{equation}
Since this can be proved by the method in \cite[Chapter~7]{Lecture},
we explain it only briefly.
A general point in the component $(\vec{k},\vec{Y})$ is 
\( 
   (E,\Phi) = (I_1(k_1 C),\Phi_1) \oplus \cdots \oplus
   (I_r(k_r C),\Phi_r)
\)
such that
\begin{aenume}
\item the support of $I_\alpha$ consists of $P_1, P_2,\dots,
P_{l(Y_\alpha)}$, contained in the exceptional curve $C$,

\item if $\xi$ is the inhomogeneous coordinate of $C = \proj^1$ and
$\eta$ is the coordinate of the fiber $\widehat{\C}^2 \cong O(-1)\to
C$, 
\[
   I_\alpha = (\xi - \xi_1, \eta^{\lambda^\alpha_1})
   \cap 
   (\xi - \xi_2, \eta^{\lambda^\alpha_2})
   \cap\cdots
\]
with $\xi_l = \xi(P_l)$.
\end{aenume}
See \cite[Figure~7.4]{Lecture}. The points $P_1, P_2,\cdots$ move in
$\proj^1$, but their order is irrelevant when the values
$\lambda^\alpha_l$ are the same. Therefore the component is isomorphic to
\begin{equation*}
   S^{Y_1}\proj^1 \times \dots \times S^{Y_r}\proj^1,
\end{equation*}
with the following notation:
For a Young diagram $Y = (\lambda_1\ge\lambda_2\ge\cdots)$, we define
$m_i = \# \{ l \mid \lambda_l = i\}$. We denote $Y =
(1^{m_1}2^{m_2}\cdots)$ in this case. We set
\begin{equation*}
   S^Y \proj^1 = S^{m_1}\proj^1 \times S^{m_2}\proj^1 \times\dots
   = \proj^{m_1}\times\proj^{m_2} \times \cdots,
\end{equation*}
where $S^m \proj^1$ is the $m$th symmetric product of $\proj^1$.

Let $(E,\Phi)$ be a fixed point in the component corresponding to
$((a_1,Y_1), \dots, (a_r, Y_r))$. Then the tangent space $T_{(E,\Phi)}
\bM(r,k,n)$ is a $\Delta\C^*\times T^r$-module. The $\Delta\C^*\times
T^r$-module structure is independent of the choice of a point, we take
the $\hT$-fixed point corresponding to $((k_1,\emptyset, Y_1), \dots,
(k_r,\emptyset, Y_r))$.  By the formula in
\thmref{thm:bM(r,k,n)weights} we have
\begin{equation*}
   T_{(E,\Phi)} \bM(r,k,n)
   = \sum_{\alpha,\beta} 
   (L_{\alpha,\beta}^{\vec{k}}(t_1,t_1) 
     + t_1^{a_\beta - a_\alpha} N^{\vec{Y}\prime}_{\alpha,\beta}(t_1)),
\end{equation*}
where $N^{\vec{Y}\prime}_{\alpha,\beta}(t_1) =
N^{\vec{Y}}_{\alpha,\beta}(1, t_1)$.
By \thmref{thm:M(r,n)weights} we have
\begin{equation}
\label{eq:N'}
\begin{split}
   N^{\vec{Y}\prime}_{\alpha,\beta}(t_1)
   & = \left(\sum_{s\in Y_\alpha} t_1^{a_{Y_\alpha}(s)+1}
  + \sum_{t\in Y_\beta} t_1^{-a_{Y_\beta}(t)}\right)
   e_\beta e_\alpha^{-1}
   .
\end{split}
\end{equation}

The following theorem is proved in a similar way as
\thmref{thm:P_t(M(r,n))}. The detail is left to the reader.
\begin{Theorem}
The Poincar\'e polynomial of $\widehat\pi^{-1}(n[0])$ is given by
\begin{equation*}
    P_t(\widehat\pi^{-1}(n[0])) 
   = \sum
    \prod_{\alpha=1}^r
      t^{2(|Y_\alpha|-l(Y_\alpha))} P_t(S^{Y_\alpha}\proj^1)
    \prod_{\alpha < \beta} t^{2(l'_{\alpha,\beta}+|Y_\alpha| +
   |Y_\beta| - n'_{\alpha,\beta})}
\end{equation*}
where the summation runs over the set
$(\vec{k},\vec{Y})$ with \eqref{eq:cond'}, and
\begin{gather*}
   l'_{\alpha,\beta}
   = 
   \begin{cases}
     \frac12 (k_\alpha - k_\beta + 1)(k_\alpha - k_\beta)
       & \text{if $k_\alpha \ge k_\beta$,} \\
     \frac12 (k_\beta - k_\alpha + 1)(k_\beta - k_\alpha) - 1
       & \text{otherwise},
   \end{cases}
\\
   n'_{\alpha,\beta}
   = 
   \begin{cases}
     (\text{$\#$ of columns of $Y_\alpha$ which are longer than
           $k_\alpha - k_\beta$})
       & \text{if $k_\alpha \ge k_\beta$,} \\
     (\text{$\#$ of columns of $Y_\beta$ which are longer than
           $k_\beta - k_\alpha - 1$})
       & \text{otherwise}.
   \end{cases}
\end{gather*}
\end{Theorem}

\begin{NB}
\begin{proof}
As in the proof of \thmref{thm:P_t(bM(r,k,n))}, we compute the
dimension of the sum of negative weight spaces with respect to
$\lambda$ for each term $L_{\alpha,\beta}^{\vec{k}}(t_1,t_1)$, $t_1^{k_\beta -
k_\alpha} N^{\vec{Y}\prime}_{\alpha,\beta}(t_1)$.

First consider $L_{\alpha,\beta}^{\vec{k}}(t_1,t_1)$.
If $\alpha = \beta$, then $L_{\alpha,\beta}^{\vec{k}} = 0$.
Assuming $\alpha < \beta$, we combine the contribution from
$L_{\alpha,\beta}^{\vec{k}}$ with one from $L_{\beta,\alpha}$.
For $L_{\alpha,\beta}^{\vec{k}}$, we count the number of terms of nonpositive
powers of $t_1$. If $k_\alpha > k_\beta$, then it is
\begin{equation*}
   \frac12 (k_\alpha - k_\beta + 1)(k_\alpha - k_\beta),
\end{equation*}
and $0$ otherwise.
For $L^{\vec{k}}_{\beta,\alpha}$, we count the number of terms of negative
powers of $t_1$. If $k_\alpha < k_\beta$, then the contribution
is
\begin{equation*}
   \frac12 (k_\beta - k_\alpha + 1)(k_\beta - k_\alpha) - 1, 
\end{equation*}
and $0$ otherwise. Thus we get $l'_{\alpha,\beta}$.

Next consider $t_1^{k_\beta - k_\alpha} N^{\vec{Y}\prime}_{\alpha,\beta}(t_1)$. If
$\alpha = \beta$, the contribution is $|Y_\alpha| - l(Y_\alpha)$.

Assuming $\alpha < \beta$, we combine the contribution from
$t_1^{k_\beta-k_\alpha}N^{\vec{Y}\prime}_{\alpha,\beta}(t_1)$ and one from
$t_1^{k_\alpha-k_\beta}N^{\vec{Y}\prime}_{\beta,\alpha}(t_1)$ as in the proof of
\thmref{thm:P_t(bM(r,k,n))}. By the symmetry
\(
   N^{\vec{Y}\prime}_{\alpha,\beta}(t_1) = N^{\vec{Y}\prime}_{\beta,\alpha}(t_1^{-1})t_1
\),
we find that the contribution from these two terms is equal to
\begin{equation*}
\begin{split}
   & |Y_\alpha| + |Y_\beta| 
   - (\text{the coefficient of $t_1^{k_\alpha-k_\beta+1}$
      in $N^{\vec{Y}\prime}_{\alpha,\beta}(t_1)$}) \\
  =\; & |Y_\alpha| + |Y_\beta| 
   - (\text{the coefficient of $t_1^{k_\beta - k_\alpha}$
      in $N^{\vec{Y}\prime}_{\beta,\alpha}(t_1)$}).
\end{split}
\end{equation*}

Suppose $k_\alpha\ge k_\beta$. By \eqref{eq:N'} we have
\begin{equation*}
\begin{split}
   & (\text{the coefficient of $t_1^{k_\alpha-k_\beta+1}$
          in $N^{\vec{Y}\prime}_{\alpha,\beta}(t_1)$}) \\
  =\; & (\text{number of columns of $Y_\alpha$ which are longer than
           $k_\alpha - k_\beta$}). 
\end{split}
\end{equation*}
On the other hand, if $k_\alpha < k_\beta$, then we have
\begin{equation*}
\begin{split}
   & (\text{the coefficient of $t_1^{k_\alpha - k_\beta + 1}$
          in $N^{\vec{Y}\prime}_{\beta,\alpha}(t_1)$}) \\
  =\; & (\text{number of columns of $Y_\beta$ which are longer than
           $k_\beta - k_\alpha - 1$}). 
\end{split}
\end{equation*}
Thus we get $n'_{\alpha,\beta}$.
\end{proof}
\end{NB}

Let us consider the generating function of Poincar\'e polynomials. We
consider the simplest case.
\begin{Corollary}\label{cor:rank2}
Assume $r=2$ and $c_1=0$.
The generating functions of Poincar\'e polynomial is
\begin{multline*}
   \left(\prod_{d=1}^\infty
        \frac1{(1-t^{4d}q^d)(1-t^{4d-2}q^d)^2(1-t^{4d-4}q^d)}\right)
\\
  \times
  \Biggl[ 
    \sum_{k \ge 0} 
    \prod_{d=1}^{2k}
        \frac{1-t^{4d-4}q^d}{1-t^{4d}q^d}
        t^{2k(2k+1)} q^{k^2}
   + \sum_{k > 0} 
     \prod_{d=1}^{2k-1}
        \frac{1-t^{4d-4}q^d}{1-t^{4d}q^d}
        t^{2k(2k+1)-2} q^{k^2}
  \Biggr].
\end{multline*}
\end{Corollary}

Comparing with the formula in \corref{cor:gen}, we get the following
identity
\begin{equation}\label{eq:Ochiai}
\begin{split}
   & \sum_{k \ge 0}\prod_{d=1}^{2k}
        \frac{1-t^{4d-4}q^d}{1-t^{4d}q^d} t^{2k(2k+1)} q^{k^2}
   + \sum_{k > 0} \prod_{d=1}^{2k-1}
        \frac{1-t^{4d-4}q^d}{1-t^{4d}q^d} t^{2k(2k+1)-2} q^{k^2} \\
  =\; & \prod_{d=1}^\infty \frac{1-t^{4d-2}q^d}{1-t^{4d}q^d}
    \sum_{k=-\infty}^\infty t^{2k(2k+1)} q^{k^2}.
\end{split}
\end{equation}
Since this identity does not involve any geometric information, it is
natural to expect to have a direct proof. Such a proof was provided
for us by Hiroyuki Ochiai. See \secref{sec:Ochiai}.

Finally let us remark that the results of this and the previous
subsections give the blowup formula for the virtual Hodge polynomials
of moduli spaces for an arbitrary projective surface $X$.
Let $H$ be an ample line bundle over $X$. For $c_1\in H^2(X,\Z)$,
$n\in\Q$, let $M_H(r,c_1,n)$ be the moduli space of $H$-stable sheaves
$E$ on $X$ with $c_1(E) = c_1$, $c_2(E) - \frac{r-1}{2r} c_1(E)^2 = n$.
We assume $\operatorname{GCD}(r,\langle c_1, H\rangle) = 1$.

Let $\bM_H(r,c_1+kC,n)$ be the moduli space of $(H-\ve C)$-stable
sheaves $E$ on $\widehat X$ with $c_1(E) = p^* c_1+kC$, $c_2(E) -
\frac{r-1}{2r} c_1(E)^2 = \Delta$, where $c_1$, $n$ is as above,
$k\in\Z$, and $\ve > 0$ is sufficiently small.

Let $e(Y;x,y)$ denote the virtual Hodge polynomial of $Y$ introduced in
\cite{virHodge}.

\begin{Theorem}\label{thm:euler_blowup}
The ratio
\begin{equation*}
    \left.\sum_n e(\bM_H(r,c_1+kC,n);x,y) q^n\right/
    \sum_n e(M_H(r,c_1,n);x,y) q^n
\end{equation*}
is independent of the surface $X$ and is given by
\begin{equation*}
    \left( \prod_{d=1}^\infty \frac1{1-(xy)^{rd}q^d} \right)^r 
    \times\! \sum_{\{\vec{k}\} = -\frac{k}r}\!
    (xy)^{\langle \vec{k},\rho\rangle}
    ((xy)^{r} q)^{(\vec{k},\vec{k})/2}.
\end{equation*}
\end{Theorem}
This result is proved as follows. From the proof of
\thmref{thm:blowup} for arbitrary surface $X$, we have a
stratification of $M_H(r,c_1,n)$ such that $\widehat\pi$ is a
fibration over each stratum. The fiber is independent of $X$, and
isomorphic to our $\widehat\pi^{-1}(0)$ defined for the framed moduli
spaces. Then properties of virtual Hodge polynomials give the above
assertion.

Finally remark that the above holds in the Grothendieck group of
varieties, if we replace $(xy)^n$ by $[\C^n]$. This generalizes
\cite{Go2} from rank $2$ to higher ranks.

\begin{Remark}
This result was obtained by the second author \cite{Yoshioka:1996}. In
fact, he assumed that the moduli spaces are nonsingular and used the
Weil conjecture to count numbers of rational points over finite
fields. He did not use the framed moduli spaces nor the morphism
$\widehat\pi$ as did in here. The above proof, under the same
assumption, was obtained in July, 1997. The authors then noticed that
W-P.~Li and Z.~Qin \cite{LQ1,LQ2,LQ3} obtained the above result for
rank $2$ case, where the universal function is given in the form
corresponding to \corref{cor:rank2}. The authors then learned the
virtual Hodge polynomials are natural language here.
\end{Remark}

\section{Nekrasov's deformed partition function}\label{sec:Nek}

Nekrasov's deformed partition function \cite{Nek}, more precisely, the
one with {\it higher order Casimir operators} turned on, can be
considered as the generating function of the {\it equivariant
homology\/} version of Donaldson invariants on $\C^2$. We give its
definition and also the one for the blowup $\widehat{\C}^2$ at the
origin in this section. We then study their relation by using the
localization theorem in the equivariant homology groups. (See
\secref{sec:BMhom}.)

As for the calculation of {\it original\/} Donaldson invariants, the
localization technique was not so useful even if we assume the base
manifold has large symmetry (say $X = \CP^2$). This was because the
fixed point sets are not isolated in general, and are still difficult
to study. The crucial difference here is the existence of the framing:
A point in $M(r,n)$ is fixed by the action given by the change of the
framing if and only if it is a direct sum of rank $1$ sheaves. The
rank $1$ sheaves are easy to study.

\subsection{Equivariant integration}\label{subsec:integ}
Before giving the definition, we explain a general setting for the
`equivariant integration' via the localization theorem. 

Let $T$ be a torus acting on an algebraic variety $N$.
Suppose that the fixed point set $N^T$ consists of a single point
$o$. Then the push-forward homomorphism for the inclusion
$\iota_0\colon \{o\} \to N$ induces an isomorphism between localized
equivariant homology groups
\begin{equation*}
   (\iota_o)_*\colon \mathcal S
   \cong H^T_*(o)\otimes_S \mathcal S
   \xrightarrow{\cong} H^T_*(N)\otimes_S \mathcal S,
\end{equation*}
where $S = H^*_T(\mathrm{pt})$ and $\mathcal S$ is its quotient
field.
Furthermore, if $f\colon M\to N$ is a $T$-equivariant proper morphism,
we can define $H^T_*(M) \to \mathcal S$ by
\begin{equation*}
   \alpha \longmapsto (\iota_0)_*^{-1} f_*\alpha.
\end{equation*}
We denote this by $\int_M$. This makes sense even when $M$ is not
necessarily compact. But it takes a value in the rational function
field $\mathcal S$. When $M$ is compact, it coincides with the usual
integration and has values in $S$.

For equivariant $K$-homology groups, we have a similar homomorphism
defined by the same formula:
\begin{equation*}
   K^G(N) \ni \alpha \longmapsto
   (\iota_0)_*^{-1} f_*\alpha \in \mathcal R,
\end{equation*}
where $\mathcal R$ is the quotient field of the representation ring of
$G$. This has a relation with equivariant Hilbert polynomials. See
\cite[\S3]{part1}.

Suppose $M$ is nonsingular. Let $M^T = \bigsqcup_i F_i$ be the
decomposition of the fixed point set $M^T$ to irreducible
components. Let $N_i$ be the normal bundle. Note that we only have
finitely many components, and each $F_i$ is compact, as $f$ is proper
and $N^T=\{o\}$. By the functoriality of the push-forward homomorphism,
we have
\begin{equation*}
   \int_M \alpha =
   \sum_i \int_{F_i} \frac1{e_T(F_i)} \iota_i^*\alpha, 
\end{equation*}
where $e_T(N_i)$ is the equivariant Euler class and $\iota_i^*$ is the
pull-back homomorphism for the inclusion $\iota_i\colon F_i\to M$
defined via the Poincar\'e duality homomorphism. Here $\int_{F_i}$ is
the usual integration as $F_i$ is compact.
When $M$ is compact, the fractional parts of each summand of the right
hand side cancel out, and the final answer is in $S$. But this does
not happen when $M$ is noncompact in general.

\subsection{Universal sheaves}

Since we need higher rank generalization of Donaldson's $\mu$-map, we
begin with the description of universal sheaves on the moduli spaces.

Over the moduli space $M(r,n)$, we have a natural vector bundle $V$,
whose fiber at $(E,\varphi)$ is
\(
   H^1(E(-\ell_\infty)).
\)
In the ADHM description in \cite[Chapter 2]{Lecture}, this is the
bundle associated with the natural principal $\GL_n(\C)$-bundle,
coming from the construction of $M(r,n)$ as a quotient space. If
$\mathcal E$ denotes a universal sheaf on $\proj^2\times M(r,n)$, we
have $V = R^1 p_{2*}(\mathcal E\otimes p_1^*(\shfO(-\linf)))$.

We also have another natural vector bundle $W$, given by the fiber at
infinity: $W = H^0(\mathcal E|_{\linf})$. This is a trivial bundle,
but nontrivial as an equivariant bundle. We also consider bundles
\(
   S^+ = \shfO_{\proj^2}(-1)\oplus \shfO_{\proj^2}(1)
\)
and
\(
   S^- = \shfO_{\proj^2}\oplus \shfO_{\proj^2}
\)
over $\proj^2$ with the $T^2$-action such that
\( 
  \ch S^+_0 = 1 + t_1^{-1} t_2^{-1},
\)
\(
   \ch S^-_0 = t_1^{-1} + t_2^{-1}
\)
on the fibers at the origin. (They are positive and negative spinor
bundles $S^+$, $S^-$ when restricted on $\R^4$.)
We then form a virtual equivariant vector bundle on $\proj^2\times
M(r,n)$ by
\begin{equation*}
   \shfO_{\proj^2}\boxtimes W + (S^- - S^+)\boxtimes V.
\end{equation*}

\begin{NB}
The formula \cite[2.29]{Nek2} seems wrong. The above was used in
\cite[the third line in p.22]{Nek2}. The above is natural in the
context of quiver varieties.
\end{NB}

By \cite[Lemma~1.8]{part1}, this virtual equivariant bundle is
isomorphic to the universal sheaf $\mathcal E$ in the equivariant
$K$-cohomology group $K_{\hT}(\proj^2\times M(r,n))$. We denote the
virtual bundle by $\mathcal E$ hereafter.

The character of the fiber of $\mathcal E$ at the fixed point $(0,\vec{Y})$
is given by
\begin{equation}\label{eq:char}
   \iota_{(0,\vec{Y})}^* \ch(\mathcal E)
   = \sum_{\alpha=1}^r e^{a_\alpha} \left(
     1
     - (1-e^{\ve_1})(1-e^{\ve_2})
     \sum_{s\in Y_\alpha} e^{l'(s)\ve_1 + a'(s)\ve_2}
   \right),
\end{equation}
where we set $e_\alpha = e^{- a_\alpha}$, $t_1 = e^{-\ve_1}$, $t_2 =
e^{-\ve_2}$ as usual.
Let $q_1 = [0 : 1 : 0]$, $q_2 = [0 : 0 : 1]$ be two other fixed points
in $\proj^2$. We have
\begin{equation*}
   \iota^*_{(q_1,\vec{Y})}\ch(\mathcal E)
   = \iota^*_{(q_2,\vec{Y})}\ch(\mathcal E)
   = \sum_{\alpha=1}^r e^{a_\alpha}.
\end{equation*}
We define
\begin{equation*}
   \ch(\mathcal E)/[\C^2]
   = \frac1{\ve_1\ve_2} \iota^*_{\{0\}\times M(r,n)}\ch(\mathcal E).
\end{equation*}
Since $\C^2$ is noncompact, the slant product $/[\C^2]$ is not defined
in the usual sense. So we define it by formally applying Bott's
formula. The homogeneous degree part of this is an element of the
localized equivariant cohomology group
$H^*_{\hT}(M(r,n))\otimes_{S}\mathcal S$, but its fractional part is a
constant in the following sense:
\begin{equation*}
  \begin{split}
    & \iota^*_{\vec{Y}}\left(\ch(\mathcal E)/[\C^2]\right)
\\
    = 
    \; &
    \iota^*_{\vec{Y}}\left(\ch(\mathcal E)/[\proj^2]\right)
    + \frac1{\ve_1(\ve_2-\ve_1)}\iota^*_{(q_1,\vec{Y})}\ch(\mathcal E)
    + \frac1{(\ve_1-\ve_2)\ve_2}\iota^*_{(q_2,\vec{Y})}\ch(\mathcal E)
\\
    =\; & \iota^*_{\vec{Y}}\left(\ch(\mathcal E)/[\proj^2]\right)
    + \frac1{\ve_1\ve_2} \sum_{\alpha=1}^r e^{a_\alpha}
   ,
  \end{split}
\end{equation*}
and each degree part of
\begin{equation*}
  \iota^*_{\vec{Y}}\left(\ch(\mathcal E)/[\proj^2]\right)
  = \frac{(1-e^{\ve_1})(1-e^{\ve_2})}{\ve_1\ve_2}
        \sum_{\alpha=1}^r e^{a_\alpha}
      \sum_{s\in Y_\alpha} e^{l'(s)\ve_1 + a'(s)\ve_2}
\end{equation*}
is a polynomial. We remark
\begin{equation*}
   \ch(\mathcal E)/[\proj^2]
   = \ch\left((S^- - S^+)\boxtimes V\right)/[\proj^2]
\end{equation*}
since $\ch\left(W\boxtimes \shfO_{\proj^2}\right)/[\proj^2] = 0$. Now
the slant product $/[\proj^2]$ in the right hand side can be replaced
by $/[\C^2]$ since $S^- - S^+$ is zero at $\linf = \proj^2\setminus
\C^2$. This observation also matches with the above formula of
$\iota^*_{\vec{Y}}\left(\ch(\mathcal E)/[\proj^2]\right)$.

Let us define the instanton part of the {\it partition function\/} by
\begin{equation*}
  \begin{split}
    &
    \Zin(\ve_1,\ve_2,\vec{a};\q,\vec{\tau}) = \sum_{n=0}^\infty
    \q^n
    \int_{M(r,n)}
   \exp\left(
     \sum_{p=1}^\infty \tau_{p}
     \ch_{p+1}(\mathcal E)/[{\C}^2] \right)
\\
   =\; & \sum_{\vec{Y}}
   \frac{\q^{|\vec{Y}|}}{\displaystyle\prod_{\alpha,\beta}
    n^{\vec{Y}}_{\alpha,\beta}(\ve_1,\ve_2,\vec{a})}
\\
   & \quad
  \times 
  \exp\left(
  \sum_{p=1}^\infty \sum_{\alpha=1}^r \tau_p \left[
    \frac{e^{a_\alpha}}{\ve_1\ve_2} \left( 1 - (1-e^{\ve_1})(1-e^{\ve_2})
     \sum_{s\in Y_\alpha} e^{l'(s)\ve_1 + a'(s)\ve_2}\right)\right]_{p-1}
  \right)
  ,
  \end{split}
\end{equation*}
where $\ch_{p+1}$ is the degree $(p+1)$-part of the Chern character,
$[\heartsuit]_{p-1}$ denotes the degree $(p-1)$-part of
$\heartsuit$, and
\(
  \int_{M(r,n)} \spadesuit
\)
means
\(
  (\iota_{0*})^{-1} \pi_*\left(\spadesuit\cap[M(r,n)]\right)
\)
with $\iota_0\colon \{0\}\to M_0(r,n)$ is the inclusion of the unique
fixed point $0\in M_0(r,n)^{\hT}$.
Furthermore, $n^{\vec{Y}}_{\alpha,\beta}(\ve_1,\ve_2,\vec{a})$ is the
equivariant Euler class of the tangent space at the fixed point
$\vec{Y}$. It is given by the explicit formula:
\begin{multline*}
n^{\vec{Y}}_{\alpha,\beta}(\ve_1,\ve_2,\vec{a})
   = \prod_{s \in Y_\alpha}
      \left( -l_{Y_\beta}(s)\ve_1 + (a_{Y_\alpha}(s)+1)\ve_2 + a_\beta 
   - a_\alpha\right)
\\
    \times\prod_{t\in Y_\beta} 
        \left((l_{Y_\alpha}(t)+1)\ve_1 -a_{Y_\beta}(t)\ve_2 + a_\beta 
   - a_\alpha\right)
.
\end{multline*}
The indeterminate $\q$ should be distinguished from $q =
e^{\pi\sqrt{-1}\tau}$ which will appear later.

Note that each coefficient of $\q^n$ has a complicated, but explicit
expression. For small $n$, it is easy to compute. (The authors wrote a
MAPLE program, which was very useful when we found our main result.) But
when $n$ increases, the number of Young diagrams becomes large, and it
becomes difficult to compute. We are interested not in individual
coefficients, but in the generating function.

By \eqref{eq:char} we have
\begin{equation*}
   \iota^*_{\vec{Y}} \left(\ch_1 (\mathcal E)/[{\C}^2]\right)
   = \frac1{\ve_1\ve_2}\sum_\alpha a_\alpha = 0,
\qquad
   \iota^*_{\vec{Y}} \left(\ch_2 (\mathcal E)/[{\C}^2]\right)
   = \frac1{2\ve_1\ve_2} \sum_\alpha a_\alpha^2
   - n 
   .
\end{equation*}
By the first equality, we did not include $\tau_0$ in $\vec{\tau}$.
Also, the second equation means that $\tau_1$ is essentially equal to
$-\log\q$ (see also \subsecref{subsec:versus}), but we use both of
them to simplify the blowup formula below.

If we set $\vec{\tau} = 0$, we get the partition function studied in
\cite{part1}, which was denoted by $Z(\ve_1,\ve_2,\vec{a};\q)$ there.
But we emphasize that this is only the {\it instanton part\/} of the
partition function. It is more natural to include the perturbative
part also. This will be done in \subsecref{subsec:Pert}. This is the
reason why we change the notation.

If we expand the exponential,
$\Zin(\ve_1,\ve_2,\vec{a};\q,\vec{\tau})$ becomes
\begin{equation*}
    \sum_{n=0}^\infty
    \q^n \sum_{\mu} \prod_i \tau_{\mu_i}
    \int_{M(r,n)}
     \prod_i \ch_{\mu_i+1}(\mathcal E)/[{\C}^2]
   ,
\end{equation*}
where $\mu = (\mu_1\ge\mu_2\ge\cdots)$ is a partition. Therefore
$\Zin(\ve_1,\ve_2,\vec{a};\q,\vec{\tau})$ is the generating function
of {\it all\/} intersection numbers of Chern classes of universal
sheaves slanted by the fundamental cycle $[\C^2]$.

The following formula will be useful later:
\begin{multline}\label{eq:pdZ}
   \left(\prod_i \pd{}{\tau_{\mu_i}}\right)
    \Zin(\ve_1,\ve_2,\vec{a};\q,\vec{\tau})
\\
 = \sum_{n=0}^\infty
    \q^n
    \int_{M(r,n)}
   \left(
     \prod_i \ch_{\mu_i+1}(\mathcal E)/[{\C}^2] \right)
     \cap
   \exp\left(
     \sum_{p=1}^\infty \tau_{p}
     \ch_{p+1}(\mathcal E)/[{\C}^2] \right)
\end{multline}
for a partition $\mu = (\mu_1\ge\mu_2\ge\cdots)$.

\subsection{Partition function on the blowup}

We consider similar partition functions on the blowup:
\begin{equation*}
  \begin{split}
   & \bZin_{c_1 = k}(\ve_1,\ve_2,\vec{a};\q,\vec{\tau},\vec{t})
\\
  =\; & \sum_{n} \q^n
  \int_{\bM(r,k,n)}
   \exp\left(
     \sum_{p=1}^\infty \left\{ t_{p}\left(\ch_{p+1}(\bE)/[C]\right)
   + 
   \tau_{p}
     \left(\ch_{p+1}(\bE)/[\widehat{\C}^2]\right)\right\}
   \right)
   ,
  \end{split}
\end{equation*}
where $\bE$ is a universal sheaf over $\widehat{\proj}^2\times
\bM(r,k,n)$ and $\ch(\bE)/[\widehat{\C}^2]$ is defined as above via
the localization (see below). And the summation runs over $n\in
\Z_{\ge 0}+\frac1{2r}k(r-k)$. Here we do not include $t_0$ though
$\ch_1(\bE)/[C]$ is not $0$ in general. In fact, it is constant $-k$.

We calculate this by using the localization formula. Recall
\(
   \{ ([z_0 : z_1 : z_2], [z : w] \in \proj^2\times\proj^1\mid 
   z_1 w = z_2 z \}.
\)
Let $p_1 = ([1:0:0], [1:0])$, $p_2 = ([1:0:0], [0:1])$, $q_1 =
([0:1:0], [1:0])$, $q_2 = ([0:0:1], [0:1])$ be the fixed points in
$\widehat{\proj}^2$. We use the same notation for the latter two
points as fixed points in $\proj^2$, since they are mapped to
corresponding points under the projection
$p\colon\widehat{\proj}^2\to\proj^2$.
The characters of the fibers over the fixed
points are given by
\begin{gather*}
  \iota^*_{(p_1, \vec{k},\vec{Y^1},\vec{Y^2})}( \ch(\bE))
  = 
  \left.\iota^*_{(0, \vec{Y^1})}(\ch(\mathcal E))
    \right|_{\substack{\ve_1 \to \ve_1\\ \ve_2\to \ve_2-\ve_1\\
        \vec{a} \to \vec{a}+\ve_1\vec{k}}},
\\
  \iota^*_{(p_2, \vec{k},\vec{Y^1},\vec{Y^2})}( \ch(\bE) )
  = 
  \left.\iota^*_{(0, \vec{Y^2})}(\ch(\mathcal E))
    \right|_{\substack{\ve_1 \to \ve_1-\ve_2\\ \ve_2\to \ve_2\\
        \vec{a} \to \vec{a}+\ve_2\vec{k}}},
\\
  \iota^*_{(q_1, \vec{k},\vec{Y^1},\vec{Y^2})}( \ch(\bE) )
  = \iota^*_{(q_2, \vec{k},\vec{Y^1},\vec{Y^2})}( \ch(\bE) )
  = \sum_{\alpha=1}^r e^{-a_\alpha}.
\end{gather*}
For the first two equalities, see the proof of \cite[2.4]{part1}. 
\begin{NB}
The Chern character is normalized so that $\ch(\det\bE) = 0$. So this
$\vec{k}$ is the shifted one, i.e., 
$\vec{k} = {}^t\!(k_1-\frac{k}r,\dots,k_r-\frac{k}r)$.
\end{NB}
The last two equalities are obvious since $p$ is isomorphism outside
the exceptional set. Therefore $\ch(\bE)/[\widehat{\C}^2]$ has the
same constant fractional part as $\ch(\mathcal E)/[\C^2]$.

Therefore we have
{\allowdisplaybreaks
\begin{equation*}
\begin{split}
   & \iota^*_{(\vec{k},\vec{Y^1},\vec{Y^2})}
     \left(\ch(\bE)/[C]\right)
\\
   =\; & \sum_{\alpha=1}^r
   \begin{aligned}[t]
     &  \left[ 
       \frac{e^{a_\alpha + \ve_1 k_\alpha} \ve_1}{\ve_1(\ve_2-\ve_1)}
       \left( 1 - (1-e^{\ve_1})(1-e^{\ve_2-\ve_1})
       \sum_{s\in Y^1_\alpha} e^{l'(s)\ve_1+a'(s)(\ve_2-\ve_1)} \right)
     \right.
     \\
     & \qquad\left.
     + \frac{e^{a_\alpha + \ve_2 k_\alpha} \ve_2}{(\ve_1-\ve_2)\ve_2}
     \left( 1 - (1-e^{\ve_1-\ve_2})(1-e^{\ve_2})
     \sum_{s\in Y^2_\alpha} e^{l'(s)(\ve_1-\ve_2)+a'(s)\ve_2} \right)
    \right]
     \end{aligned}
\\
   =\; & \ve_1
   \left.\iota_{\vec{Y^1}}^* \left(\ch(\mathcal E)/[{\C}^2]\right)
     \right|_{\substack{\ve_1 \to \ve_1\\ \ve_2 \to \ve_2-\ve_1
         \\ \vec{a} \to \vec{a}+\ve_1\vec{k}}}
   + \ve_2
   \left.\iota_{\vec{Y^2}}^* \left(\ch(\mathcal E)/[{\C}^2]\right)
     \right|_{\substack{\ve_1 \to \ve_1-\ve_2\\ \ve_2 \to \ve_2
         \\ \vec{a} \to \vec{a}+\ve_2\vec{k}}}
\end{split}
\end{equation*}
and}
{\allowdisplaybreaks
\begin{equation*}
  \iota^*_{(\vec{k},\vec{Y^1},\vec{Y^2})}
     \left(\ch(\bE)/[\widehat{\C}^2]\right)  
   = 
   \begin{aligned}[t]
   & \left.
   \iota^*_{\vec{Y^1}}
     \left(\ch(\mathcal E)/[{\C}^2]\right)
   \right|_{\substack{\ve_1 \to \ve_1\\ \ve_2 \to \ve_2-\ve_1
         \\ \vec{a} \to \vec{a}+\ve_1\vec{k}}}
     + 
     \left.\iota_{\vec{Y^2}}^* \left(\ch(\mathcal E)/[{\C}^2]\right)
     \right|_{\substack{\ve_1 \to \ve_1-\ve_2\\ \ve_2 \to \ve_2
         \\ \vec{a} \to \vec{a}+\ve_2\vec{k}}}
    .
   \end{aligned}
\end{equation*}
\begin{NB}
\begin{equation*}
\begin{split}
   & \iota^*_{(\vec{k},\vec{Y^1},\vec{Y^2})}
     \left(\ch(\bE)/[C]\right)
\\
   =\; & \ve_1
   \left.\iota_{\vec{Y^1}}^* \left(\ch(\mathcal E)/[{\Bbb P}^2]\right)
     \right|_{\substack{\ve_1 \to \ve_1\\ \ve_2 \to \ve_2-\ve_1
         \\ \vec{a} \to \vec{a}+\ve_1\vec{k}}}
   + \ve_2
   \left.\iota_{\vec{Y^2}}^* \left(\ch(\mathcal E)/[{\Bbb P}^2]\right)
     \right|_{\substack{\ve_1 \to \ve_1-\ve_2\\ \ve_2 \to \ve_2
         \\ \vec{a} \to \vec{a}+\ve_2\vec{k}}}
+\sum_{\alpha=1}^r
\frac{e^{\ve_1 k_\alpha}-e^{\ve_2 k_\alpha}}{\ve_2-\ve_1}e^{a_\alpha}
\end{split}
\end{equation*}
and
\begin{equation*}
  \iota^*_{(\vec{k},\vec{Y^1},\vec{Y^2})} \left(\ch(\bE)/[\bp]\right)  
   = 
   \begin{aligned}[t]
   & \left.
   \iota^*_{\vec{Y^1}}
     \left(\ch(\mathcal E)/[{\proj}^2]\right)
   \right|_{\substack{\ve_1 \to \ve_1\\ \ve_2 \to \ve_2-\ve_1
         \\ \vec{a} \to \vec{a}+\ve_1\vec{k}}}
     + 
     \left.\iota_{\vec{Y^2}}^* \left(\ch(\mathcal E)/[{\Bbb P}^2]\right)
     \right|_{\substack{\ve_1 \to \ve_1-\ve_2\\ \ve_2 \to \ve_2
         \\ \vec{a} \to \vec{a}+\ve_2\vec{k}}}
     \\
     & \qquad + \sum_{\alpha=1}^r
     \left(\frac{\ve_2 e^{\ve_1 k_\alpha}-\ve_1 e^{\ve_2 k_\alpha}}
       {\ve_1\ve_2(\ve_2-\ve_1)}
       -\frac{1}{\ve_1\ve_2}\right)e^{a_\alpha}.
   \end{aligned}
\end{equation*}
\end{NB}
Note} also that the equivariant Euler class of the tangent space
$(\vec{k}, \vec{Y^1},\vec{Y^2})$ is given by
\begin{equation*}
   \prod_{\alpha,\beta} 
     l^{\vec{k}}_{\alpha,\beta}(\ve_1,\ve_2,\vec{a})
     \;
     n^{\vec{Y}^1}_{\alpha,\beta}(\ve_1,\ve_2-\ve_1,
       \vec{a} + \ve_1 \vec{k})
     \;
     n^{\vec{Y}^2}_{\alpha,\beta}(\ve_1-\ve_2,\ve_2,
        \vec{a} + \ve_2\vec{k}),
\end{equation*}
where
\begin{equation}\label{eq:s}
\begin{split}
   & l^{\vec{k}}_{\alpha,\beta}(\ve_1,\ve_2,\vec{a}) = 
   s^{k_\alpha-k_\beta}(\ve_1,\ve_2,a_\beta-a_\alpha),
\\
   & s^k(\ve_1,\ve_2,x) = 
   \begin{cases}
     {\displaystyle
     \prod_{\substack{i,j\ge 0\\i+j \le k-1}}}
          (-i\ve_1 -j\ve_2 + x)
       & \text{if $k > 0$}, \\
     {\displaystyle
     \prod_{\substack{i,j\ge 0\\i+j \le -k-2}}}
          \left((i+1)\ve_1 + (j+1)\ve_2 + x\right)
       & \text{if $k < -1$}, \\
     1 & \text{$k=0$ or $-1$}.
   \end{cases}
\end{split}
\end{equation}
Therefore
\begin{equation}\label{eq:blowup_inst}
\begin{split}
    & \bZin_{c_1=k}(\ve_1,\ve_2,\vec{a};\q,\vec{\tau},\vec{t})
\\
   =\; & \sum_{\{ \vec{k}\} = - \frac{k}r}
    \frac{\q^{\frac12(\vec{k},\vec{k})}}
   {\displaystyle 
     \prod_{\alpha,\beta}
     l^{\vec{k}}_{\alpha,\beta}(\ve_1,\ve_2,\vec{a})}
   \Zin(\ve_1,\ve_2-\ve_1,\vec{a}+\ve_1\vec{k};\q,\vec{\tau}+\ve_1\vec{t})
 \\
   & \qquad\qquad\qquad\qquad\qquad\qquad  \times
   \Zin(\ve_1-\ve_2,\ve_2,\vec{a}+\ve_2\vec{k};\q,\vec{\tau}+\ve_2\vec{t}),
\end{split}
\end{equation}
where $\vec{\tau}+\ve_1\vec{t}$ means $(\tau_1+\ve_1 t_1, \tau_2+\ve_1
t_2, \cdots)$, and
This is a generalization of the blowup formula. (In our previous
paper \cite{part1}, we only consider the case $\vec{t} = (t,0,0,\dots)$.)

\begin{NB}
If a reader do not like the integration over noncompact spaces as
$/[\C^2]$, $/[\widehat{\C}^2]$, we can define the partition function by
\begin{equation*}
   {}_{\mathrm c}\Zin
   (\ve_1,\ve_2,\vec{a};\q,\vec{\tau}) = \sum_{n=0}^\infty
   \q^n
   \int_{M(r,n)}
   \exp\left(
     \sum_{p=1}^\infty \tau_{p}
     \ch_{p+1}(\mathcal E)/[{\proj}^2] \right)
   .
\end{equation*}
Then the blowup formula becomes
\begin{equation*}
\begin{split}
    & {}_{\mathrm c}\bZin_{c_1=k}
    (\ve_1,\ve_2,\vec{a};\q,\vec{\tau},\vec{t})
\\
   =\; & \sum_{\{ \vec{k}\} = - \frac{k}r}
   \frac{\q^{\frac12(\vec{k},\vec{k})}
     \displaystyle \exp\left(\sum_{\alpha=1}^N
     \sum_{p=1}^\infty
     \left(
     \left[
    \frac{e^{\ve_1 k_\alpha}-e^{\ve_2 k_\alpha}}{\ve_2-\ve_1}e^{a_\alpha}
    \right]_{p} t_p
     +
     \left[
     \left(\frac{\ve_2 e^{\ve_1 k_\alpha}-\ve_1 e^{\ve_2 k_\alpha}}
       {\ve_1\ve_2(\ve_2-\ve_1)}
       -\frac{1}{\ve_1\ve_2}\right)e^{a_\alpha}
      \right]_{p-1} \tau_p\right)
    \right)}
   {\displaystyle 
     \prod_{\alpha\in\Delta} l^{\vec{k}}_{\alpha}(\ve_1,\ve_2,\vec{a})}\\
  & \quad \quad  \times
   {}_{\mathrm c}Z
     (\ve_1,\ve_2-\ve_1,\vec{a}+\ve_1\vec{k};\q,\vec{\tau}+\ve_1\vec{t})
   {}_{\mathrm c}Z
     (\ve_1-\ve_2,\ve_2,\vec{a}+\ve_2\vec{k};\q,\vec{\tau}+\ve_2\vec{t}).
\end{split}
\end{equation*}
Now it is better to multiply ${}_{\mathrm c}\Zin$ by
\(
   \exp\left(\sum_{p=1}^\infty \sum_{\alpha} \frac{\tau_p}{\ve_1\ve_2}
     \frac{a_\alpha^{p+1}}{(p+1)!}\right)
\)
by hand, so that the enumerator in the blowup formula is absorbed as
\begin{equation*}
  \begin{split}
    & \widetilde\gamma_{\ve_1,\ve_2-\ve_1}
   (a_\alpha+\ve_1 k_\alpha;\vec{\tau}+\ve_1\vec{t})
   + \widetilde\gamma_{\ve_1-\ve_2,\ve_2}
   (a_\alpha+\ve_2 k_\alpha;\vec{\tau}+\ve_2\vec{t})
   - \widetilde\gamma_{\ve_1,\ve_2}
   (a_\alpha;\vec{\tau})
\\
   =\; & \sum_{p=1}^\infty \left\{
     \frac{\tau_p + \ve_1 t_p}{\ve_1(\ve_2-\ve_1)}
      \frac{(a_\alpha+\ve_1 k_\alpha)^{p+1}}{(p+1)!}
     + \frac{\tau_p + \ve_2 t_p}{(\ve_1-\ve_2)\ve_2}
      \frac{(a_\alpha+\ve_2 k_\alpha)^{p+1}}{(p+1)!}
     - \frac{\tau_p}{\ve_1\ve_2}\frac{a_\alpha^{p+1}}{(p+1)!}
     \right\}
\\
   =\; & \sum_{p=1}^\infty \left\{
    \left[
    \frac{e^{\ve_1 k_\alpha}-e^{\ve_2 k_\alpha}}{\ve_2-\ve_1}e^{a_\alpha}
    \right]_{p} t_p
     +
     \left[
     \left(\frac{\ve_2 e^{\ve_1 k_\alpha}-\ve_1 e^{\ve_2 k_\alpha}}
       {\ve_1\ve_2(\ve_2-\ve_1)}
       -\frac{1}{\ve_1\ve_2}\right)e^{a_\alpha}
      \right]_{p-1} \tau_p
     \right\},
  \end{split}
\end{equation*}
where
\begin{equation*}
   \widetilde\gamma_{\ve_1,\ve_2}(x;\vec{\tau})
  =  \sum_{p=1}^\infty \frac{x^{p+1}}{(p+1)!}\frac{\tau_p}{\ve_1\ve_2}.
\end{equation*}
But this modification is the same as one defined by using $\ch(\mathcal
E)/[\C^2]$.

Let us record the first few terms:
\begin{gather*}
  \left[
    \frac{e^{\ve_1 k_\alpha}-e^{\ve_2 k_\alpha}}{\ve_2-\ve_1}e^{a_\alpha}
    \right]_{0} = - k_\alpha 
\quad
  \left[
    \frac{e^{\ve_1 k_\alpha}-e^{\ve_2 k_\alpha}}{\ve_2-\ve_1}e^{a_\alpha}
    \right]_{1}
  = - \left(k_\alpha a_\alpha
    + \frac{k_\alpha^2}2(\ve_1+\ve_2)\right), 
\\
     \left[
     \left(\frac{\ve_2 e^{\ve_1 k_\alpha}-\ve_1 e^{\ve_2 k_\alpha}}
       {\ve_1\ve_2(\ve_2-\ve_1)}
       -\frac{1}{\ve_1\ve_2}\right)e^{a_\alpha}
      \right]_{0}
     = - \frac12 k_\alpha^2
\end{gather*}
\end{NB}

\subsection{Adding perturbation term}\label{subsec:Pert}

We define the {\it full\/} partition function by
\begin{equation*}
    Z(\ve_1,\ve_2,\vec{a};\q,\vec{\tau})
    = \exp\left[-
      \sum_{\alpha\neq\beta} \gamma_{\ve_1,\ve_2}(a_\alpha-a_\beta;\Lambda)
      \right]\, \Zin(\ve_1,\ve_2,\vec{a};\q,\vec{\tau}),
\end{equation*}
where $\gamma_{\ve_1,\ve_2}$ is as in \secref{sec:perturb} and
$\Lambda = \q^{\frac1{2r}}$ as usual. The first term is the
perturbation term of the partition function.

\begin{NB}
We must put a minus sign here.
\end{NB}

We choose the branch of $\log$ in the perturbative term as explained
in the paragraph after \propref{prop:F_beh}. See \eqref{eq:double}.

We have
{\allowdisplaybreaks
\begin{equation*}
  \begin{split}
  & \sum_{\alpha\neq\beta}
      \gamma_{\ve_1,\ve_2-\ve_1}(a_\alpha-a_\beta+\ve_1(k_\alpha-k_\beta);
      \Lambda 
      )
      + 
      \gamma_{\ve_1-\ve_2,\ve_2}(a_\alpha-a_\beta+\ve_2(k_\alpha-k_\beta);
      \Lambda 
      )
\\
  =\; &
  \begin{aligned}[t]
  & \sum_{\alpha\neq\beta} \Biggl[
  \gamma_{\ve_1,\ve_2}(a_\alpha-a_\beta;\Lambda)
  + \log s^{k_\beta-k_\alpha}(\ve_1,\ve_2,a_\alpha-a_\beta)
  - \frac{(k_\alpha-k_\beta)^2}2 \log\Lambda
  \Biggr]
  \end{aligned}
  \end{split}
\end{equation*}
by \eqref{eq:pert_shift}.} Therefore
{\allowdisplaybreaks
\begin{equation}\label{eq:product}
  \begin{split}
    & Z(\ve_1,\ve_2-\ve_1,\vec{a}+\ve_1\vec{k};\q,\vec{\tau}+\ve_1\vec{t})\,
      Z(\ve_1-\ve_2,\ve_2,\vec{a}+\ve_2\vec{k};\q,\vec{\tau}+\ve_2\vec{t})
\\
  = 
  & 
    \begin{aligned}[t]
      &\exp\left[-
      \sum_{\alpha\neq\beta}
      \gamma_{\ve_1,\ve_2}(a_\alpha-a_\beta;\Lambda)
      \right]
      \prod_{\alpha,\beta} \frac{\Lambda^{(k_\beta-k_\alpha)^2/2}}
      {l^{\vec{k}}_{\beta,\alpha}(\ve_1,\ve_2,\vec{a})}
     \\
     & 
     \times
      \Zin(\ve_1,\ve_2-\ve_1,\vec{a}+\ve_1\vec{k}
      ;\q,\vec{\tau}+\ve_1\vec{t}
      )\,
      \Zin(\ve_1-\ve_2,\ve_2,\vec{a}+\ve_2\vec{k}
      ;\q,\vec{\tau}+\ve_2\vec{t}
      ).
  \end{aligned}
  \end{split}
\end{equation}
Then} the blowup formula \eqref{eq:blowup_inst} is
simplified as
\begin{equation}\label{eq:blowup_formula}
  \begin{split}
   & \widehat Z^{c_1=k}(\ve_1,\ve_2,\vec{a};\q,\vec{\tau},\vec{t})
\\
   = & \sum_{\{ \vec{k}\} = -\frac{k}r} \!
   Z(\ve_1,\ve_2-\ve_1,\vec{a}+\ve_1\vec{k}
   ;\q,\vec{\tau}+\ve_1\vec{t}
   )\,
   Z(\ve_1-\ve_2,\ve_2,\vec{a}+\ve_2\vec{k}
   ;\q,\vec{\tau}+\ve_2\vec{t}
   ),
  \end{split}
\end{equation}
where
\begin{equation*}
  \widehat Z^{c_1=k}(\ve_1,\ve_2,\vec{a};\q,\vec{\tau},\vec{t}) 
  = \exp\left[-
      \sum_{\alpha\neq\beta} \gamma_{\ve_1,\ve_2}(a_\alpha-a_\beta;\Lambda)
      \right]\,
      \bZin_k(\ve_1,\ve_2,\vec{a};\q,\vec{\tau},\vec{t})
    .
\end{equation*}

\begin{Remark}
As we saw, the blowup equations are simplified if we include the
perturbation part. This was pointed out to us by N.~Nekrasov. Probably
this is already enough for the reason for the perturbation term. But
it has a geometric meaning as a regularization of the Euler class of
an `infinite rank' vector bundle. (See \cite[\S3.10]{Nek} and
\subsecref{subsec:coord} below.)
\end{Remark}

\subsection{$\tau_1$ versus $\log\q$}\label{subsec:versus}

Let $\vec{\tau}_1 = (\tau_1,0,0,\cdots)$ be a vector with the first
entry only. We have
\begin{equation*}
   \Zin(\ve_1,\ve_2,\vec{a};\q,\vec{\tau}+\vec{\tau}_1)
   = \exp\left[
      \frac{\tau_1}{2\ve_1\ve_2} \sum_\alpha a_\alpha^2
      \right]\, 
   \Zin(\ve_1,\ve_2,\vec{a};\q e^{-\tau_1},\vec{\tau}).
\end{equation*}
On the other hand, we have
\begin{equation*}
\begin{split}
&     \exp\left[-
      \sum_{\alpha\neq\beta}
      \gamma_{\ve_1,\ve_2}(a_\alpha-a_\beta;
        \Lambda\exp(-\frac{\tau_1}{2r}))\right]
      \left/
          \exp\left[-
      \sum_{\alpha\neq\beta}
      \gamma_{\ve_1,\ve_2}(a_\alpha-a_\beta; \Lambda)\right]\right.
\\
      =\; &
      \exp\left[-
      \sum_{\alpha\neq\beta}
      \frac{\tau_1}{2r}\left\{
        \frac{(a_\alpha-a_\beta)^2}{2\ve_1\ve_2}
       +
        \frac{(a_\alpha-a_\beta)(\ve_1+\ve_2)}{2\ve_1\ve_2}
       +
        \frac{\ve_1^2+\ve_2^2+3\ve_1\ve_2}{12\ve_1\ve_2}
        \right\}
      \right]
\\
    =
    \; &
    \exp\left[
      -\frac{\tau_1}{2\ve_1\ve_2} \sum_\alpha a_\alpha^2
      \right]
    \exp\left[
      - \frac{\tau_1(r-1)(\ve_1^2+\ve_2^2+3\ve_1\ve_2)}{24\ve_1\ve_2}
      \right].
\end{split}
\end{equation*}
by \eqref{eq:e^u}. Therefore
\begin{equation*}
   Z(\ve_1,\ve_2,\vec{a};\q,\vec{\tau}+\vec{\tau}_1)
   = \exp\left[
     - \frac{\tau_1(r-1)(\ve_1^2+\ve_2^2+3\ve_1\ve_2)}{24\ve_1\ve_2}
      \right]\, 
   Z(\ve_1,\ve_2,\vec{a};\q e^{-\tau_1},\vec{\tau}).
\end{equation*}
In particular, we have
\begin{equation}\label{eq:versus}
  (\pd{}{\tau_1})^N Z(\ve_1,\ve_2,\vec{a};\q,\vec{\tau})
  = \left(-\frac{(r-1)(\ve_1^2+\ve_2^2+3\ve_1\ve_2)}{24\ve_1\ve_2}
         - \q\pd{}{\q}\right)^N
   Z(\ve_1,\ve_2,\vec{a};\q,\vec{\tau})
\end{equation}
for $N\in\Z_{\ge 0}$.

\section{The blowup equation and Nekrasov's conjecture}\label{sec:blowupeq}

The blowup formula \eqref{eq:blowup_formula} equates the unknown
function $\widehat Z^{c_1=k}$ to the unknown $Z$. It is useless unless
we know either or their independent relation. We do not have such
knowledge so far in general, but we do know something when we restrict
to the subspace $\vec{\tau} = 0$. This will be given in this
section. An application is a solution of Nekrasov's conjecture:
\(
  \left.\ve_1\ve_2 \log
    Z(\ve_1,\ve_2,\vec{a};\q,0)\right|_{\ve_1,\ve_2=0}
\)
is equal to the Seiberg-Witten prepotential $\mathcal
F_0(\vec{a};\Lambda)$ introduced in \secref{sec:SW}.


Let
\begin{equation*}
  F(\ve_1,\ve_2,\vec{a};\q,\vec{\tau})
  = \ve_1\ve_2 \log Z(\ve_1,\ve_2,\vec{a};\q,\vec{\tau}).
\end{equation*}
Here the logarithm is defined as follows: We first separate this into
the perturbative part and the instanton part
\begin{equation*}
  F(\ve_1,\ve_2,\vec{a};\q,\vec{\tau})
  = - \ve_1\ve_2 \sum_{\alpha\neq\beta}
     \gamma_{\ve_1,\ve_2}(a_\alpha-a_\beta;\Lambda)
     + \ve_1\ve_2 \log \Zin(\ve_1,\ve_2,\vec{a};\q,\vec{\tau}).
\end{equation*}
We denote the second part by $\Fin(\ve_1,\ve_2,\vec{a};\q,\vec{\tau})$.
It has the form
\begin{equation*}
  \begin{aligned}[t]
  & \sum_{p=1}^\infty \sum_{\alpha=1}^r
  \tau_p \left[e^{a_\alpha}\right]_{p+1}
  + \ve_1\ve_2\log \left[
    \sum_{n=0}^\infty \q^n
    \int_{M(r,n)}
   \exp\left(
     \sum_{p=1}^\infty \tau_{p}
     \ch_{p+1}(\mathcal E)/[{\C}^2]
   \right)
  \right].
  \end{aligned}
\end{equation*}
Since the summation in the last part starts with $1$, its logarithm
makes sense as a formal power series in $\q$.

\subsection{Gap from the dimension counting}

In this subsection we derive a differential equation from a simple
geometric consideration, which is well-known in the context of
Donaldson invariants.

\begin{Lemma}\label{lem:c_1blowup}
Let $0 < k < r$. Then
\begin{equation*}
   \bZin_{c_1=k}(\ve_1,\ve_2,\vec{a};\q,\vec{\tau},\vec{t})
   = O((\vec{\tau},\vec{t})^{k(r-k)}).
\end{equation*}
\end{Lemma}

Here we put $\deg\tau_p = p-1$, $\deg t_p = p$ and
$O((\vec{\tau},\vec{t})^N)$ means that it is a sum of monomials of
degree greater than or equal to $N$. When it is a function only in
$\vec{\tau}$ (or $\vec{t}$), we simply denote by $O(\vec{\tau}^N)$ (or
$O(\vec{t}^N)$). This convention will be used in what follows.

\begin{proof}
Consider the projective morphism
$\widehat\pi\colon \bM(r,k,n)\to M_0(r,n-\frac1{2r}k(r-k))$.
If $x\in H^{2d}_{\hT}(\bM(r,k,n))$, we have
\begin{equation*}
   \widehat\pi_*\left(x \cap[\bM(r,k,n)]\right)
   \in H^{\hT}_{2\dim\bM(r,k,n)-2d}(M_0(r,n-\frac1{2r}k(r-k))),
\end{equation*}
and this is $0$ if
\begin{equation*}
   \dim\bM(r,k,n)- d > \dim M_0(r,n-\frac1{2r}k(r-k))
   \Longleftrightarrow
   k(r-k) > d.
\end{equation*}
In the definition of the partition function on the blowup, we have
$\ch_{p+1}(\bE)/[C] \in H^{2p}_{\hT}$,
$\ch_{p+1}(\bE)/[\widehat{\proj}^2] \in H^{2(p-1)}_{\hT}$. The degrees
exactly match with the above definition of the degrees of $t_p$ and
$\tau_q$.
\end{proof}

Similarly we consider the $k=0$ case. The morphism $\widehat\pi\colon
\bM(r,0,n)\to M_0(r,n)$ is an isomorphism outside the inverse image of
the closure of $\{0\}\times\Mreg(r,n-1)$. Furthermore we have
$\widehat\pi^*(\mathcal E) \cong \bE$ there. Since
$\operatorname{codim}(\{0\}\times\Mreg(r,n-1)) = 2r$, the same argument
as above shows
\begin{equation*}
   \bZin_{c_1=0}(\ve_1,\ve_2,\vec{a};\q,\vec{\tau},\vec{t})
   = Z(\ve_1,\ve_2,\vec{a};\q,\vec{\tau})
   + O((\vec{\tau},\vec{t})^{2r}).
\end{equation*}

Combined with the blowup formula in the previous subsection, we get
\begin{align}
   &
   \sum_{\{ \vec{k}\} = 0}
   Z(\ve_1,\ve_2-\ve_1,\vec{a}+\ve_1\vec{k}
   ;\q,\vec{\tau}+\ve_1\vec{t}
   )\,
   Z(\ve_1-\ve_2,\ve_2,\vec{a}+\ve_2\vec{k}
   ;\q,\vec{\tau}+\ve_2\vec{t}
   )
   \label{eq:blowup_formula2}
\\
   & \qquad\qquad\qquad\qquad\qquad\qquad\qquad\qquad\qquad
   =
   Z(\ve_1,\ve_2,\vec{a};\q,\vec{\tau}) + O((\vec{\tau},\vec{t})^{2r}),
   \notag
\\
  &
   \sum_{\{ \vec{k}\} = -\frac{k}r}\!\!\!
   Z(\ve_1,\ve_2-\ve_1,\vec{a}+\ve_1\vec{k}
   ;\q,\vec{\tau}+\ve_1\vec{t}
   )\,
   Z(\ve_1-\ve_2,\ve_2,\vec{a}+\ve_2\vec{k}
   ;\q,\vec{\tau}+\ve_2\vec{t}
   )
   \label{eq:blowup_formula3}
\\
   & \qquad\qquad\qquad\qquad\qquad\qquad\qquad\qquad\qquad
   =
   O((\vec{\tau},\vec{t})^{k(r-k)}) \qquad (0 < k < r).
   \notag
\end{align}
We call these {\it blowup equations}.

\subsection{Recursive structure}\label{subsec:recursive}
In this subsection we illustrate the power of the blowup equations:
They determine $Z(\ve_1,\ve_2,\vec{a};\q,\vec{\tau})$ up to
$O(\vec{\tau}^{2r-3})$.

We introduce two auxiliary functions:
\begin{gather*}
   F_a(\vec{a}) = F(\ve_1,\ve_2-\ve_1,\vec{a};\q,\vec{\tau})
   ,
  \qquad
  F_b(\vec{a}) = F(\ve_1-\ve_2,\ve_2,\vec{a};\q,\vec{\tau})
  .
\end{gather*}
We suppress the $\ve_1$, $\ve_2$-dependence in the notation.

We divide \eqref{eq:blowup_formula2} for $c_1 = k = 0$ by
\(
   Z(\ve_1,\ve_2-\ve_1,\vec{a};\q,\vec{\tau})
   Z(\ve_1-\ve_2,\ve_2,\vec{a};\q,\vec{\tau}),
\)
expand with respect to the variables $\vec{t}$, and take the
coefficients of $t_1$ and $t_1^2$:
\begin{equation}\label{eq:recursive}
\begin{split}
  &
  \begin{aligned}[t]
    & \sum_{\{ \vec{k}\} = 0}
  \frac1{\ve_2-\ve_1}
  \left( \pd{}{\tau_1} F_a(\vec{a}+\ve_1\vec{k})
    -
   \pd{}{\tau_1} F_b(\vec{a}+\ve_2\vec{k})
   \right)
   \\ &\quad
   \times\exp\left[
   \frac1{\ve_2-\ve_1}
   \left(
   \frac{ F_a(\vec{a}+\ve_1\vec{k})
     - F_a(\vec{a}) }{\ve_1}
   -
   \frac{ F_b(\vec{a}+\ve_2\vec{k}) 
    - F_b(\vec{a}) }{\ve_2}
    \right)\right]
  \end{aligned}
\\
  =
  \; &  
  O(\vec{\tau}^{2r-1}),
\\
    &
   \begin{aligned}[t]
    & 
    \sum_{\{ \vec{k}\} = 0}
    \Biggl[
  \left\{ \frac1{\ve_2-\ve_1}
  \left( \pd{}{\tau_1} F_a(\vec{a}+\ve_1\vec{k})
    -
   \pd{}{\tau_1} F_b(\vec{a}+\ve_2\vec{k})
   \right)
   \right\}^2
   \\ & \qquad\qquad
   + \frac1{\ve_2-\ve_1}\left(
     {\ve_1}\pd{^2}{\tau_1^2} F_a(\vec{a}+\ve_1\vec{k})
     - \ve_2\pd{^2}{\tau_1^2} F_b(\vec{a}+\ve_2\vec{k})
   \right)
   \Biggr]   
   \\ &\quad
   \times\exp\left[
   \frac1{\ve_2-\ve_1}
   \left(
   \frac{ F_a(\vec{a}+\ve_1\vec{k})
    - F_a(\vec{a}) }{\ve_1}
   -
   \frac{ F_b(\vec{a}+\ve_2\vec{k})
    - F_b(\vec{a}) }{\ve_2}
    \right)\right]
  \end{aligned}
\\
  =
  \; &  
  O(\vec{\tau}^{2r-2}).
\end{split}
\end{equation}

\begin{Theorem}
\textup{(1)} The solution $F(\ve_1,\ve_2,\vec{a};\q,\vec{\tau})$ of
the equations \eqref{eq:recursive} is unique up to
$O(\vec{\tau}^{2r-2})$. In fact, \eqref{eq:recursive} determine the
coefficients of $\q^n$ in $\Fin(\ve_1,\ve_2,\vec{a};\q,\vec{\tau})$
recursively up to $O(\vec{\tau}^{2r-2})$. 
  
\textup{(2)} The coefficients of monomials in $\vec{\tau}$ of degree
$<2r-2$ in $F(\ve_1,\ve_2,\vec{a};\q,\tau)$ is regular at $\ve_1 =
\ve_2 = 0$.
\end{Theorem}

\begin{proof}
We prove the assertions by the induction on the power of $\q$.
Suppose that the coefficient of $\q^m$ in
$F(\ve_1,\ve_2,\vec{a};\q,\vec{\tau})$ are determined for $m <
n$. We show that the coefficients of $\q^{n}$ in
$F_a$, $F_b$, and hence in $F(\ve_1,\ve_2,\vec{a};\q,\vec{\tau})$ are
determined from \eqref{eq:recursive}.

Let us separate the terms with $\vec{k} = 0$. The remaining terms with
$\vec{k}\neq 0$ are divisible by $\q$ by \eqref{eq:product} (recall
$\prod_{\alpha,\beta} \Lambda^{(k_\beta-k_\alpha)^2/2} =
\q^{\frac12(\vec{k},\vec{k})}$). Then the equations are written as
\begin{equation*}
\begin{split}
  &
    \pd{}{\tau_1} F_a(\vec{a}) - \pd{}{\tau_1} F_b(\vec{a})
    = \q\times \text{known up to order $n-1$},
\\
  &
    \ve_1\pd{^2}{\tau_1^2} F_a(\vec{a})
     - \ve_2\pd{^2}{\tau_1^2} F_b(\vec{a})
     = \q\times \text{known up to order $n-1$}.
\end{split}
\end{equation*}
After noticing that $\pd{}{\tau_1}$ is essentially equal to
$\q\pd{}{\q}$ as \subsecref{subsec:versus}, the above equations gives
a system of linear equations on the coefficients of $\q^n$ in
$F_a(\vec{a},\ve_1\vec{t}')$,$F_b(\vec{a},\ve_2\vec{t}')$. This system
is uniquely solvable since the determinant of
\(
\left(
\begin{smallmatrix}
n & -n
\\
\ve_1 n^2 & -\ve_2 n^2  
\end{smallmatrix}
\right)
\)
is nonzero.

Furthermore, the right hand sides divided by $\ve_1 - \ve_2$ are
regular at $\ve_1,\ve_2 = 0$ if $F(\ve_1,\ve_2,\vec{a};\q,\tau)$ is
regular. Again by the induction, we get the second assertion.
\end{proof}


\subsection{Contact term equations as limit of blowup equations}
In this subsection we study the specialization of the differential 
equation \eqref{eq:blowup_formula2} at $\ve_1 = \ve_2 = 0$.

Let
\begin{equation}\label{eq:expansion}
   F(\ve_1,\ve_2,\vec{a};\q,\vec{\tau})
   = 
   \begin{aligned}[t]
   &
   \Fz(\vec{a};\q,\vec{\tau}) + (\ve_1+\ve_2) H(\vec{a};\q,\vec{\tau})
\\
   &\qquad
   + (\ve_1+\ve_2)^2 G(\vec{a};\q,\vec{\tau})
   + \ve_1\ve_2 \Fo(\vec{a};\q,\vec{\tau}) + \cdots,
   \end{aligned}
\end{equation}
where we consider terms up to $O(\vec{\tau}^{2r-2})$.

By the exactly same argument as in \cite[6.1]{part1}, we have
\(
  \Zin(\ve_1,-2\ve_1,\vec{a};\q,\vec{\tau}) 
  = \Zin(2\ve_1,-\ve_1,\vec{a};\q,\vec{\tau})
\)
up to $O(\vec{\tau}^{2r-2})$. In particular, $H(\vec{a};\q,\vec{\tau})$
up to $O(\vec{\tau}^{2r-2})$ comes only from the perturbative term:
\begin{equation*}
\begin{split}
   & H(\vec{a};\q,\vec{\tau}) = 
   \frac12 \sum_{\alpha < \beta}
   (a_\alpha-a_\beta)\log(-1) 
   = \pi \sqrt{-1} \langle\vec{a}, \rho\rangle.
\end{split}
\end{equation*}
See \eqref{eq:double}.

The first part of the following was proved by \cite{part1} and
independently by \cite{NO}.
\begin{Theorem}\label{thm:main}
\textup{(1)} $\Fz(\vec{a};\q,0)$ is equal to the
Seiberg-Witten prepotential $\mathcal F_0(\vec{a};\Lambda)$ with $\q =
\Lambda^{2r}$.

\textup{(2)} For $p=2,\dots,r$, let $c_p$ be the $p$th power sum
in $z_1,\dots,z_r$ multiplied by $\frac{(-\sqrt{-1})^p}{p!}$ given in
\eqref{eq:c}. We have
\[
   \left.\pd{F_0}{\tau_{p-1}}\right|_{\vec{\tau} = 0}
   = c_p.
\]
\end{Theorem}

\begin{proof}
As the name of this subsection suggests, we prove the assertion by
studying limit of blowup equations.

We have
{\allowdisplaybreaks
\begin{equation*}
\begin{split}
   & \frac{\Fz(\vec{a}+\ve_1\vec{k}; \q, \ve_1 \vec{t})}
   {\ve_1(\ve_2-\ve_1)}
    +\frac{\Fz(\vec{a}+\ve_2\vec{k}; \q, \ve_2\vec{t})}
    {(\ve_1-\ve_2)\ve_2}
\\
  & \qquad\qquad\qquad =
  \frac{1}{\ve_1\ve_2}\Fz
  - \left[
    \frac{\partial^2 \Fz
    }{\partial\tau_p\partial\tau_q} \frac{t_p t_q}2
    +
    \frac{\partial^2 \Fz
    }{\partial\tau_p\partial a^l}t_p k^l
    +
    \frac{\partial^2 \Fz
    }{\partial a^l\partial a^m}
      \frac{k^l k^m}{2}
  \right]
  + \cdots,
\\
   & \frac{\ve_2 H(\vec{a}+\ve_1\vec{k}; \q, \ve_1\vec{t})}
   {\ve_1(\ve_2-\ve_1)}
    +\frac{\ve_1 H(\vec{a}+\ve_2\vec{k}; \q, \ve_2\vec{t})}
   {(\ve_1-\ve_2)\ve_2}
\\
  & \qquad\qquad =
  \frac{\ve_1+\ve_2}{\ve_1\ve_2} H
  + \left[\frac{\partial H
    }{\partial\tau_p} t_p
    +
    \frac{\partial H
    }{\partial a^l} k^l
  \right]
  + \cdots
  = \frac{\ve_1+\ve_2}{\ve_1\ve_2} H
  + \pi\sqrt{-1}\langle\vec{k},\rho\rangle  + \cdots
  ,
\\
   & \frac{\ve_2^2 G(\vec{a}+\ve_1\vec{k};\q,\ve_1\vec{t})}
  {\ve_1(\ve_2-\ve_1)}
    +\frac{\ve_1^2 G(\vec{a}+\ve_2\vec{k};\q,\ve_2\vec{t})}
  {(\ve_1-\ve_2)\ve_2}
  =
  \frac{(\ve_1+\ve_2)^2-\ve_1\ve_2}{\ve_1\ve_2} G
  + \cdots,
\\
   & \Fo(\vec{a}+\ve_1\vec{k}; \q,\ve_1\vec{t})
    + \Fo(\vec{a}+\ve_2\vec{k}; \q,\ve_2\vec{t})
   = 
  2 \Fo
  + \cdots.
\end{split}
\end{equation*}
Here} and throughout the proof, $F_0$, $H$, $G$, $F_1$ and their
derivatives in the right hand side are all restriction to $\vec{\tau}
= 0$.

We divide both hand sides of \eqref{eq:blowup_formula2} by
$Z(\ve_1,\ve_2,\vec{a};\q,\vec{\tau})$, set $\vec{\tau} = 0$, and take
limit $\ve_1, \ve_2 \to 0$:
\begin{multline}\label{eq:limit}
   \sum_{\{ \vec{k}\} = 0}
   \exp\left[
       - \frac{\partial^2 \Fz
    }{\partial\tau_p\partial\tau_q} \frac{t_p t_q}2
    -
    \frac{\partial^2 \Fz
    }{\partial\tau_p\partial a^l}t_p k^l
    -
    \frac{\partial^2 \Fz
    }{\partial a^l\partial a^m}
      \frac{k^l k^m}{2}
    \right]
\\
   \times (-1)^{\langle \vec{k},\rho\rangle}
   \exp\left(\Fo - G \right)
   =
   1 + O(\vec{t}^{2r}),
\end{multline}
where the summation symbol over $p,q,l,m$ are omitted. Logically
speaking, we only show $F$ is regular up to $O(\vec{t}^{2r-2})$ at this
moment. Thus the higher order terms may diverge in the limit $\ve_1 =
\ve_2 = 0$. Therefore this equation should be understood that the left
hand side is equal to $1$ if we set all higher order terms to be zero.

Let
\begin{equation}\label{eq:tau}
   \tau_{kl} = - \frac1{2\pi\sqrt{-1}}
   \frac{\partial^2F_0}{\partial a^k\partial a^l}.
\end{equation}
This is symmetric and positive definite for $\q$ small. So we consider
the corresponding theta function $\Theta_E$ with the characteristic
\(
  E = {}^t\!
   \begin{pmatrix}
   \frac12 & \frac12 & \cdots
   \end{pmatrix}
\)
in the notation for the root system of type $A_{r-1}$.

Comparing the constant term of \eqref{eq:limit}, we get
\begin{equation}\label{eq:1}
   \exp(G - F_1) = \Theta_E(0|\tau).
\end{equation}
Comparing the coefficients of $t_p t_q$ with $p+q\le 2r-1$, we get
\begin{equation}\label{eq:contact}
  0 =
   - \frac{\partial^2 \Fz
    }{\partial\tau_p\partial\tau_q}
  + 
  \frac1{\pi\sqrt{-1}}
  \sum_{l,m}
  \frac{\partial^2 \Fz
    }{\partial\tau_p\partial a^l}
  \frac{\partial^2 \Fz
    }{\partial\tau_q\partial a^m}
  \pd{}{\tau_{lm}} \log\Theta_E(0|\tau).
%
\end{equation}
By \eqref{eq:versus} we have
\(
   \pd{}{\tau_1} F_0 = - \q \pd{}{\q} F_0
   = -\frac1{2r}\Lambda\pd{}{\Lambda} F_0
\)
and
\(
   \frac{\partial^2}{\partial\tau_1^2} F_0
   = \frac1{4r^2}\left(\Lambda\pd{}{\Lambda}\right)^2 F_0.
\)
Therefore the equation with $p=q=1$ is nothing but the contact term
equation in \corref{cor:contact}.

When we consider the contact term equation as the differential
equation for $F_0$, it has similar recursive structure as the blowup
equation studied in \subsecref{subsec:recursive}. The coefficients of
$\q^n$ in the instanton part of $F_0$ are determined from lower
coefficients. In particular, the solution is unique if the
perturbative part is given. Since $F_0$ and the Seiberg-Witten
prepotential $\mathcal F_0$ have the same perturbative part, we have
(1) of \thmref{thm:main}.

Let us prove the second assertion. The case $p=2$ is nothing but the
renormalization group equation in \propref{prop:renorm}. We substitute
$p=p-1$, $q=2$ in \eqref{eq:contact}:
\begin{equation*}
  \Lambda\pd{}{\Lambda} \left(\pd{F_0}{\tau_{p-1}}\right)
  =
  \frac{2r}{\pi\sqrt{-1}}
  \sum_{l,m}
  \pd{u_2}{a^m}
  \pd{}{a^l}\left(\pd{F_0}{\tau_{p-1}}\right)
  \pd{}{\tau_{lm}} \log\Theta_E(0|\tau).
\end{equation*}
If we expand the both hand sides into the power series in $\q$ (plus
the perturbative part), the equation determines the coefficients
recursively. The point here is the observation that $\pd{}{\tau_{lm}}
\log\Theta_E(0|\tau)$ is divisible by $\q$. In particular, the
solution is unique if the perturbative part is given.

The perturbative part of $\pd{F_0}{\tau_{p-1}}$ is given by
\begin{equation*}
   \left[e^{a_\alpha}\right]_{p}
   = \sum_\alpha \frac{a_\alpha^{p}}{p!}.
\end{equation*}
This is equal to the perturbative part of $c_p$. This shows our
assertion.
\end{proof}

\section{Fintushel-Stern's blowup formula}

\subsection{}
In this and next subsections we assume that derivatives of
$F(\ve_1,\ve_2,\vec{a};\q,\vec{\tau})$ up to the second order are
regular at $\ve_1 = \ve_2 = 0$, when restricted to $\vec{\tau} = 0$.

As we did in \eqref{eq:limit} we can derive the limit of the blowup
formula \eqref{eq:blowup_formula} as
\begin{equation}\label{eq:limitBlowup}
  \begin{split}
   & 
   \lim_{\ve_1,\ve_2\to 0}
   \frac{\widehat Z^{c_1=k}(\ve_1,\ve_2,\vec{a};\q,0,\vec{t})}
   {Z(\ve_1,\ve_2,\vec{a};\q,0)}
\\
   = & \sum_{\{ \vec{k}\} = -\frac{k}r} \!
   \exp\left[
       - \frac{\partial^2 \Fz
    }{\partial\tau_p\partial\tau_q} \frac{t_p t_q}2
    -
    \frac{\partial^2 \Fz
    }{\partial\tau_p\partial a^l}t_p k^l
    -
    \frac{\partial^2 \Fz
    }{\partial a^l\partial a^m}
      \frac{k^l k^m}{2}
    \right]
\\
    &
    \qquad\qquad
   \times (\sqrt{-1})^{2\langle \vec{k},\rho\rangle}
   \exp\left(\Fo - G \right),
  \end{split}
\end{equation}
where the summation symbol is omitted as before. And we restrict
functions to the subspace $\vec{\tau} = 0$ also as before.

We define the {\it contact term\/} by
\begin{equation*}
    T_{p,q}(\vec{a};\q)
    = \frac12 \frac{\partial^2 \Fz
    }{\partial\tau_{p}\partial\tau_{q}}(\vec{a};\q,0).
\end{equation*}
We also set $t_r = t_{r+1} = \cdots = 0$. Then the right hand side of
\eqref{eq:limitBlowup} is rewritten as
\begin{equation}\label{eq:factor}
  \exp\left(-\sum_{p,q=1}^{r-1} T_{p,q}{t_p t_q}\right)
  \frac{\Theta_{E_k}
    \left.\left(\displaystyle\frac{\sqrt{-1}}{2\pi}
    \sum_{p=1}^{r-1} \frac{d c_{p+1}}{d\vec{a}} t_p\right|\tau\right)}
  {\Theta_E(0|\tau)}
\end{equation}
where
\begin{equation*}
   \frac{d c_{p+1}}{d\vec{a}} = {}^t\!
   \begin{pmatrix}\displaystyle
     \pd{c_{p+1}}{a^1} & \cdots & \displaystyle\pd{c_{p+1}}{a^{r-1}}
   \end{pmatrix},
\end{equation*}
and $E_k$ is the characteristic given in
\subsecref{subsec:Riemann}. We have used \thmref{thm:main}(2) and
\eqref{eq:1}.

This is a generalization of \cite[7.1]{part1}. Note that this
expression for $r=2$  coincides with the blowup formula derived from
the $u$-plane integral in \subsecref{subsec:u-plane_integral}. The
identification of the contact term follows from \eqref{eq:contact_term}.
\begin{NB}
\begin{equation*}
   T = T_{1,1}
   = \frac1{32} \left(\Lambda\frac{\partial}{\partial\Lambda}\right)^2 F_0
   = - \frac18 \Lambda\frac{\partial}{\partial\Lambda} u
   = \frac1{24} \left(\frac{du}{da}\right)^2 E_2(\tau) - \frac16 u.
\end{equation*}
\end{NB}

For a physical derivation of a higher rank generalization, see
\cite{LNS1}. 

\subsection{A reformulation}\label{subsec:reform}
We reformulate the blowup formula in the previous subsection in a form
which does not involve the limit. It also provides an interpretation
of \thmref{thm:main}(2) which does not involve the limit.
The goal is to express the formula in $H_*^T$ instead of
$H_*^{\hT}$. (Recall $\hT = \C^*\times\C^*\times T$.)

Let $j^!\colon H_*^{\hT}(M_0(r,n)) \to H_*^T(M_0(r,n))$ be the
homomorphism given by the restriction of the action. This can be
defined via the pull-back homomorphism with respect to a locally
trivial fibration with fiber $\hT/T$
\begin{equation*}
   j\colon M_0(r,n)\times_T U \to M_0(r,n)\times_{\hT} U,
\end{equation*}
where $U$ is a $\hT$-variety as in \secref{sec:BMhom}.

Recall that we have made an identification
\begin{equation*}
   \prod_n \left( H_*^{\hT}(M_0(r,n))\otimes_S \mathcal S\right)
   \Lambda^{2rn}
   \cong 
   \prod_n \mathcal S \Lambda^{2rn}
\end{equation*}
via $(\iota_0)_*$. The multiplication of $\Lambda$ in the right hand
side is identified with the push-forward homomorphism $i_{n,n+1*}$
of the natural embedding $i_{n,n+1}\colon M_0(r,n)\to
M_0(r,n+1)$. The identification follows from the commutativity of the
diagram
\begin{equation*}
\begin{CD}
   H^{\hT}_*(\mathrm{pt}) @>\iota_{0*}>> H^{\hT}_*(M_0(r,n))
\\
   @|                                     @VV{i_{n,n+1*}}V
\\
   H^{\hT}_*(\mathrm{pt}) @>\iota_{0*}>> H^{\hT}_*(M_0(r,n+1))
\end{CD}
\end{equation*}
In particular, the multiplication of $\Lambda^{2r}$ makes sense as
operators on $\prod_n H^{\hT}_*(M_0(r,n))\Lambda^{2rn}$. It does makes
sense also on $\prod_n H^{T}_*(M_0(r,n))\Lambda^{2rn}$, and two
homomorphisms commute with $j^!$.

We consider
\begin{equation*}
\begin{split}
  & \Zin(\vec{a};\Lambda)
   = \sum_{n} \Lambda^{2rn} \pi_*[M(r,n)]
   ,
\\
  & \bZin_{c_1=k}(\vec{a};\Lambda,\vec{t})
  = \sum_{n} \Lambda^{2rn}
  \widehat\pi_* \left[\exp\left(
     \sum_{p=1}^\infty t_{p}\left(\ch_{p+1}(\bE)/[C]\right)
   \right)
   \cap [\bM(r,k,n)]\right].
\end{split}
\end{equation*}
These are the formal sums of elements in $H^{T}_*(M_0(r,n))$. They are
the pull-backs of the corresponding elements in $H^{\hT}_*(M_0(r,n))$
via $j^!$.

\begin{Lemma}
Let $R = \{f(\vec{a},\ve_1,\ve_2) \mid f(\vec{a},0,0) \ne 0 \}$.
Then
$H^{\widetilde{T}}_*(M_0(r,n))_R:=
H^{\widetilde{T}}_*(M_0(r,n)) \otimes_{S(T)} S(T)_R$
is a torsion free $S(\widetilde{T})_R$-module.
\end{Lemma}

\begin{proof} 
By the localization theorem, we have
\begin{equation*}
H^{\widetilde{T}}_*(M_0(r,n))_R \cong
H^{\widetilde{T}}_*(M_0(r,n)^T)_R \cong
H^{{\Bbb C}^* \times {\Bbb C}^*}_*(M_0(r,n)^T) \otimes_{\Bbb C} 
 S(T)_R.
\end{equation*}
Since $M_0(r,n)^T = S^n \C^2$, $H^{{\Bbb C}^* \times {\Bbb
C}^*}_*(M_0(r,n)^T)$ is a torsion free ${\Bbb C}[\ve_1,\ve_2]$-module.
\end{proof}

Therefore the blowup formula in the previous section can be
restricted:
\begin{equation}\label{eq:blowupeq}
   \bZin_{c_1=k}(\vec{a};\Lambda,\vec{t})
   =
   \eqref{eq:factor}
   \times
   \Zin(\vec{a};\Lambda).
\end{equation}
This is an equality in the formal power series in $\Lambda^{2r}$ and
$\vec{t}$ with values in $H^T_*(M_0(r,n))\otimes_{S(T)}\mathcal
S(T)$. Let us emphasize again that the multiplication of
$\Lambda^{2r}$ is $i_{n,n+1*}$.

Recall $\Zin(\ve_1,\ve_2,\vec{a};\q,\vec{\tau}) =
\exp(\Fin(\ve_1,\ve_2,\vec{a};\q,\vec{\tau})/\ve_1\ve_2)$. By
\eqref{eq:pdZ} we have
\begin{equation*}
\begin{split}
   & \sum_{n=0}^\infty
    \Lambda^{2rn}
    \int_{M(r,n)}
     \prod_i \left( \ch_{\mu_i+1}(\mathcal E)/[0]\right)
 =  \left.\left(\prod_i \ve_1\ve_2 \pd{}{\tau_{\mu_i}}\right)
    \Zin(\ve_1,\ve_2,\vec{a};\Lambda^{2r},\vec{\tau})\right|_{\vec{\tau}=0}
\\
   =\; &
   \left(\prod_i
    \pd{\Fin}{\tau_{\mu_i}}(\ve_1,\ve_2,\vec{a};\Lambda^{2r},0)\right)
   \Zin(\ve_1,\ve_2,\vec{a};\Lambda^{2r},0)
\end{split}
\end{equation*}
where $[0]$ is the fundamental class of the origin. We can make a
restriction:
\begin{equation}\label{eq:simple_type}
   \sum_{n=0}^\infty
    \Lambda^{2rn}
    \pi_*\left(
     \prod_i \left(\ch_{\mu_i+1}(\mathcal E)/[0]\right)
   \cap [M(r,n)]
   \right)
     = 
     \left(\prod_i c_{\mu_i+1}\right)
   \Zin(\vec{a};\Lambda^{2r})
\end{equation}
thanks to \thmref{thm:main}(2), where we assume all $\mu_i\le r-1$.

Note that this formula explains the meaning of $c_p = \partial
F_0/\partial \tau_{p-1}$ without taking limit. It is just
multiplication of $\ch_p(\mathcal E)/[0]$. Also it formally looks like
\eqref{eq:expect_u} when $r=2$, $\mu_i = 1$.

We suppose
\begin{enumerate}
\item The factor \eqref{eq:factor} is in
$\C[c_2,\dots,c_r,\Lambda^{2r}] [[t_1,\dots,t_{r-1}]]$.
\end{enumerate}
Let us denote it by $B^{c_1=k}(\vec{c},\Lambda,\vec{t})$. Note that
this is a conjecture on theta functions. It seems that this was proved
in \cite{EGM}. But we do not quite check the detail.
Then \eqref{eq:blowupeq} becomes
\begin{equation}
\label{eq:blowupeq2}
\begin{split}
   & \sum_{n} \Lambda^{2rn}
   \widehat\pi_* \left[\exp\left(
     \sum_{p=1}^\infty t_{p}\left(\ch_{p+1}(\bE)/[C]\right)
   \right)
   \cap [\bM(r,k,n)]\right]
\\
   &\qquad\qquad\qquad\qquad\qquad =
   \sum_{n} \Lambda^{2rn} \pi_*\left(
   B^{c_1=k}(\ch(\mathcal E)/[0],\Lambda,\vec{t})\cap [M(r,n)]\right).
\end{split}
\end{equation}
We conjecture
\begin{enumerate}
\setcounter{enumi}{1}
\item \eqref{eq:blowupeq2} holds for moduli spaces $M_H(r,c_1,n)$,
$\bM_H(r,c_1+kC,n)$ for an arbitrary projective surface $X$.
\end{enumerate}

Note also that \eqref{eq:simple_type} explains the meaning
Kronheimer-Mrowka's simple type condition if we have the same formula
for an arbitrary surface:
\begin{equation*}
   \sum_{n=0}^\infty
    \Lambda^{4n}
    \pi_*\left(
     \left(\ch_{2}(\mathcal E)/[0]\right)^2
   \cap [M_H(2,c_1,n)]
   \right)
     = 
   4\Lambda^2 \sum_{n=0}^\infty
    \Lambda^{4n}
    \pi_*
   [M_H(2,c_1,n)].
\end{equation*}
It is equivalent to $c_2^2 = u_2^2 = 4\Lambda^2$. It means that the
Seiberg-Witten curve is singular.

\subsection{Rank $2$ case}
We assume $r = 2$. In this subsection we give an explicit expression
for $B^{c_1=k}(u,\Lambda,t)$ in terms of Weierstrass
$\sigma$-functions as in \cite{FS}. This is an exercise in elliptic
functions and involves no geometry.

The regularity assumption made in the previous subsections is true as
we consider the derivative with respect to $\tau_1$.

\subsubsection{Weierstrass functions}

Let $\wp(z)$ be the Weierstrass $\wp$-function
with the period ${\Bbb Z}\omega+{\Bbb Z}\omega'$, 
where $\omega'/\omega=\tau$, $\Im \tau>0$.
Then the associated elliptic curve is given by
\begin{equation*}
\begin{split}
y^2&=4(x-e_1)(x-e_2)(x-e_3)\\
&=4x^3-g_2x-g_3.
\end{split}
\end{equation*}
where $e_1=\wp(\omega/2)$,
$e_2=\wp(-\omega/2-\omega'/2)$ and
$e_3=\wp(\omega'/2)$.

Let $\sigma(z)$ be the Weierstrass $\sigma$-function.
We have the following expansion:
\begin{equation*}
\begin{split}
\wp(z) &=\sum_{n \geq 0}c_n(g_2,g_3)z^{2n-2}=
\frac{1}{z^2}+\frac{g_2}{2^2 \cdot 5}z^2+\frac{g_3}{2^2\cdot 7}z^4+\cdots\\
\sigma(z)&= \sum_{n \geq 0}c_n'(g_2,g_3)z^{2n+1}=
z-\frac{g_2}{2^4\cdot 3\cdot 5}z^5-
\frac{g_3}{2^3\cdot 3\cdot 5\cdot 7}z^7-\cdots
\end{split}
\end{equation*}
where $c_n, c_n'$ are weighted homogeneous polynomials of
degree $n$ with
$\deg g_2=2$ and $\deg g_3=3$.
Let $\sigma_i(z)$, $i=1,2,3$ be three more sigma functions 
associated to $e_i$. 
We assign $\deg e_i:=1$.
Since $\sigma_i(z)^2=\sigma(z)^2(\wp(z)-e_i)$,
$\sigma_i(z)$ also has an expansion
\begin{equation}
\sigma_i(z)=1-\frac{e_i}{2}z^2+\sum_{n \geq 2}c_n''(e_i,g_2,g_3)z^{2n}
\end{equation}
where $c_n''$ are weighted homogeneous polynomials of
degree $n$.

In terms of modular forms, we have the following expressions:
\begin{equation}\label{formula:1}
\begin{split}
\sigma(\omega z)&= \omega e^{\omega \eta z^2}
\frac{\theta_{11}(z|\tau)}{\theta_{11}'(0|\tau)}=
-\frac{\omega}{\pi \theta_{10}(0|\tau)\theta_{00}(0|\tau)}
e^{\omega \eta z^2}\frac{\theta_{11}(z|\tau)}{\theta_{01}(0|\tau)}
\\
\sigma_3(\omega z)&=e^{\omega \eta z^2}
\frac{\theta_{01}(z|\tau)}{\theta_{01}(0|\tau)},
\end{split}
\end{equation}
where $\eta=\zeta(\omega/2)$ is given by
\begin{equation*}
\eta \omega=\frac{\pi^2}{6}E_2(\tau).
\end{equation*}

We also have
\begin{equation}\label{formula:2}
\begin{split}
e_1=&\frac{1}{3}\left(\frac{\pi}{\omega}\right)^2
(\theta_{00}^4+\theta_{01}^4),\\
e_2=&\frac{1}{3}\left(\frac{\pi}{\omega}\right)^2
(\theta_{10}^4-\theta_{01}^4),\\
e_3=&-\frac{1}{3}\left(\frac{\pi}{\omega}\right)^2
(\theta_{10}^4+\theta_{00}^4).
\end{split}
\end{equation}

\begin{equation}\label{formula:3}
\begin{split}
g_2 &=-4
\left(\frac{\pi}{\omega}\right)^4\left(-\frac{1}{3}(\theta_{10}^4+\theta_{00}^4)^2+
(\theta_{10}^4\theta_{00}^4)\right),\\
g_3 &=\frac{1}{27}\left(\frac{\pi}{\omega}\right)^6
(\theta_{10}^4+\theta_{00}^4)\left(8(\theta_{10}^4+\theta_{00}^4)^2-36
(\theta_{10}^4\theta_{00}^4)\right).
\end{split}
\end{equation}

\subsubsection{}
We write $u = u_2$, $a = a_2$ as in \subsecref{subsec:rank2}.
We set $\omega=\pi \theta_{00} \theta_{10}/\Lambda$. Then
\begin{equation*}
\begin{split}
g_2 &=4
\left(\frac{1}{3}u^2-\Lambda^4\right)\\
g_3 &=-\frac{1}{27}
u\left(8u^2-36\Lambda^4\right)\\
e_3&=\frac{u}{3}.
\end{split}
\end{equation*}
Hence we get
\begin{equation}\label{eq:FS}
\begin{split}
e^{-T_{1,1}t^2}
\frac{\theta_{01}(\frac{\sqrt{-1}}{2 \pi}\frac{du}{da}t|\tau)}
{\theta_{01}(0|\tau)}
&=e^{\frac{u}{6}t^2}\sigma_3(t),\\
e^{-T_{1,1}t^2}
\frac{\theta_{11}(\frac{\sqrt{-1}}{2 \pi}\frac{du}{da}t|\tau)}
{\theta_{01}(0|\tau)}
&=e^{\frac{u}{6}t^2}\sigma(t)\Lambda.
\end{split}
\end{equation}
This checks the conjecture~(1) in the previous subsection.

Since $\ch_2(\mathcal E) = -c_2(\mathcal E)$, we put $x = -u$. Then
$x$ corresponds to the insertion of the point class $\mu(p) =
c_2(\mathcal E)/[0]$ in \eqref{eq:blowupeq2}. The above \eqref{eq:FS}
exactly coincides with the functions $B(x,t)$, $S(x,t)$ appeared in
Fintushel-Stern's blowup formula for Donaldson invariants \cite{FS}.

This checks the conjecture~(2) in a weak sense, i.e., if we multiply
products of $\mu(S)$ and integrate, then the equality
holds. Conversely if we can prove (2) directly, it gives a new proof
of Fintushel-Stern's blowup formula. The proof of (2) probably
requires more detailed study of the map $\widehat\pi$.

\section{Gravitational corrections}

As we mentioned in Introduction, Nekrasov asserts that higher order
terms in \eqref{eq:expansion} are {\it gravitational corrections\/} to
the gauge theory \cite[\S4]{Nek}.
In some cases these are some known quantities, which are really
related genus $g$ curves, e.g., Gromov-Witten invariants with domain
genus $g$.
We are still far away from verifying this conjecture in full
generality. But we have some nontrivial examples, which we review in
this section.

\subsection{Genus $1$ part}\label{subsec:genus 1}
The result of this subsection is based on discussions with
N.~Nekrasov.

We determine the coefficients $G$, $F_1$ in the expansion
\eqref{eq:expansion} for $\vec{\tau} = 0$. These terms are considered
as genus $1$ gravitational corrections as we said. Since we are only
interested in $\vec{\tau} = 0$ case, we omit $\vec{\tau}$ from the
notation.

We consider the blowup equation \eqref{eq:blowup_formula3} for $c_1 =
kC$ ($k\neq 0$) with $\vec{\tau} = \vec{t} = 0$:
\begin{equation*}
   0 = \sum_{\{\vec{k}\} = -\frac{k}r}
   Z(\ve_1,\ve_2-\ve_1,\vec{a}+\ve_1\vec{k};\q)
        Z(\ve_1-\ve_2,\ve_2,\vec{a}+\ve_2\vec{k};\q).
\end{equation*}
As in the derivation of \eqref{eq:limit} we have
{\allowdisplaybreaks
\begin{equation*}
\begin{split}
  0 &= \sum_{\{\vec{k}\} = -\frac{k}r}\!\!
  \begin{aligned}[t]
  & 
  \exp\Biggl[
  - \frac{\partial^2 \Fz}{\partial a^l\partial a^m}
  \frac{k^l k^m}2
  + \frac{\partial H}{\partial a^l}k^l
    \\
  & \qquad + 
  (\ve_1+\ve_2)\left\{
    -\frac{\partial^3 \Fz}{\partial a^l\partial a^m\partial a^n}
  \frac{k^l k^m k^n}{3!}\right.
  + \left.
    \frac{\partial (G+ \Fo)}{\partial a^l} k^l
  \right\} + \cdots
  \Biggr]
  \end{aligned}
\\
  &= \sum_{\{\vec{k}\} = -\frac{k}r} \!\!
  \begin{aligned}[t]
  &
  \sqrt{-1}^{2\langle\vec{k},\rho\rangle}
  \exp\Biggl[
  \pi\sqrt{-1}\tau_{lm} k^l k^m
    \\
  & \qquad + 
  (\ve_1+\ve_2)\left\{
    -\frac{\partial^3 \Fz}{\partial a^l\partial a^m\partial a^n}
  \frac{k^l k^m k^n}{3!}\right.
  + \left.
    \frac{\partial (G+\Fo)}{\partial a^l} k^l
  \right\} + \cdots
  \Biggr],
  \end{aligned}
\end{split}
\end{equation*}
where} $\tau_{lm}$ is the period of the Seiberg-Witten curve as
before \eqref{eq:period}. 
 
Setting $\ve_1 = \ve_2 = 0$, we get
\begin{equation*}
  0 = \sum_{\{\vec{k}\} = -\frac{k}r}\!\!
  \sqrt{-1}^{2\langle\vec{k},\rho\rangle}
  \exp\left(\pi\sqrt{-1}\tau_{lm} k^l k^m\right) = 
  \Theta_{E_k}
  (\vec{0}|\tau)
  .
\end{equation*}

\begin{NB}
This equality is trivial when $r=2$, $k=1$.
\begin{equation*}
  \sum_{\{\vec{k}\} = -\frac12} 
  \sqrt{-1}^{2\langle\vec{k},\rho\rangle}
  \exp\left[ \pi\sqrt{-1}\tau k^2 \right]
  = \sum_{n\in\Z} \sqrt{-1}^{2(n+\frac12)}
  \exp\left[ \pi\sqrt{-1}\tau (n+\frac12)^2 \right]
\end{equation*}
vanishes because of the symmetry $n+\frac12 \leftrightarrow -(n+\frac12)
= (-1-n)+\frac12$. 
\end{NB}

Next we take the coefficient of $\ve_1+\ve_2$ in the above to get
\begin{equation*}
  0 =
\sum_{\{\vec{k}\} = -\frac{k}r} \!\!
  \begin{aligned}[t]
  &
  \sqrt{-1}^{2\langle\vec{k},\rho\rangle}
  \exp\Biggl[
  \pi\sqrt{-1}\tau_{lm} k^l k^m \Biggr]
    \\
  & \qquad \times
  \left\{
    2\pi\sqrt{-1}\pd{\tau_{mn}}{a^l}\frac{k^l k^m k^n}{3!}\right.
  + \left.
    \frac{\partial (G+\Fo)}{\partial a^l} k^l
  \right\},
  \end{aligned}
\end{equation*}
i.e.,
\begin{equation*}
   \sum_l
   \frac{\partial (G+\Fo)}{\partial a^l}
    \pd{}{\xi^l}\Theta_{E_k}
   (\vec{0}|\tau)
   +
   \frac13 \frac{\partial^2}{\partial a^l\partial\xi^l}
   \Theta_{E_k}
   (\vec{0}|\tau)
   = 0.
\end{equation*}
We believe this equation with $k=1,\dots,r-1$ determine
$\frac{\partial (G+\Fo)}{\partial a^l}$ for $l=2,\dots,r$. But we do
not know the required identities for the theta functions, as far as
the authors are concerned. So we assume $r=2$ from now. Then the
equation is
\begin{equation}\label{eq:3}
   \frac{\partial (G + \Fo)}{\partial a}
   = - \frac1{3} \frac{\partial}{\partial a} \log\left(
   \pd{}{\xi}
   \Theta_{E_1}
   (\vec{0}|\tau)
   \right).
\end{equation}

We now switch to the notation in \subsecref{subsec:g=1}. From
(\ref{eq:1}, \ref{eq:3}), we get
\begin{equation*}
    \exp\left(G-\Fo\right) = \theta_{01}(0,\tau), \quad
    \exp\left(G+\Fo\right) = C\theta_{11}'(0,\tau)^{-1/3}
\end{equation*}
for some constant $C$ independent of $a$. Therefore
\begin{equation*}
   \exp(2\Fo) = C\theta_{11}'(0,\tau)^{-1/3}\theta_{01}(0,\tau)^{-1}.
\end{equation*}
By Jacobi's triple product identity (see e.g., \cite[Chap.~I,
\S14]{Mumford}) we have
\begin{equation*}
   \theta_{01}(0,\tau) = \prod_{d=1}^\infty
   \left[(1-q^{2d})(1-q^{2d-1})^2\right].
\end{equation*}
By Jacobi's derivative formula (see \cite[Chap.~I, \S13]{Mumford}) we have
\begin{equation*}
    \theta_{11}'(0,\tau) = -2\pi q^{\frac14}
    \prod_{d=1}^\infty (1-q^{2d})^3.
\end{equation*}
Therefore we get
\begin{equation*}
   \exp \Fo = C' q^{-\frac1{24}} \prod_{d=1}^\infty (1-q^d)^{-1}
   = \frac{C'}{\eta(\frac\tau{2})}
\end{equation*}
for some constant $C'$ independent of $a$. A priori, $C'$ may depend
on $\q$ (or $\Lambda$), but in fact, it does not as follows. Let us
define degrees of variables by
\begin{equation*}
   \deg \ve_1 = \deg \ve_2 = \deg a_\alpha = 1, \qquad
   \deg \q = 2r \quad(\deg \Lambda = 1).
\end{equation*}
(This definition applies for arbitrary $r$, not necessarily $2$.)  By
the definition, $Z$ has degree $0$. Therefore $\Fz$ has degree $2$,
while $G$ and $\Fo$ have degree $0$. Then $\tau$ has degree $0$, and
hence so is $\eta(\frac\tau{2})$. Therefore $C'$ has degree $0$. Since
it is independent of $a$, it means that it is also independent of $\q$
(or $\Lambda$). Therefore $C'$ can be computed by studying the
expansions of $\Fo$ and $\eta$ in $a/\Lambda$:
\begin{equation*}
  \begin{split}
    & \Fo = \frac1{12} \left\{ \log\left(\frac{2a}{\Lambda}\right)
    + \log\left(\frac{-2a}{\Lambda}\right) \right\} + \cdots
  = \frac16 \log\left(\frac{2\sqrt{-1}a}{\Lambda}\right) + \cdots,
\\
    & \log \frac1{\eta(\frac{\tau}2)}
    = - \frac1{24} \log q - \sum_{d=1}^\infty \log (1 - q^d)
    = - \frac1{24}\pi\sqrt{-1}\tau + \cdots.
  \end{split}
\end{equation*}
Furthermore, we have
\begin{equation*}
   \tau = \frac{\sqrt{-1}}\pi 4\log\left(\frac{2\sqrt{-1}a}{\Lambda}
   \right) + \cdots.
\end{equation*}
Therefore we get $C' = 1$, and hence
\begin{equation*}
   \Fo = - \log \eta(\frac\tau2), \qquad
   G = \log \left[
     q^{-\frac1{24}} \prod_{d=1}^\infty (1-q^{2d-1})
     \right].
\end{equation*}

It is better to make the following combination:
\begin{equation*}
   F_1 = A - \frac23 B, \qquad
   G = \frac13 B.
\end{equation*}
Then
\begin{equation*}
\begin{split}
   \exp A &= \exp \left(F_1 + 2G\right)
   = q^{-\frac18}\prod_{d=1}^\infty \frac{1-q^{2d-1}}{1-q^{2d}}
   = \left(-\frac{2\pi\theta_{01}}
     {\theta_{11}'}\right)^{\frac12}
   = \left(\frac2
     {\theta_{00}\theta_{10}}\right)^{\frac12},
\\
   \exp B & = \exp\left(3 G\right)
   = q^{-\frac18} \prod_{d=1}^\infty (1 - q^{2d-1})^3
   = \left(\frac{-2\pi \theta_{01}^3}{\theta_{11}'}\right)^{\frac12}
   = \left(\frac{2\theta_{01}^2}
     {\theta_{00}\theta_{10}}\right)^{\frac12}
\\
   &= \left(\frac{4\theta_{01}^4}
     {\theta_{00}^2\theta_{10}^2}\right)^{\frac14}
   = \left(\frac{4 (\theta_{00}^4 -
   \theta_{10}^4)}{\theta_{00}^2\theta_{10}^2}\right)^{\frac14}
   .
\end{split}
\end{equation*}
Comparing with \eqref{eq:rank2}, we find
\begin{equation}\label{eq:AB}
   \exp A = \left(\frac{\sqrt{-1}}\Lambda\frac{du}{da}\right)^{\frac12},
\qquad
   \exp B = \left(\frac{4 (u^2 - 4\Lambda^4)}{\Lambda^4}\right)^{\frac18}.
\end{equation}
Note that the last expression is given by the quantum discriminant
\eqref{eq:desc}:
\begin{equation*}
   \Delta = 2^{12}\Lambda^8(u^2 - 4\Lambda^4).
\end{equation*}
Therefore
\begin{equation*}
   \ve_1 \ve_2 F_1 + (\ve_1 + \ve_2)^2 G
   = \ve_1\ve_2 \log \left(\frac{\sqrt{-1}}\Lambda\frac{du}{da}\right)^{\frac12}
   +
   \frac{\ve_1^2+\ve_2^2}3\log
   \left(\frac{\Delta}{2^{10}\Lambda^{12}}\right)^{\frac18}.
\end{equation*}
Comparing with \eqref{eq:u-plane_int}, this suggests the following
formula for the equivariant Euler number and signature for $\C^2$:
\begin{equation*}
   \chi(\C^2) = \ve_1\ve_2, \qquad
   \sigma(\C^2) = \frac{\ve_1^2+\ve_2^2}3.
\end{equation*}
This is natural from the following formal computation:
\begin{equation*}
\begin{gathered}
   \chi(\C^2) = c_2(\C^2), \qquad
   \sigma(\C^2) = \frac13 \left(c_1(\C^2)^2 - 2 c_2(\C^2)\right),
\\
   c_1(\C^2) = \ve_1 + \ve_2, \quad c_2(\C^2) = \ve_1\ve_2.
\end{gathered}
\end{equation*}

Nekrasov conjectures that \eqref{eq:AB} holds higher rank case also if
we replace $\frac{du}{da}$ by $\det(\pd{u_p}{a^i})$.

\subsection{Coordinate rings of symmetric products}\label{subsec:coord}
The next example is related to the {\it perturbative part\/} of the
$K$-theory version of the partition function in \secref{sec:Nek}. (See
\cite[A.0.3]{NO}.) It fits with geometric engineering quite well (see
below), but it looks like an accident if we understand it as a purely
mathematical statement. The authors learned the result from
\cite[\S A.0.3]{NO}.

Let us consider the $n$th symmetric product $S^n(\C^2)$ of the affine
plane $\C^2$. We define an action of the two torus $T^2 =
\C^*\times\C^*$ on $\C^2$ by
\begin{equation*}
   (x, y) \longmapsto (t_1 x, t_2 y), \qquad (t_1,t_2)\in T^2.
\end{equation*}
We also have an induced action on $S^n(\C^2)$.

The coordinate ring $H^0(S^n(\C^2), \shfO)$, that is the ring of
polynomial functions on $S^n(\C^2)$, is a $T^2$-module. We consider
its character
\begin{equation*}
  \ch H^0(S^n(\C^2), \shfO)
  = \sum_{m,n\ge 0} t_1^m t_2^n \dim H^0(S^n(\C^2), \shfO)_{m,n},
\end{equation*}
where
\begin{equation*}
   H^0(S^n(\C^2), \shfO)_{m,n} = \left.\left\{ f\in H^0(S^n(\C^2),
   \shfO)\, \right|
   (t_1, t_2)\cdot f = t_1^m t_2^n f \right\}
\end{equation*}
is a simultaneous eigenspace, i.e., a weight space. The character is
called {\it Hilbert series\/} sometimes. It is standard in algebraic
geometry to show that
\begin{enumerate}
\item each weight space is finite-dimensional, and hence the character
is well-defined as a formal sum,
\item the character is a rational function in $t_1^\pm$, $t_2^\pm$. 
\end{enumerate}
In fact, in this case, an explicit answer can be written down:
\begin{Proposition}
The generating function of the character is given by
\begin{equation}\label{eq:chHilb}
   \sum_{n=0}^\infty \q^n \ch H^0(S^n(\C^2), \shfO)
   = \exp\left(\sum_{d=1}^\infty \frac{\q^d}{d(1-t_1^d)(1-t_2^d)}\right).
\end{equation}
\end{Proposition}
We refer \cite[\S3]{part1} for the proof. But we recommend the reader
to write down the proof by himself/herself since it is a nice exercise
on a treatment of generating functions.

The generating function \eqref{eq:chHilb} is the rank $1$ version of
Nekrasov's deformed partition function (for the $K$-theory
version). Let us expand it into a formal power series in $\hbar$,
after putting $t_1 = \exp \hbar$, $t_2 = \exp(-\hbar)$:
\begin{equation}\label{eq:expand_ch}
\begin{split}
   & \log\left(\sum_{n=0}^\infty \q^n \ch H^0(S^n(\C^2), \shfO)\right)
    \left|_{\substack{t_1 = e^{\hbar}\\ t_2 =
      e^{-\hbar}}}\right.
   =
  \sum_{d=1}^\infty \frac{\q^d}
  {d(1-e^{\hbar d})(1-e^{-\hbar d})}
\\
  =\; &
   - 
   \sum_{d=1}^\infty \frac{\q^d}{d^3} \hbar^{-2}
   - \frac1{12}\log\left(1-\q\right) \hbar^0
   + \sum_{g \ge 2} \frac{B_{2g}}{2g(2g-2)!}
     \sum_{d=1}^\infty d^{2g-3} \q^d \hbar^{2g-2},
\end{split}
\end{equation}
where $B_{2g}$ is the $2g$th Bernoulli number as in
\secref{sec:perturb}. The series start with $\hbar^{-2}$ and have only
even powers of $\hbar$. In the next subsection we will see that
\eqref{eq:expand_ch} is equal to the generating function of certain
Gromov-Witten invariants. Then $g$ will be identified with the genus
of the domain curve, and hence $2-2g$ with the Euler number.

\subsection{Gromov-Witten invariants of the resolved
conifold}\label{subsec:conifold}

Let $X$ be the total space of the rank $2$ vector bundle $E =
\shfO(-1)\oplus \shfO(-1)$ over $\proj^1$. This space is called the
{\it resolved conifold\/} in physics.
The local Gromov-Witten invariants for target $X$ is defined as
follows: Let $M_{g,n}(\proj^1,d)$ be the moduli space of stable
maps for target $\proj^1$ from a genus $g$ curve with $n$ marked
points with degree $d$. We consider $d > 0$ case only, that is the
stable map is not constant. We have the diagram
\begin{equation*}
   M_{g,0}(\proj^1,d) \xleftarrow{\rm forget}
   M_{g,1}(\proj^1,d) \xrightarrow{\rm eval} \proj^1,
\end{equation*}
where $\rm forget$ is the map given by forgetting the marked point,
and ${\rm eval}$ is the map given by taking the image of the marked
point under the stable map. We consider a vector bundle $\mathbb E =
R^1{\rm forget}_*{\rm eval}^*E$. Let $c_{\rm top}(\mathbb E)$ be its
top Chern class. Then the local Gromov-Witten invariant is defined by
\begin{equation*}
   C(g,d) = 
   \int_{M_{g,0}(\proj^1,d)^{\operatorname{vir}}} c_{\rm top}(\mathbb E),
\end{equation*}
where $M_{g,0}(\proj^1,d)^{\operatorname{vir}}$ is the virtual
fundamental class. This is a rational number. This is a local
contribution to the {\it global\/} Gromov-Witten invariant of a
Calabi-Yau 3-fold of multiple covers of a fixed rational curve with
the normal bundle $E = \shfO(-1)\oplus \shfO(-1)$. Thus the formula
for $C(g,d)$ is important in the Gromov-Witten theory. The genus $g=0$
case is known as the Aspinwall-Morrison formula. The complete answer
is given by
\begin{Theorem}[\cite{AM,M,Vo} for $g=0$, \cite{GrP} for $g=1$,
\cite{FP} for $g\ge 2$]
\begin{equation*}
   C(0,d) = \frac1{d^3}, \quad
   C(1,d) = \frac1{12d}, \quad
   C(g,d) = \frac{(-1)^{g-1}B_{2g}}{2g(2g-2)!} d^{2g-3}.
\end{equation*}
\end{Theorem}

Comparing with \eqref{eq:expand_ch}, we get
\begin{equation}\label{eq:GWgen}
   \log\left(\sum_{n=0}^\infty \q^n \ch H^0(S^n(\C^2), \shfO)\right)
    \left|_{\substack{t_1 = e^{\hbar}\\ t_2 =
      e^{-\hbar}}}\right.
    =
    \sum_{d=1}^\infty \sum_{g=0}^\infty
    C(g,d) \q^d (i\hbar)^{2g-2}.
\end{equation}
This is our first example of the assertion that the gauge theory
partition function is identified with generating functions of
Gromov-Witten invariants.

It is worthwhile mentioning that the exponential of the right hand
side is the generating function of Gromov-Witten invariants whose
domain curves are not necessarily {\it connected}. This does not make
sense in the gauge theory side, but somehow related to the recursive
structure among the symmetric products $S^n \C^2$ for various $n$.


\subsection{The $r=1$ case and Gromov-Witten invariants for $\proj^1$}
In the main body of the paper, the case $r=1$ was excluded as the
Seiberg-Witten geometry does not make sense. However Nekrasov's
partition function does make sense $r=1$ also. The $K$-theory version
with $\vec{\tau} = 0$ is what we already saw in
\subsecref{subsec:coord}. This is because the moduli space $M(1,n)$ is
nothing but the Hilbert scheme of points, and it is known that the
higher direct image sheaves for $M(1,n)\to S^n\C^2$ vanish and
$\ch H^0(S^n\C^2, \shfO)$ can be given by Atiyah-Bott formula for
$M(1,n)$. (See \cite[\S3]{part1} for more detail.)
When $\vec{\tau}\neq 0$, it was studied in \cite{Haiman}. (His main
result is a positivity property, which is not studied in this
paper. But it is natural to conjecture a similar property for higher
rank cases also.)

The homology version of the partition function, i.e., our
$Z(\ve_1,\ve_2;\q,\vec{\tau})$, has the presentation by the Fock space
when it is restricted to $\ve_2 = -\ve_1$. Then comparing with the
presentation for the Gromov-Witten invariants for $\proj^1$ \cite{OP1},
one gets
\begin{equation*}
   \log Z(\ve_1,-\ve_1;\q,\vec{\tau})
   = 
   \;
   \begin{minipage}[c]{.65\textwidth}
   the generating function of the Gromov-Witten invariants for
   $\proj^1$,
   \end{minipage}
\end{equation*}
where $\ve_1$ is mapped to an indeterminate for the domain genus, and
$\vec{\tau}$ to those for gravitational descendants. This remarkable
observation was done by \cite{LQW} and \cite{Nek2} independently. We
refer the precise statement and the proof to the original papers.

\subsection{Geometric Engineering}\label{subsec:engineer}

The geometric engineering of Katz-Klemm-Vafa \cite{KKV} realizes
$4$-dimen\-sional $\mathcal N=2$ supersymmetric gauge theories as
limits of type IIA string theory compactified on certain noncompact
Calabi-Yau $3$-folds. Mathematically it poses the following conjecture:
\begin{equation}\label{eq:Math_Engineer}
\fbox{
\begin{minipage}{.85\textwidth}
Partition functions in $4$-dimensional $\mathcal N=2$ supersymmetric
gauge theories are equal to limits of generating functions of local
Gromov-Witten invariants for certain noncompact Calabi-Yau $3$-folds.
\end{minipage}}
\end{equation}
The noncompact Calabi-Yau $3$-fold is chosen according to the gauge
theory to realize. Typically it is an ALE space fibration over
$\proj^1$. Recall that an ALE space is the minimal resolution of a
simple singularity and contains a configuration of $\proj^1$'s
intersecting as one of $ADE$ Dynkin diagrams. The group of the gauge
theory is the corresponding $ADE$ group.

This statement is rather striking since it seems difficult to compare
two types of moduli spaces directly, i.e., moduli spaces of stable
maps and instanton moduli spaces. Moreover, we must sum up over
degrees for Gromov-Witten invariants as we will see below.

For simplicity, we restrict ourselves to rank $2$ case. (In higher
rank cases, the above naive definition of the local Gromov-Witten
invariants must be modified as the base space become singular.)
As we have already seen in \subsecref{subsec:conifold}, the $K$-theory
version of the partition function is more natural here as we do not
need to take a limit. We define
\begin{equation*}
\begin{gathered}
  Z_K(\hbar,a;\beta)
   = \exp\left[-
      \gamma^K_{\hbar}(2a;\beta) - \gamma^K_{\hbar}(-2a;\beta)
      \right]\, \Zin_K(\hbar,{a};\beta),
\\
  \Zin_K(\hbar,{a};\beta) = 
  \sum_{n=0}^\infty \beta^{2rn} \sum_i (-1)^i
  \ch H^i(M(r,n), \shfO),
\\
  t_1 = e^{\beta\hbar}, \quad t_2 = e^{-\beta\hbar}, \quad
  e_1 = e^{-\beta a}, \quad e_2 = e^{\beta a},
\end{gathered}
\end{equation*}
where $\ch$ is the character of $\hT$-module and the one-dimensional
$\hT$-modules (in the notation \ref{not:module}) are replaced as
indicated. The characters of the cohomologies have expressions in
terms of Young diagrams $\vec{Y}$ by the localization formula for the
$K$-theory:
\begin{equation*}
   \Zin_K(\hbar,{a};\beta)
   = \sum_{\vec{Y}} \beta^{2r|\vec{Y}|}
\end{equation*}
(See \cite[\S3]{part1} for detail.) The perturbative term
is given by
\begin{equation*}
  \gamma^K_{\hbar}(x;\beta) =
  \sum_{d=1}^\infty \frac{e^{-\beta d x}}
  {d(1-e^{\hbar d})(1-e^{-\hbar d})}.
\end{equation*}
This is equal to the one in \cite[\S A.0.3]{NO} up to polynomials in
$a$.

Note that we do not include higher Casimir operators. This is not for
the brevity. We (at least the authors) do not know what are
counterparts in the Gromov-Witten theory.

In the rank $2$ case, the noncompact Calabi-Yau $3$-fold is supposed
to be the canonical bundle $K_{\mathbf F_n}$ of the Hirzebruch surface
$\mathbf F_n$. And it is conjectured that the Gromov-Witten invariants
are {\it essentially\/} independent of $n$. So we further  restrict to
the case $n=0$, i.e., the $K_{\mathbf F_0}$ of $\mathbf F_0 =
\proj^1\times\proj^1$. (For the actual calculation, it is necessary to
assume $\mathbf F_n$ is a toric variety in order to apply the
localization technique.)

We define the local Gromov-Witten invariants as in the case of the
resolved conifold. The degree of maps is a pair of integers $(n,d)$
corresponding to the base and fiber respectively. (Although we have a
symmetry exchanging two factors, we break it and consider one is base
and the other is fiber. This is automatic for other $\mathbf F_n$.)
Let us introduce two parameters $q_b$, $q_f$ respectively.

Now the geometric engineering asserts that $\log Z_K(\hbar,a;\beta)$
is equal to the generating function of Gromov-Witten invariants, under
a suitable identification of parameters. The parameter $\hbar$ should
count the genus of domain curves as in
\subsecref{subsec:conifold}.

Next we match the degree $n$ for the base with the instanton number
$n$. Although this identification is quite natural, we do not have any
mathematically rigorous justification of this statement. Anyway we
should identify $q_b$ with $\beta^{2r}$. In fact, the analysis below
gives the exact answer:
\begin{equation*}
   q_b = \left(\frac{\beta}2\right)^{2r}.
\end{equation*}

Let us consider the Gromov-Witten invariants for $n=0$. In this case,
we only have multiple covers of the fiber. It is given by
\eqref{eq:GWgen} multiplied by $-2$. Since $n=0$ means zero instanton
number in the gauge theory side, it should be equal to the
perturbation term of the gauge theory, namely:
\begin{equation*}
   - \gamma_{\hbar}^K(2a|\beta) - \gamma_{\hbar}^K(-2a|\beta)
 \overset{?}{=}
   -2 \sum_{d=1}^\infty \frac{q_f^d}
   {d(1-e^{\hbar d})(1-e^{-\hbar d})}
\end{equation*}
Since the left hand side is equal to $-2\gamma_{\hbar}^K(2a|\beta)$ up
to a polynomial in $a$, thanks to the inversion formula for the
polylogarithms, this (up to a polynomial in $a$) follows from what we
observed in \subsecref{subsec:conifold} when we equate the parameters
as
\begin{equation*}
   q_f = e^{-2\beta a}.
\end{equation*}
Note that for this identification, we must sum up the Gromov-Witten
invariants for various degrees on fibers (and various genus). This is
also true for $n$-instanton corrections.

For $n > 0$, the genus $0$ Gromov-Witten invariants were calculated
using the local mirror symmetry for $X = K_{\mathbf F_0}$ in
\cite{KKV}. And the limit of their generating function was identified
with the Seiberg-Witten prepotential $F_0$. In fact, they identify the
limit of the local mirror of $X$ with the Seiberg-Witten curve. Note
that genus $0$ case is enough to identify the parameters.

Recently Iqbal+Kashani-Poor identify $\log Z_K$ with the generating
function of all genus Gromov-Witten invariants by using the large $N$
duality \cite{IK}. In fact, they identify the expression of $\log Z_K$
via Young diagrams $\vec{Y}$ with Jones-Witten invariants for the Hopf
link. (They assume certain combinatorial identities which are proved
in more recent paper \cite{EK}.)

\appendix
\section{The root system of type $A_{r-1}$}\label{sec:root}

Let $Q$ be the coroot lattice of type $A_{r-1}$:
\begin{equation*}
   Q = \left.\left\{ \vec{k} = (k_1,\dots,k_r)\in\Z^r\,\right|
     \textstyle \sum_\alpha k_\alpha = 0  \right\}.
\end{equation*}
We take simple coroots
\begin{equation*}
   \alpha_i^\vee = (0,\dots, 0, \overset{i}{1},
   \overset{i+1}{-1},0,\dots,0),
   \qquad (i=1,\dots,r-1).
\end{equation*}
We can write
\begin{equation*}
   Q \ni \vec{k} = \sum_i k^i \alpha_i^\vee.
\end{equation*}
For a given $k\in\Z$, elements $\vec{k}\in \Z^r$ with $\sum_\alpha
   k_\alpha = k$ are identified
\begin{equation*}
   \left\{ \vec{l} = (l_1,\dots,l_r)\in\Q^r \left|\,
   \textstyle  \sum_\alpha l_\alpha = 0,
  \forall \alpha\; l_\alpha \equiv -\frac{k}r \mod \Z
  \right\}\right..   
\end{equation*}
This is a subset of the coweight lattice
\(
  P = \{ \vec{l} = (l_1,\dots,l_r)\in\Q^r \mid \sum_\alpha l_\alpha = 0,
  \exists k\in\Z\; \forall \alpha\; l_\alpha \equiv -\frac{k}r \mod \Z
  \}.
\)
There exists a homomorphism $P\to \Z/r\Z$ by taking the fractional
part of $l_\alpha$. It can be identified with the natural quotient
homomorphism $P\to P/Q$. We denote it by $\vec{l}\mapsto \{ \vec{l}\}$.
Hereafter we identify $\vec{l}$ with $\vec{k}$ and denote both by
$\vec{k}$. We write $\vec{k} = \sum_i k^i \alpha_i^\vee$ in either
case $k=0$, $\neq 0$. But $k^i$ may be rational in the latter case.
Let $(\ ,\ )$ be the standard inner product on $Q$. The
Killing form $B_{\SU(r)}$ of $\SU(r)$ satisfies $B_{\SU(r)} = 2r(\ ,\ 
)$.
The following formulas are useful later:
\begin{equation}
\begin{gathered}
   \frac1{2r}   
   \sum_{\alpha, \beta} (k_\alpha - k_\beta)(a_\alpha - a_\beta)
   =
   (\vec{k},\vec{a})
   = \sum_{ij} C_{ij} a^i k^j,
\\
    \frac1{2r}\sum_{\alpha, \beta} (k_\alpha - k_\beta)^2
    =
    (\vec{k}, \vec{k})
    = \sum_{i,j} C_{ij} k^i k^j,
\\
  \sum_{\alpha < \beta} \frac{k_\alpha-k_\beta}2
  =
  \langle\vec{k},\rho\rangle
  = \sum_i k^i.
\end{gathered}
\label{eq:useful}\end{equation}
Here $C_{ij}$ is the Cartan matrix, and
$\rho$ is the half of the sum of positive roots, as usual.

\section{Theta functions}

We give definitions and some properties of Riemann theta functions.

\subsection{Riemann Theta functions}\label{subsec:Riemann}

Let $Q = \Z^g$. Let $\tau = (\tau_{\alpha\beta})$ be a symmetric
$g\times g$ complex matrix whose imaginary part is positive
definite. For $\vec{\mu}, \vec{\nu}\in \C^g$, we define the theta
function with characteristic
\( \left[
\begin{smallmatrix}
 \vec{\mu} \\ \vec{\nu}
\end{smallmatrix}
\right]
\)
by
\begin{multline*}
\Theta
\begin{bmatrix}
 \vec{\mu} \\ \vec{\nu}
\end{bmatrix}
(\vec{\xi}|\tau)
\\
= 
\sum_{\vec{k}\in Q} \exp\left(
     \pi\sqrt{-1} \sum_{\alpha,\beta}\tau_{\alpha\beta}
     (k_\alpha + \mu_\alpha) (k_\beta + \mu_\beta)
     + 2\pi\sqrt{-1}\sum_\alpha
     (k_\alpha+\mu_\alpha)(\xi_\alpha + \nu_\alpha)\right).
\end{multline*}
When
\(
\left[
\begin{smallmatrix}
 \vec{\mu} \\ \vec{\nu}
\end{smallmatrix}
\right]
= 
\left[
\begin{smallmatrix}
 0 \\ 0
\end{smallmatrix}
\right],
\)
we simply denote it by $\Theta$. We have
\begin{equation*}
   \Theta
\begin{bmatrix}
 \vec{\mu} \\ \vec{\nu}
\end{bmatrix}
(\vec{\xi}|\tau)
   = \exp\left(\pi\sqrt{-1}\, {}^t\!\vec{\mu}\tau\vec{\mu}
     + 2\pi\sqrt{-1}\,{}^t\!\vec{\mu}(\vec{\xi}+\vec{\nu})\right)
   \Theta(\vec{\xi}+\tau\vec{\mu}+\vec{\nu}|\tau).
\end{equation*}
The theta function is quasi-periodic with respect to the lattice
$\Z^g\oplus \tau\Z^g$. It satisfies the heat equation:
\begin{equation*}
   \frac{\partial^2}{\partial \xi_\alpha\partial \xi_\beta}
     \Theta
\begin{bmatrix}
 \vec{\mu} \\ \vec{\nu}
\end{bmatrix}
(\vec{\xi}|\tau)
   = 4\pi\sqrt{-1} \pd{}{\tau_{\alpha\beta}}\Theta
\begin{bmatrix}
 \vec{\mu} \\ \vec{\nu}
\end{bmatrix}
(\vec{\xi}|\tau).
\end{equation*}

When $\vec{\mu},\vec{\nu}\in\frac12\Z^g$, 
\(
\left[
\begin{smallmatrix}
 \vec{\mu} \\ \vec{\nu}
\end{smallmatrix}
\right]
\)
is called a {\it half-integer characteristic}. The set of half-integer
characteristics are divided into two, {\it odd\/} or {\it even},
according to whether $\Theta
\left[
\begin{smallmatrix}
 \vec{\mu} \\ \vec{\nu}
\end{smallmatrix}
\right]
(\vec{\xi}|\tau)
$
is an odd or even function. 

In the main body of the paper, we use the $A_{r-1}$-lattice
\[
  \left.\left\{ \vec{k} = (k_1,\dots, k_r)\in \Z^r \,\right|
      \textstyle \sum_\alpha k_\alpha = 0 \right\}.
\]
(See \secref{sec:root}.)
This is identified with $\Z^{r-1}$ by taking $(k_2,\dots, k_r)$. Then
we apply the above convention, i.e., the suffix runs
$\alpha=2,\dots,r$.

A theta function with a particular half-integer even characteristic
appears often in this paper:
\begin{equation}\label{eq:characteristic}
   \Theta_E(\vec{\xi}|\tau)
   = \Theta
\begin{bmatrix}
   \vec{0} \\ \vec{\nu}
\end{bmatrix}
(\vec{\xi}|\tau)
   ,
\qquad
\nu_2 = \frac12, \nu_3 = 0, \nu_4 = \frac12, \nu_5 = 0, \cdots 
.
\end{equation}
We denote this characteristic by $E$ and the corresponding theta
function by $\Theta_E$.

We also use the notation for the root system of Lie algebra of type
$A_{r-1}$. Then
\begin{equation*}
\begin{split}
   & \frac12 k_2 + \frac12 k_4 + \cdots \equiv
   -\frac12 k_2 - \frac22 k_3 - \cdots - \frac{r-1}2 k_r
\\
   =\; &\frac{r-1}4 k_1 + \frac{r-3}4 k_2 + \dots + \frac{1-r}4 k_r
   = \sum_{\alpha < \beta} \frac{k_\alpha - k_\beta}4
   = \frac12 \sum_i k^i,
\end{split}
\end{equation*}
where $\equiv$ means the equality modulo $\Z^{r-1}$. The last equality
is \eqref{eq:useful}. Therefore the characteristic is
\begin{equation*}
   {}^t\!
   \begin{pmatrix}
   \frac12 & \frac12 & \cdots
   \end{pmatrix}
\end{equation*}
in this notation.

We also use the theta function where the summation range is replaced
by
\[
  \left.\left\{ \vec{k} = (k_1,\dots, k_r)\in \Q^r \,\right|
    \textstyle \sum_\alpha k_\alpha = 0, 
    \forall \alpha\; l_\alpha \equiv -\frac{k}r \mod \Z
    \right\}
\]
for a fixed $k\in\Z$. It is
\(
    \Theta\begin{bmatrix}
      \vec{\alpha_k} \\ \vec{\Delta}
    \end{bmatrix}
\)
with
\(
   \vec{\alpha}_k = 
   \frac{k}r\; {}^t
   \begin{pmatrix}
     1 & 2 & \cdots & r-1
   \end{pmatrix}.
\)
We denote by $E_k$ this characteristic.

\subsection{}\label{subsec:g=1}

When $g=1$, we use the following notation for the theta functions:
\begin{equation*}
  \begin{split}
  & \theta_{00}(z,\tau) = \sum_{n\in\Z} q^{n^2} w^{2n},
\quad
  \theta_{01}(z,\tau) = \sum_{n\in\Z} (-1)^n q^{n^2} w^{2n},
\\
  & \theta_{10}(z,\tau) = \sum_{n\in\Z} q^{(n+\frac12)^2} w^{2n+1},
\quad
  \theta_{11}(z,\tau) = \sqrt{-1}\sum_{n\in\Z} (-1)^n
  q^{(n+\frac12)^2} w^{2n+1},
  \end{split}
\end{equation*}
where
\begin{equation*}
    q = \exp(\pi\sqrt{-1}\tau), \quad
    w = \exp(\pi\sqrt{-1}z).
\end{equation*}
This is the same as \cite{Mumford}. We have
\begin{equation*}
   \Theta_E = \theta_{01}, \qquad
   \Theta_{E_1} = \theta_{11}.
\end{equation*}

\subsection{Riemann surfaces and theta functions}

Let $C$ be a compact Riemann surface of genus $g$. Let $K_C$ be its
canonical bundle. We choose and fix a symplectic basis $A_1,\dots,
A_g, B_1, \dots, B_g$ of $H_1(C,\Z)$ so that $A_\alpha \cdot A_\beta =
0 = B_\alpha\cdot B_\beta$, $A_\alpha\cdot B_\beta =
\delta_{\alpha\beta}$ for $\alpha,\beta = 1,\dots, g$. We then have a
basis $\omega_1,\dots, \omega_g$ of holomorphic differentials
$H^0(C,K_C)$ such that $\int_{A_\alpha}\omega_\beta =
\delta_{\alpha\beta}$. The period matrix of $C$ is defined by
\begin{equation*}
   \tau_{\alpha\beta} = \int_{B_\beta} \omega_\alpha.
\end{equation*}
It is symmetric and its imaginary part is positive-definite.

Using the period matrix $(\tau_{\alpha\beta})$ of the Riemann surface
$C$, we consider the associated theta function
\(
   \Theta
   \left[
   \begin{smallmatrix}
     \vec{\mu} \\ \vec{\nu}
   \end{smallmatrix}
   \right]
   (\vec{\xi}|\tau)
\)
as in the previous section. We consider it as a multi-valued function
(or a section of a line bundle) on the Jacobian variety
\(
   J(C) = H^0(C,K_C)^*/H_1(C,\Z).
\)
Here $H^0(C,K_C)$ is identified with $\C^g$ by the basis
$\omega_1,\dots, \omega_g$.

We choose a base point $P_0$ in $C$. Then we have the {\it Abel-Jacobi map\/}
\begin{equation*}
  C\ni P \longmapsto \int_{P_0}^P \omega \in J(C); \qquad
  \omega\in H^0(C,K_C).
\end{equation*}
We denote it by $\mathbf A$.
It extends a map from $J_g(C)$ the Picard variety of divisor classes
(linear equivalence classes) of degree $g$ divisors. When $g=0$, it is
independent of the base point and we have an isomorphism $J_0(C) \cong
J(C)$.

Riemann's theorem (\cite[Theorem~1.1]{Fay}, \cite[Chap.~II,
3.1]{Mumford}) says that
there exists a vector $\vec{K}\in\C^g$ 
\begin{NB}
(denoted by $\vec{\Delta}$ in \cite{Mumford})
\end{NB}
such that for all $\vec{e}\in\C^g$ the composition
\(
   \Theta
   \circ(\mathbf A+\vec{e})
\)
either vanishes identically, or has $g$ zeroes $Q_1,\dots Q_g$ such
that $\sum_{k=1}^g \mathbf A(Q_k) + \vec{e} \equiv \vec{K} \mod
H_1(C,\Z)$. The vector is called the {\it Riemann constant}.
The vector $\vec{K}$ depends on the choice of the symplectic basis of
$H_1(C,\Z)$ and the base point $P_0$.
However, if we denote it by $\vec{K}_{P_0}$, then $\Delta = (g-1)P_0 +
\vec{K}_{P_0}\in J_{g-1}(C)$ is independent of $P_0$. See
\cite[VI.3.7]{FK}, \begin{NB}(denoted by $\Delta$)\end{NB}
\cite[Chap.~II, 3.11, 3.18]{Mumford}. \begin{NB}(denoted by $D_0$)\end{NB}
In fact, we have
\(
   \Theta = W_{g-1} - \Delta,
\)
where $\Theta$ is considered as a divisor in $J(C)$ as a zero set
$\Theta = 0$, and $W_{g-1} = \{ x_1 + \dots + x_{g-1} \mid x_\alpha\in
C\}$.

The set $\Sigma$ of divisor classes $D\in J_{g-1}(C)$ such that $2D =
K_C$ is called the set of {\it theta characteristics}. The above
$\Delta$ is an example. The set $\Sigma$ is bijective to
$\frac12 H_1(C,\Z)/H_1(C,\Z) \cong (\Z/2\Z)^{2g}$ (considered as a
subset in $J(C)$) via $D\mapsto D-\Delta$. We identify $\Sigma$ with a
characteristic for the theta function by the further identification
\(
  J(C) \cong \C^g/\Z^g\oplus \tau\Z^g
\)
given by the choice of cycles.

\subsection{Green functions}

Let $\phi\colon\widetilde{C} \to C$ be the universal covering of $C$.
We take a nonsingular odd theta characteristic $D \in J_{g-1}(C)$ (i.e.,
the theta divisor is smooth at $D$ \cite[IIIb,Lem. 1]{Mumford})
and let 
$\delta= \left[
   \begin{smallmatrix}
     \vec{\mu} \\ \vec{\nu}
   \end{smallmatrix}
   \right]$ be the corresponding half integer characteristic.
By Riemann's theorem,
\begin{equation}
\zeta(x)=\sum_{\alpha=1}^g 
\frac{\partial\Theta_{\delta}
}{\partial \xi_\alpha}(0)\omega^{\alpha}(x)
\end{equation}
is a section of $K_C$ which vanishes on $D$.
Since $D$ is a nonsingular odd characteristic, 
$H^0(C,K_C(-D)) \cong H^0(C,{\cal O}_C(D)) \cong{\Bbb C}$,
and hence 
$\sqrt{\zeta(x)}$ is a section of ${\cal O}_C(D)$.
We define a {\it prime form\/} by
\begin{equation}
E(x,y):=\frac{
\Theta_\delta
 \left(\int_x^y\vec{\omega}|\tau \right)}
{\sqrt{\zeta(x)}\sqrt{\zeta(y)}}.
\end{equation}
This is a holomorphic differential form on
$\widetilde{C} \times \widetilde{C}$.
It is also regarded as a holomorphic section of
a line bundle on $C \times C$:
Let $\pi_i\colon C \times C \to C$, $i=1,2$ be projections and
$\mu\colon C \times C \to J(C)$ be the map
$(x,y) \mapsto y-x$.
Then $E(x,y)$ is a section of
$\pi_1^*{\cal O}_C(-\Delta) \otimes \pi_2^* {\cal O}_C(-\Delta) \otimes
\mu^*({\cal O}_{J(C)}(\Theta))$, where 
$\Theta$ is the theta divisor.
\begin{NB}
Let $s(x)$ be the section of ${\cal O}_C(P_0)$. 
Then $s(x)\sqrt{\zeta(x)}$ defines a section of 
${\cal O}_C(D+P_0)$.
We note that 
${\cal O}_C(D+P_0)$ is the pull-back of $T^*_\delta({\cal O}_{J(C)}(\Theta))$ 
by $C \hookrightarrow J(C)$.
Then we see that $E'(x,y):=E(x,y)/(s(x)s(y))$ is a meromorphic section of
$\pi_1^*{\cal O}_C(-\Delta-P_0) \otimes \pi_2^* {\cal O}_C(-\Delta-P_0) \otimes
\mu^*({\cal O}_{J(C)}(\Theta))$. 
For example, 
\begin{equation}
\begin{split}
&E'(x,y+\tau \vec{m}+\vec{n})/E'(x,y)\\
=&\exp \left(-\pi \sqrt{-1} {}^t\vec{m} \tau \vec{m}-
2\pi \sqrt{-1} {}^t\vec{m}\int_x^y \vec{\omega} \right)/ 
\exp \left(-\pi \sqrt{-1} {}^t\vec{m} \tau \vec{m}-
2\pi \sqrt{-1} {}^t\vec{m}\int_{P_0}^y \vec{\omega} \right). 
\end{split}
\end{equation}
Similar relation holds for $E'(x+\tau \vec{m}+\vec{n},y)/E'(x,y)$,
since $\delta$ is a half characteristic.
\end{NB}

We pick up some properties of $E(x,y)$.
\begin{enumerate}
\item
$E(x,y)=0 \Longleftrightarrow \phi(x)=\phi(y)$.
\item
$E(x,y)$ has a first order zero along the diagonal
$\Delta_C \subset C \times C$ and locally
$E(x,y)=\frac{x-y}{\sqrt{dx}\sqrt{dy}}(1+O((x-y)^2))$.
\item
$E(x,y)=-E(y,x)$.
\begin{NB}
I removed (4).
If $x$ is moved by an $A$-cycle, then
$E(x,y)$ is invariant.
If $x$ is moved by $B_\alpha$, then
$E(x,y)$ is multiplied by $\exp \left(-\pi \sqrt{-1} \tau_{\alpha \alpha}+
2\pi\sqrt{-1}\int_x^y \omega_\alpha \right)$ as a weight
$(-1/2,-1/2)$ differential form
(i.e., as a meromorphic section of $\pi_1^* {\cal O}_C(-\Delta) \otimes 
\pi_2^* {\cal O}_C(-\Delta)$).

\end{NB}
\end{enumerate}
\begin{NB}
These properties characterize $E(x,y)$.
(It is determined by the location of zero+(4) up to constant,
and (2) determines the scalar).  
\end{NB}

Let
\begin{equation}
W(z_1,z_2)=\partial_{z_1}\partial_{z_2} \log E(z_1,z_2).
\end{equation}
This is a well-defined meromorphic 2-form on $C \times C$ and it is used to
construct differentials of the 2nd kind.
For $c \in {\Bbb C}^g$ with $\Theta(c) \ne 0$,
we set
\begin{equation}
\Psi_c(z_1,z_2)=\frac{\Theta_c(\int_{z_2}^{z_1}\vec{\omega}|\tau)}
{\Theta_c(0)E(z_1,z_2)}.
\end{equation}
It is called the {\it Szeg\"{o} kernel}.

By Fay's trisecant identity (\cite[p. 34, formula 45]{Fay} or
\cite[IIIb,2]{Mumford}),
we have
\begin{equation}\label{eq:Fay}
   \Psi_c^2(z_1,z_2) = W(z_1,z_2)
   + \sum_{\alpha,\beta}
     \omega^\alpha(z_1)\omega^\beta(z_2)
     \frac{\partial^2}{\partial \xi_\alpha\partial \xi_\beta}
     \log\Theta_c(0|\tau).
\end{equation}
for a half integer characteristic $c$
(\cite[Cor. 2.12 formula 38]{Fay}, \cite[IIIb,3 (2)]{Mumford}).

\begin{NB}
Use
$\Theta(c+\vec{z})\Theta(c-\vec{z})=
\exp \left(-2\pi \sqrt{-1} {}^t\!\vec{\mu}\tau\vec{\mu}
- 4\pi\sqrt{-1}\,{}^t\!\vec{\mu} \vec{\nu}\right)
\Theta_c(\vec{z})\Theta_{-c}(\vec{z})
=\frac{\Theta(c)^2}{\Theta_c(0)^2}\Theta_c(\vec{z})\Theta_{-c}(\vec{z})$.

\end{NB}
\begin{NB}
In hep-th/9802007, (A2) seems to be wrong.
$\prod_{i<j}E(\xi_i,\xi_j)E(\zeta_i,\zeta_j)$ seems to be
$\prod_{i<j}E(\xi_i,\xi_j)E(\zeta_j,\zeta_i)$ 
(cf Kawamoto et al, Prop. 6.14).
So (A.3) should be modified (there is also typos there).
the correct form is  
\begin{equation}
\begin{split}
&\theta_e(\vec{0})\theta_e(\vec{\xi_1}+\vec{\xi_2}
-\vec{\zeta_1}-\vec{\zeta_2})E(\xi_1,\xi_2)E(\zeta_1,\zeta_2)\\
=&-\theta_e(\vec{\xi_1}-\vec{\zeta_1})\theta_e(\vec{\xi_2}
-\vec{\zeta_2})E(\xi_2,\zeta_1)E(\xi_1,\zeta_2)+
\theta_e(\vec{\xi_2}-\vec{\zeta_1})\theta_e(\vec{\xi_1}
-\vec{\zeta_2})E(\xi_1,\zeta_1)E(\xi_2,\zeta_2)
\end{split}
\end{equation}
This coincides with \cite{Mumford} IIIb sec. 2.
So the sign of the second line of (A.5) is $-$.
There is another mistake of sign in the computation of limit
$\zeta_1 \to \zeta_2$.
Therefore the final result is correct.
\end{NB}

\subsection{Hyperelliptic curves}

Let $Q(z)$ be a polynomial of degree $2g+2$. Let $C = \{ y^2 = Q(z)\}$
be the corresponding hyperelliptic curve of genus $g$. We denote by
$\iota\colon C\to C$ be the involution. Let $\{ Q_1,\dots, Q_{2g+2}\}$
be the set of branched points, i.e., the roots of $Q(z) = 0$.
We choose cycles $A_\alpha$, $B_\alpha$ as in Figure~\ref{fig:SW}
where we replace as
\begin{equation*}
 z_1^+\rightarrow Q_1, \quad
 z_1^-\rightarrow Q_2, \quad
 z_2^-\rightarrow Q_3, \quad
 z_2^+\rightarrow Q_4, \quad \cdots.
\end{equation*}
Then the choice of cycles is exactly the same as \cite[p.12
Example]{Fay} with shifting the numbering by $1$, i.e., $A_2 \rightarrow
A_1$, $B_2\rightarrow B_1$, etc. If we choose $Q_1$ for the base point
of the Abel-Jacobi map, we have (\cite[p.14]{Fay},
\cite[Chap.~IIIa.5]{Mumford})
\begin{equation*}
   \vec{K} =
   \tau\, {}^t\!
   \begin{pmatrix}
     \frac12 & \frac12 & \cdots &\frac12
   \end{pmatrix}
   +
   {}^t\!
   \begin{pmatrix}
     \frac{1}2 & \frac{2}2 & \cdots &\frac{g}2
   \end{pmatrix}
   .
\end{equation*}

Let $L$ be the divisor class of degree $2$ containing $P+\iota P$ for
$P\in C$. Then the set of theta characteristics is bijective to the
set of subsets $T\subset \{ Q_1,\dots, Q_{2g+2}\}$ with $\# T \equiv
(g+1)\mod 2$ modulo the equivalence relation $T\sim T^c$ by
\begin{equation*}
   T \longmapsto \sum_{P\in T} P + \frac{g-1-\# T}2 L.
\end{equation*}
Under this correspondence, the vector $\Delta = (g-1)Q_1 + \vec{K}$
is mapped to the $\{ Q_1, Q_3, \dots, Q_{2g+1}\}$.

When $\# T = (g+1)$, the corresponding Szeg\"o kernel is given
explicitly by
\begin{equation}\label{eq:Szego}
   \Psi_c(z_1,z_2)
   = \frac12\left(
     \sqrt[4]{\frac{\psi(z_1)}{\psi(z_2)}}
     +
     \sqrt[4]{\frac{\psi(z_2)}{\psi(z_1)}}
     \right)
     \frac{\sqrt{dz_1 dz_2}}{z_1-z_2},
\end{equation}
where
\(
   \psi(z) = \prod_{Q_\alpha\in T} (z - Q_\alpha)
   \times \prod_{Q_\beta\in T^c} (z - Q_\beta)^{-1}.
\)
See \cite[p.12 Example]{Fay}.

\begin{NB}
For the Seiberg-Witten curve, we have $\psi(z) = (P(z) -
2\Lambda^r)/(P(z) + 2\Lambda^r)$ for $c = E$. Then 
\begin{equation*}
\Psi_E^2(z_1,z_2)=\frac{y(z_1) y(z_2)+P(z(z_1))P(z(z_2))-4\Lambda^{2r}}
{2y(z_1)y(z_2)}
\frac{dz(z_1)dz(z_2)}{(z(z_1)-z(z_2))^2}.
\end{equation*}
\end{NB}

\section{Equivariant Borel-Moore homology}\label{sec:BMhom}

We use equivariant Borel-Moore homology in this paper. For the usual
Borel-Moore homology, see e.g., \cite[\S B.2]{Fulton}. As we only use
the Borel-Moore homology, we denote it by $H_*(\ )$. If we do not
specify the coefficients, we mean the complex coefficients.

The following properties are crucial.
\begin{aenume}
\item If $X$ is nonsingular, $H_i(X)$ is isomorphic to the ordinary
cohomology group $H^{2\dim X-i}(X)$.

\item For an irreducible algebraic variety $X$, its fundamental class
$[X]\in H_{2\dim X}(X)$ is defined.
\item For a proper continuous map $f\colon X\to Y$, the push-forward
homomorphism $f_*\colon H_*(X)\to H_*(Y)$ is defined.
\item If $U\subset X$ is open with complement $Y = X\setminus U$, we
have the long exact sequence
\begin{equation*}
  \cdots \to H_i(Y)\xrightarrow{\iota_*} H_i(X) \xrightarrow{j^*}
  H_i(U) \to H_{i-1}(Y) \to \cdots,
\end{equation*}
where $\iota\colon Y\to X$, $j\colon U\to X$ are inclusions, and $j^*$
is the restriction homomorphism.
\end{aenume}

For an equivariant Borel-Moore homology, we use the one given in
\cite{Lu:GH}, but we shift the degree so that the fundamental class
$[X]$ has degree $2\dim X$. This definition is the same as \cite{EG}.

Let us recall the definition briefly.
Let $G$ be a linear algebraic group acting on an algebraic variety $X$.
(Everything is over $\C$.) We have a finite dimensional approximation
of the classifying space $EG\to BG$, i.e., for any $n$, there exists a
smooth irreducible variety $U$ with $G$-action such that
\begin{aenume}
\item The quotient $U \to U/G$ exists and is a principal $G$-bundle.
\item $H^i(U) = 0$ for $i=1,\dots,n$.
\end{aenume}
We then define
\begin{equation*}
   H^G_n(X) = H_{n-2\dim G+2\dim U}(X\times_G U).
\end{equation*}
Here $U$ is smooth, in particular $\dim U$ makes sense. One can show
that this is independent of the choice of $U$, using the double
fibration argument.

Note that $H^G_n(X) = 0$ if $n > 2\dim X$, but $H^{G}_n(X)$ may be
nonzero for $n<0$. ($X$ is pure dimensional.)

On the other hand, we define the equivariant {\it co}-homology as
\begin{equation*}
   H_G^n(X) = H^n(X\times_G U),
\end{equation*}
where $H^n(\ )$ is the ordinary cohomology. This coincides with the
usual definition. It is a graded ring.
We have the Poincar\'e duality isomorphism
\begin{equation*}
   H_G^n(X) \cong H^G_{2\dim X - n}(X)
\end{equation*}
when $X$ is nonsingular.

As a projection $X\times_G U\to U/G$ is flat, $H^G_*(X)$ has a
structure of a $H_G^*(\mathrm{pt})$-module.

Suppose that $G$ is reductive. Then $H_G^*(\mathrm{pt})$ is isomorphic
to $S^*(\mathfrak h^*)^W$, where $\mathfrak h$ is a Cartan subalgebra,
$S^*(\mathfrak h^*)$ is the symmetric algebra of its dual, and $W$ is
the Weyl group. We denote this by $S$ or $S(G)$.

Let $T$ be a torus acting on $X$. Let $X^T$ be the fixed point set and
$\iota\colon X^T\to X$ be the inclusion. We have the push-forward
homomorphism
\(
   \iota_*\colon H_*^T(X^T)\to H_*^T(X).
\)
Since $T$ acts trivially on $X^T$, we have
\(
   H_*^T(X^T) = H_*(X^T)\otimes_\C S
\)
The localization theorem (see \cite{AB}) says that $\iota_*$ becomes an
isomorphism after tensoring the quotient field $\mathcal S$ of $S$.

When $X$ is nonsingular, the inverse of $\iota_*$ can be explicitly
given. Let $X^T = \bigsqcup F_i$ be the decomposition to irreducible
components. Each $F_i$ is nonsingular. Let $N_i$ be the normal bundle
of $F_i$ in $X$. Then we have
\begin{equation*}
   \left(\iota_*\right)^{-1} = \sum_i \frac1{e_T(N_i)} \iota_i^*,
\end{equation*}
where $e_T(N_i)$ is the equivariant Euler class and $\iota_i^*$ is the
pull-back homomorphism for the inclusion $\iota_i\colon F_i\to X$
defined via the Poincar\'e duality homomorphism.

\section{The proof of \eqref{eq:Ochiai} by Hiroyuki Ochiai}
\label{sec:Ochiai}

Let
\begin{equation*}
   (a)_k = (a,q)_k = (1-a)(1-aq)\cdots (1-aq^{k-1}),
\qquad
   (a)_\infty = \prod_{d=0}^\infty (1 - aq^d).
\end{equation*}


We start with the formula, \cite[p.16,(7.2)]{Ramanujan}.

We substitute $d = aq/c$, then we have
$$
\sum_{k=0}^\infty
\frac{(a)_k (e)_k (f)_k (1-a q^{2k})(-a/qef)^k q^{k(k+3)/2}}%
     {(aq/e)_k (aq/f)_k (q)_k (1-a)}
= \frac{(aq)_\infty (aq/ef)_\infty}{(aq/e)_\infty (aq/f)_\infty}.
$$
($q$-hypergeometric part vanishes since $(1)_k = 0$ for $k\ge1$.)
We put $e=-(aq)^{1/2}$, 
Then we have 
$$
\sum_{k=0}^\infty
\frac{(a)_k (f)_k (1-a q^{2k})(a^{1/2} q^{-3/2}/ f)^k q^{k(k+3)/2}}%
     {(aq/f)_k (q)_k (1-a)}
=\frac{(aq)_\infty (-(aq)^{1/2}/f)_\infty}{(-(aq)^{1/2})_\infty (aq/f)_\infty}.
$$
Finally, we put $f=a q^{-1/2}$.
Then we have
$$
\sum_{k=0}^\infty
\frac{(a)_k (a q^{-1/2})_k (1-a q^{2k})(a^{-1/2}q^{-1})^k q^{k(k+3)/2}}%
     {(q^{3/2})_k (q)_k}
=\frac{(a)_\infty (-a^{-1/2}q)_\infty}{(-(aq)^{1/2})_\infty (q^{3/2})_\infty}.
$$
Using Jacobi triple product identity 
$(q)_\infty (-a^{1/2})_\infty (-a^{-1/2}q)_\infty
= \sum_{l = -\infty}^\infty (q/a)^{l/2} q^{l^2/2}$,
the right hand side is
$$
\frac{(a)_\infty}{(-a^{1/2}q^{1/2})_\infty (-a^{1/2})_\infty (q)_\infty
(q^{3/2})_\infty} \sum_{l = -\infty}^\infty (q/a)^{l/2} q^{l^2/2}.
$$
Using $(b)_\infty (b q^{1/2})_\infty = (b, q^{1/2})_\infty$,
and $(b^2)_\infty = (b, q^{1/2})_\infty (-b,q^{1/2})_\infty$,
we get
$$
\sum_{k=0}^\infty
\frac{(a q^{-1/2}, q^{1/2})_{2k} (1-a q^{2k})(a^{-1/2}q^{-1})^k q^{k(k+3)/2}}%
     {(q, q^{1/2})_{2k}}
= \frac{(a^{1/2}, q^{1/2})_\infty}{(q,q^{1/2})_\infty} 
\sum_{l = -\infty}^\infty (q/a)^{l/2} q^{l^2/2}.
$$
Also
$$
\sum_{k=0}^\infty
\frac{(a q^{-1/2}, q^{1/2})_{2k} (1-a q^{2k})(a^{-1/2}q^{-1})^k q^{k(k+3)/2}}%
     {(q^{1/2}, q^{1/2})_{2k+1}}
= \frac{(a^{1/2}, q^{1/2})_\infty}{(q^{1/2},q^{1/2})_\infty} 
\sum_{l = -\infty}^\infty (q/a)^{l/2} q^{l^2/2}.
$$
Now we substitute $a \mapsto q^2 t^2$ and  $q \mapsto q^4 t^2$.
Then
$$
\sum_{k=0}^\infty
\frac{(t, q^2 t)_{2k} (1- q^{8k+2} t^{4k+2}) q^{k(2k+1)} t^{k^2} }%
     {(q^2 t, q^2 t)_{2k+1}}
= \frac{(qt, q^2 t)_\infty}{(q^2 t,q^2 t)_\infty} 
\sum_{l = -\infty}^\infty q^{l(2l+1)} t^{l^2}.
$$
Using the identity
$$
1- q^{8k+2} t^{4k+2} =
(1-q^{4k} t^{2k+1}) q^{4k+2} t^{2k+1} + (1-q^{4k+2} t^{2k+1}),
$$
we see the left hand side is
$$
\sum_{k=0}^\infty
\left(
\frac{(t, q^2 t)_{2k}}{(q^2 t, q^2 t)_{2k}}  q^{k(2k+1)} t^{k^2} +
\frac{(t, q^2 t)_{2k+1}}{(q^2 t, q^2 t)_{2k+1}}  q^{(k+1)(2k+3)-1} t^{(k+1)^2}
\right).
$$
This is the end of the proof.

\section{Perturbation term}\label{sec:perturb}
\subsection{One parameter version}
Let
\begin{equation*}
   \gamma_{\hbar}(x;\Lambda)
   = \left.\frac{d}{ds}\right|_{s=0}
   \frac{\Lambda^s}{\Gamma(s)}\int_0^\infty
    \frac{dt}t t^s \frac{e^{-tx}}{(e^{\hbar t}-1)(e^{-\hbar t}-1)},
\end{equation*}
where $\Gamma(s)$ is the Gamma function
\begin{equation*}
   \Gamma(s) = \int_0^\infty \frac{dt}t t^s e^{-t}.
\end{equation*}
The integral in the right hand side converges when $\Re(s) > 2$. The
analytic continuation can be done by the standard procedure using the
Taylor expansion of the integrand. (See below.)

If we formally expand as
\begin{equation*}
  \frac1{(e^{\hbar t}-1)(e^{-\hbar t}-1)}
  = \sum_{m,n\ge 0} e^{\hbar( m - n) t},
\end{equation*}
we get
\begin{equation*}
    \gamma_{\hbar}(x;\Lambda)
    \overset{\text{formally}}{=}
    \sum_{m,n\ge 0} \log\left(\frac{x-\hbar(m-n)}{\Lambda}\right).
\end{equation*}
Thus $\gamma_{\hbar}(x;\Lambda)$ is a regularization of the right hand
side.

We introduce Bernoulli numbers by
\begin{equation*}
  \frac{t}{e^t-1} = \sum_{n=0}^\infty \frac{B_n}{n!} t^n.
\end{equation*}
We have $B_0 = 1$, $B_1 = -\frac12$, $B_2 = \frac16$, $B_{2k+1} = 0$
for $k\ge 1$. Note
\begin{equation*}
    \frac1{(e^t-1)(e^{-t}-1)}
    = \frac{d}{dt} \frac{1}{e^t-1}
    = -\frac1{t^2} + \sum_{g=1}^\infty \frac{B_{2g}}{2g(2g-2)!} t^{2g-2}.
\end{equation*}
Then
\begin{equation}\label{eq:expand}
  \begin{split}
   & \gamma_\hbar(x;\Lambda)
\\
   =\; & \left.\frac{d}{ds}\right|_{s=0}
   \Biggl[
   \begin{aligned}[t]
     &
     - (\frac{x}{\hbar})^2(\frac{\Lambda}{x})^s \frac{\Gamma(s-2)}{\Gamma(s)}
     + \frac{B_2}2 (\frac{\Lambda}{x})^s
     \\
     & \qquad\qquad
     + \sum_{g=2}^\infty \frac{B_{2g}}{2g(2g-2)!} (\frac{\hbar}x)^{2g-2}
      \frac{\Gamma(s+2g-2)}{\Gamma(s)}
   \Biggr]
   \end{aligned}
\\
  = \; &
  \hbar^{-2}
  \left\{\frac12 x^2 \log\left(\frac{x}{\Lambda}\right) - \frac34 x^2\right\}
  - \frac1{12}\log\left(\frac{x}{\Lambda}\right)
  + \sum_{g=2}^\infty \frac{B_{2g}}{2g(2g-2)} (\frac{\hbar}x)^{2g-2}.
  \end{split}
\end{equation}
We have the difference equation
\begin{equation*}
  \gamma_\hbar(x+\hbar;\Lambda)
  + \gamma_\hbar(x-\hbar;\Lambda)
  - 2\gamma_\hbar(x;\Lambda) = \log\left(\frac{x}{\Lambda}\right).
\end{equation*}
In fact, the left hand side is equal to
\begin{equation*}
  \left.\frac{d}{ds}\right|_{s=0}
   \frac{\Lambda^s}{\Gamma(s)}\int_0^\infty
    \frac{dt}t t^s e^{-tx}
    \frac{e^{t\hbar} + e^{-t\hbar} - 2}{(e^{\hbar t}-1)(e^{-\hbar t}-1)}
    = - \left.\frac{d}{ds}\right|_{s=0}
    \left(\frac{\Lambda}{x}\right)^s.
\end{equation*}

We have
\begin{equation*}
\begin{split}
   & \gamma_\hbar\left(x + \frac\hbar2;\Lambda\right)
   - \gamma_\hbar\left(x - \frac\hbar2;\Lambda\right)
   = \left.\frac{d}{ds}\right|_{s=0}
   \frac{\Lambda^s}{\Gamma(s)}\int_0^\infty
    \frac{dt}t t^s e^{-tx}
    \frac{e^{-\frac{t\hbar}2} - e^{\frac{t\hbar}2}}
    {(e^{\hbar t}-1)(e^{-\hbar t}-1)}
\\
   =\; &
   \left.\frac{d}{ds}\right|_{s=0}
   \frac{\Lambda^s}{\Gamma(s)}\int_0^\infty
    \frac{dt}t t^s \frac{e^{-t(x+\frac\hbar2)}}{1-e^{-\hbar t}}
   =
   \left.\frac{d}{ds}\right|_{s=0}
   \left(\frac{\Lambda}{\hbar}\right)^s
   \zeta\left(s,\frac{x}{\hbar} + \frac12\right)
\\
   =\; &
   \log\left(\frac\Lambda{\hbar}\right)
   \zeta\left(0,\frac{x}{\hbar} + \frac12\right)
   + \zeta'\left(0,\frac{x}{\hbar} + \frac12\right)
   =
   \log\left(\frac1{\sqrt{2\pi}}
     \left(\frac{\hbar}\Lambda\right)^{\frac{x}{\hbar}}
     \Gamma\left(\frac{x}{\hbar} + \frac12 \right)
   \right)
\end{split}
\end{equation*}
where $\zeta(s,a)$ is the Hurwitz zeta function. And at the final
equality, we have used the Lerch formula (see \cite[XV\S]{Lang}).

We have
\begin{equation*}
   \gamma_{\hbar}(x;\Lambda) + \gamma_{\hbar}(-x;\Lambda)
   = 2\gamma_{\sqrt{-1}\hbar}(\sqrt{-1}x;\Lambda).
\end{equation*}
This can be seen from the expansion \eqref{eq:expand}.

\subsection{Two parameter version}
Let us introduce a generalization of $\gamma_\hbar(x;\Lambda)$:
\begin{equation*}
   \gamma_{\ve_1,\ve_2}(x;\Lambda)
   = \left.\frac{d}{ds}\right|_{s=0}
   \frac{\Lambda^s}{\Gamma(s)}\int_0^\infty
    \frac{dt}t t^s \frac{e^{-tx}}{(e^{\ve_1 t}-1)(e^{\ve_2 t}-1)}.
\end{equation*}
This is formally equal to
\begin{equation*}
    \sum_{m,n\ge 0} \log\left(\frac{x-m\ve_1-n\ve_2}{\Lambda}\right).
\end{equation*}

The difference equation is
\begin{equation*}
  \gamma_{\ve_1,\ve_2}(x-\ve_1;\Lambda)
  + \gamma_{\ve_1,\ve_2}(x-\ve_2;\Lambda)
  - \gamma_{\ve_1,\ve_2}(x;\Lambda)
  - \gamma_{\ve_1,\ve_2}(x-\ve_1-\ve_2;\Lambda)
  = \log\left(\frac{x}{\Lambda}\right).
\end{equation*}

Let $k$ be an integer. We have
\begin{equation*}
  \begin{split}
   & \gamma_{\ve_1,\ve_2-\ve_1}(x+\ve_1k;\Lambda)
    + \gamma_{\ve_1-\ve_2,\ve_2}(x+\ve_2k;\Lambda)
\\
  =\; &
  \left.\frac{d}{ds}\right|_{s=0}
   \frac{\Lambda^s}{\Gamma(s)}\int_0^\infty
    \frac{dt}t t^s e^{-tx} \left\{
      \frac{e^{-t\ve_1 k}}{(e^{\ve_1 t}-1)(e^{(\ve_2-\ve_1)t}-1)}
      + \frac{e^{-t\ve_2 k}}{(e^{(\ve_1-\ve_2)t}-1)(e^{\ve_2t}-1)}
      \right\}
\\
  =\; &
  \left.\frac{d}{ds}\right|_{s=0}
   \frac{\Lambda^s}{\Gamma(s)}\int_0^\infty
    \frac{dt}t t^s e^{-tx}
    \frac{(1-e^{-\ve_2t})e^{-k\ve_1t} - (1-e^{-\ve_1t})e^{-k\ve_2t}}
   {(e^{\ve_1 t}-1)(e^{\ve_2 t}-1)(e^{-\ve_1t}-e^{-\ve_2t})}.
  \end{split}
\end{equation*}
We claim that this is equal to
\begin{equation}\label{eq:pert_shift}
   \gamma_{\ve_1,\ve_2}(x;\Lambda) + 
   \log s^{-k}(\ve_1,\ve_2,x) - \frac{k(k-1)}2 \log\Lambda,
\end{equation}
where $s^{-k}(\ve_1,\ve_2,x)$ is given by \eqref{eq:s}.
If $k=0$ or $1$, this is obvious. Suppose that $k \ge 2$. Then the above is
equal to
\begin{equation*}
  \sum_{\substack{l,m\ge 0\\ l+m=k-1}}
    \gamma_{\ve_1,\ve_2}(x+l\ve_1+m\ve_2;\Lambda)
  - 
  \sum_{\substack{l,m\ge 1\\ l+m=k}}
     \gamma_{\ve_1,\ve_2}(x+l\ve_1+m\ve_2;\Lambda).
\end{equation*}
On the other hand, \eqref{eq:pert_shift} is equal to
\begin{multline*}
   \gamma_{\ve_1,\ve_2}(x;\Lambda) + 
   \sum_{\substack{l,m\ge 0\\ l+m\le k-2}} \log\left( 
     \frac{x + (l+1)\ve_1 + (m+1)\ve_2}\Lambda\right).
\\
 = \left(
  \sum_{l= 0, m = 0}
  + 
  \sum_{\substack{l\ge 0, m\ge 1 \\ l+m\le k-1}}
  + 
  \sum_{\substack{l\ge 1, m\ge 0 \\ l+m\le k-1}}
  - 
  \sum_{\substack{l\ge 1, m\ge 1 \\ l+m\le k}}
  - 
  \sum_{\substack{l\ge 0, m\ge 0 \\ l+m\le k-2}}
  \right)
  \gamma_{\ve_1,\ve_2}(x + l\ve_1 + m\ve_2;\Lambda)
\\
 = \left(
  \sum_{l= 0, m = 0}
  + 
  \sum_{\substack{l\ge 0, m\ge 0 \\ l+m= k-1}}
  - \sum_{\substack{l\ge 0, m = 0 \\ l+m\le k-1}}
  - \sum_{\substack{l\ge 1, m\ge 1 \\ l+m = k}}
  + \sum_{\substack{l\ge 1, m= 0 \\ l+m\le k-1}}
  \right)
  \gamma_{\ve_1,\ve_2}(x + l\ve_1 + m\ve_2;\Lambda)
\end{multline*}
by the difference equation. Thus we get the assertion when $k \ge 2$.

Similarly we have
\begin{multline*}
  \gamma_{\ve_1,\ve_2-\ve_1}(x+\ve_1k;\Lambda)
    + \gamma_{\ve_1-\ve_2,\ve_2}(x+\ve_2k;\Lambda)
\\
  =  \gamma_{\ve_1,\ve_2}(x;\Lambda) + 
   \sum_{\substack{l,m\ge 0\\ l+m\le -k-1}} \log\left( 
     \frac{x - l\ve_1 - m\ve_2}\Lambda\right)
\end{multline*}
for $k\le -1$. This is nothing but the assertion.

\subsection{Expansion}

Let us define $c_n$ ($n=0,1,2,\dots$) by
\begin{equation*}
   \frac1{(e^{\ve_1 t}-1)(e^{\ve_2 t}-1)}
   = \sum_{n=0}^\infty \frac{c_n}{n!} t^{2-n}.
\end{equation*}
We have
\begin{equation*}
   c_0 = \frac1{\ve_1\ve_2}, \quad
   c_1 = -\frac{\ve_1+\ve_2}{2\ve_1\ve_2}, \quad
   c_2 = \frac{\ve_1^2 + \ve_2^2 + 3\ve_1\ve_2}{6\ve_1\ve_2}, \quad
   \cdots.
\end{equation*}
Then
\begin{equation}\label{eq:pert_expand}
\begin{split}
   &\gamma_{\ve_1,\ve_2}(x;\Lambda)
   = \left.\frac{d}{ds}\right|_{s=0}
   \frac{\Lambda^s}{\Gamma(s)} \sum_{n=0}^\infty c_n
   \int_0^\infty \frac{dt}t t^{n+s-2} e^{-tx}
\\
   =\; & \left.\frac{d}{ds}\right|_{s=0}
   \left(\frac{\Lambda}{x}\right)^s \sum_{n=0}^\infty c_n
   x^{2-n} \frac{\Gamma(n+s-2)}{\Gamma(s)}
\\
   =\; & 
   \begin{aligned}[t]
     & \frac1{\ve_1\ve_2}\left\{
       - \frac12 x^2 \log\left(\frac{x}{\Lambda}\right) + \frac34 x^2\right\}
     + \frac{\ve_1+\ve_2}{2\ve_1\ve_2} \left\{
       - x \log\left(\frac{x}{\Lambda}\right) + x\right\}
     \\
     & \qquad
     - \frac{\ve_1^2 + \ve_2^2 + 3\ve_1\ve_2}{12\ve_1\ve_2}
     \log\left(\frac{x}{\Lambda}\right)
  + \sum_{n=3}^\infty \frac{c_n x^{2-n}}{n(n-1)(n-2)}.
   \end{aligned}
\end{split}
\end{equation}
In particular, we have
\begin{align}
   & \gamma_{\ve_1,\ve_2}(x;\Lambda e^u)
   = \gamma_{\ve_1,\ve_2}(x;\Lambda)
   + u\left\{
     \frac{x^2}{2\ve_1\ve_2} + \frac{x(\ve_1+\ve_2)}{2\ve_1\ve_2}
     + \frac{\ve_1^2 + \ve_2^2 + 3\ve_1\ve_2}{12\ve_1\ve_2}
   \right\}, \label{eq:e^u}
\\
   &
\begin{aligned}[c]
   & \gamma_{\ve_1,\ve_2}(x;\Lambda)
     + \gamma_{\ve_1,\ve_2}(-x;\Lambda)
\\
   & \quad = 
     \frac2{\ve_1\ve_2}\left\{
       - \frac12 x^2 \log\left(\frac{\sqrt{-1}x}{\Lambda}\right)
       + \frac34 x^2\right\}
     + \frac{\ve_1+\ve_2}{2\ve_1\ve_2} \pi\sqrt{-1}x
\\
     & \qquad
     - \frac{\ve_1^2 + \ve_2^2 + 3\ve_1\ve_2}{6\ve_1\ve_2}
     \log\left(\frac{\sqrt{-1}x}{\Lambda}\right)
  + \sum_{g=2}^\infty \frac{2 c_{2g} x^{2-2g}}{2g(2g-1)(2g-2)}.
   \end{aligned} \label{eq:double}
\end{align}

\end{document}